\documentclass[12pt,reqno]{amsart}
\usepackage{a4wide}
\usepackage{amsmath,amsfonts,amsthm,cite,amssymb}
\newcommand{\middle}{{}}
\usepackage[dvips]{graphicx}
\usepackage{subfigure}

\usepackage{times,txfonts}
\newfont{\cyr}{wncyr10}
\usepackage{titlesec}
\titleformat{\section}
{\bfseries\scshape\centering}
{\thesection.}{.5em}{}
\titleformat{\subsection}
{\rmfamily\bfseries}
{\thesubsection}{.5em}{}
\titleformat{\subsubsection}
{\rmfamily\bfseries}
{\thesubsubsection}{.5em}{}
\numberwithin{equation}{section}
\renewcommand{\thefootnote}{\fnsymbol{footnote}}

\DeclareMathOperator{\CS}{CS}

\renewcommand{\Im}{\operatorname{Im}}
\DeclareMathOperator{\Li}{Li_2}

\DeclareMathOperator{\Vol}{Vol}
\DeclareMathOperator{\length}{Length}
\DeclareMathOperator{\torsion}{Torsion}
\DeclareMathOperator{\Isom}{Isom}
\newcommand{\I}{\mathrm{i}}
\newcommand{\E}{\mathrm{e}}

\newcommand{\Lob}{\text{\cyr L}}
\begin{document}

\baselineskip 16pt
\parskip 8pt



\title[Generalized Volume Conjecture \& the A-polynomials]{
  Generalized Volume Conjecture
  and the A-Polynomials:
  \\
  The Neumann--Zagier Potential Function
  as a Classical Limit of Quantum Invariant
}

\author[K. Hikami]{Kazuhiro \textsc{Hikami}}
\email{\texttt{hikami@phys.s.u-tokyo.ac.jp}}
\urladdr{http://gogh.phys.s.u-tokyo.ac.jp/{\textasciitilde}hikami/}
\date{April 4, 2006}

\address{Department of Physics, Graduate School of Science,\\
  University of Tokyo,\\
  Hongo 7--3--1, Bunkyo, Tokyo 113--0033, Japan.
}



\begin{abstract}
  We study quantum invariant
  $Z_\gamma(\mathcal{M})$ for cusped
  hyperbolic 3-manifold $\mathcal{M}$.
  We construct this  invariant based on oriented
  ideal
  triangulation of $\mathcal{M}$ by assigning   to each tetrahedron
  the quantum dilogarithm
  function, which is introduced by Faddeev in studies of the modular
  double of the quantum group.
  Following Thurston and Neumann--Zagier, we deform a complete hyperbolic
  structure of $\mathcal{M}$,
  and correspondingly we  define  quantum
  invariant $Z_\gamma(\mathcal{M}_u)$.
  This quantum invariant is shown to give
  the Neumann--Zagier
  potential function
  in the classical limit
  $\gamma\to 0$,
  and the A-polynomial can be derived from the potential function.
  We explain our construction by
  taking  examples of 3-manifolds such as complements of
  hyperbolic knots and punctured torus bundle over the circle.
\end{abstract}



\maketitle
\renewcommand{\thefootnote}{\arabic{footnote}}

\section{Introduction}

Since the quantum invariant of knots/links and 3-manifolds
as a generalization of the Jones polynomial~\cite{Jones85} is
constructed by Witten~\cite{EWitt89a} by use of the Chern--Simons path
integral, studies on quantum invariants have been much developed.
Recently
geometrical
interpretations of the quantum invariants
have received  interests  since
Kashaev  observed an  intriguing relationship~\cite{Kasha96b}
between
the hyperbolic
volume and his knot invariant,
which  is later identified with a specific value of the
$N$-colored Jones polynomial
$J_{\mathcal{K}}
\left(
  N ; \E^{2 \pi \I/N}
\right)$~\cite{MuraMura99a}
(here the $N$-colored Jones polynomial is normalized to be
$J_{\text{unknot}}(N;q)=1$).
Namely the hyperbolic volume of the knot complement
$S^3 \setminus \mathcal{K}$
is conjectured to dominate the asymptotics of the invariant
$J_\mathcal{K}\left(N ; \E^{2 \pi \I  /  N}\right)$ in the
large-$N$ limit
$N\to \infty$,
\begin{equation}
  \lim_{N\to \infty} \frac{2 \, \pi}{N}
  \log \left| J_\mathcal{K}
    \left( N; \E^{2 \pi \I / N} \right) \right|
  = \Vol(S^3 \setminus \mathcal{K})
\end{equation}
This ``volume conjecture'' is generalized for other values (near the
$N$-th root of unity) of the $N$-colored Jones polynomial~\cite{SGuko03a}
(see also Refs.~\citen{HMuraka04b,HMuraka06a}), and a relationship with the
A-polynomial is conjectured;
when we define $b$ by
\begin{equation}
  \label{Gukov_conjecture}
  b
  = -
  \frac{\mathrm{d}}{\mathrm{d} a} \,
  \lim_{
    \substack{
      N, k \to \infty
      \\
      N/k = a}}
  \frac{1}{k} \,
  \log J_{\mathcal{K}} \left( N; \E^{2 \pi \I/k} \right)
\end{equation}
the pair
$\left(\E^b, - \E^{\I a} \right)$ is a zero locus of the A-polynomial for
knot $\mathcal{K}$.
This is checked numerically for twist knots~\cite{KHikami04c}.

The A-polynomial is  defined as an algebraic curve
of eigenvalues of the $SL(2; \mathbb{C})$ representation of the
boundary torus of knot $\mathcal{K}$~\cite{CCGLS94a}
(see also Ref.~\citen{CoopLong98a}).
This can be computed from triangulation of the knot
complement $\mathcal{M}=S^3 \setminus \mathcal{K}$ into  ideal tetrahedra in the
hyperbolic space $\mathbb{H}^3$ once
the fundamental group has a irreducible representation $\rho$ into
$PSL(2;\mathbb{C})$ which
is identified with the
orientation-preserving isometries of $\mathbb{H}^3$;
\begin{equation*}
  \rho:
  \pi_1(\mathcal{M}) \to PSL(2; \mathbb{C}) \simeq
  \Isom^+(\mathbb{H}^3)
\end{equation*}
Up to conjugation, the meridian $\mu$ and the longitude of the
boundary torus of $\mathcal{K}$ have
\begin{gather*}
  \rho(\mu)
  =
  \begin{pmatrix}
    m & *
    \\
    0 & 1/m
  \end{pmatrix}
  \\[2mm]
  \rho(\lambda)
  =
  \begin{pmatrix}
    \ell & *
    \\
    0 & 1/ \ell
  \end{pmatrix}
\end{gather*}

Another geometrical aspect of the $N$-colored Jones polynomial 
$J_\mathcal{K}(N; q)$ as a
relationship with the A-polynomial
is
proposed as
``AJ conjecture''~\cite{GaroTQLe03a};
recursion relation of the colored Jones polynomial with respect to
$N$ is conjectured to be related to the
A-polynomial $A_\mathcal{K}(\ell, m)$.
This conjecture is proved for the torus knots~\cite{KHikami04a} and
the 2-bridge knots~\cite{TQLe04a}.

In this paper,
we introduce quantum invariant for cusped hyperbolic manifold $\mathcal{M}$ and
its deformation $\mathcal{M}_u$ {\`{a}} la
Thurston~\cite{WPThurs80Lecture}
following
Refs.~\citen{Hikam00d,KHikami01c,KHikami01e}, and study a classical
limit thereof.
Based  on  a triangulation of cusped  3-manifold $\mathcal{M}$, 
we define quantum invariant $Z_\gamma(\mathcal{M}_u)$
by assigning Faddeev's quantum dilogarithm function to  each oriented
ideal tetrahedron.
Originally Kashaev introduced  his invariant
$J_\mathcal{K}\left(
  N ; \E^{2 \pi \I/N}
\right)$ for triangulated
3-manifolds, although the $R$-matrix construction is developed
subsequently~\cite{Kasha95}.
He  studied Faddeev's  quantum dilogarithm function when $q$ is a root
of
unity~\cite{FaddKash94,Kasha94a}
(see also Ref.~\citen{BazhaReshe95}),
and assigning the quantum dilogarithm function to ideal tetrahedron
he defined invariant~\cite{Kasha95,RKasha95b}.
In this sense,
our invariant $Z_\gamma(\mathcal{M})$ for $\mathcal{M}$ with complete
hyperbolic structure 
can be regarded as  a non-compact
$U_q(sl(2;\mathbb{R}))$
analogue of the Kashaev invariant,
and an asymptotic behavior of $Z_\gamma(S^3 \setminus \mathcal{K})$ in
the limit
$\gamma\to 0$ should coincide with that of the Kashaev invariant
$J_{\mathcal{K}}\left(N; \E^{2 \pi \I/N}\right)$
in
the limit $N\to\infty$.

This paper is organized as follows.
In Section~\ref{sec:dilog} we recall definitions of quantum
dilogarithm function, and
we discuss properties of this function.
In Section~\ref{sec:invariant},
we shall reveal that the
three-dimensional hyperbolic geometry naturally arises from the
quantum dilogarithm function in the
classical limit $\gamma\to 0$
as was clarified in Ref.~\citen{Hikam00d}.
Then we define the quantum invariant
$Z_\gamma(\mathcal{M}_u)$ for a deformation of complete hyperbolic
cusped 3-manifold $\mathcal{M}$.
Based on a triangulation of $\mathcal{M}$, we construct the invariant
by assigning quantum dilogarithm function to oriented ideal
tetrahedron.
We discuss that the Neumann--Zagier potential function  appears
in a
classical limit of $Z_\gamma(\mathcal{M})$.
In Section~\ref{sec:example}, we take several examples of cusped
hyperbolic manifolds such as complements of hyperbolic knots and
punctured torus bundle over the circle,
and explain our assertion in detail.
We shall also give a list for other manifolds in Appendix.
The last section is devoted to conclusions and discussions.

\section{Quantum Dilogarithm Function}
\label{sec:dilog}
We define a function  $\Phi_\gamma(\varphi)$ by an integral form
following Ref.~\citen{LFadd99b}.
We set
$\gamma \in \mathbb{R}$, and for 
$ | \Im \varphi | < \pi$, we define
\begin{equation}
  \label{Faddeev_integral}
  \Phi_\gamma(\varphi)
  =
  \exp
  \left(
    \int\limits_{\mathbb{R}+\I \, 0}
    \frac{\E^{- \I \, \varphi \, x}}
    {
      4 \sinh ( \gamma \, x) \, \sinh(\pi \, x)
      } \,
    \frac{\mathrm{d} x}{x}
  \right) 
\end{equation}
The Faddeev integral~\eqref{Faddeev_integral}, which we call the
quantum dilogarithm function, is also related to the
double sine function~\cite{KharLebeSeme01a,TShinta77a,NKurok91a},
the hyperbolic gamma function~\cite{Ruijse97a,EBarn01a}
and the quantum exponential function~\cite{SLWoron00a}.
We see that 
the integral  $\Phi_\gamma(\varphi)$ has a duality,
\begin{equation}
  \Phi_{\frac{\pi^2}{\gamma}} ( \varphi)
  =
  \Phi_{\gamma} \left(
    \frac{\gamma}{\pi} \, \varphi
  \right) 
\end{equation}
and that it satisfies the inversion relation,
\begin{equation}
  \label{inversion_Phi}
  \Phi_\gamma( \varphi ) \cdot \Phi_\gamma( - \varphi)
  =
  \exp
  \left(
    - \frac{1}{2 \, \I \, \gamma}
    \Bigl(
    \frac{\varphi^2}{2} + \frac{\pi^2 + \gamma^2}{6}
    \Bigr)
  \right) 
\end{equation}
The Faddeev integral satisfies the difference equations;
\begin{subequations}
  \begin{align}
    \frac{\Phi_\gamma(\varphi + \I \, \gamma)}
    {    \Phi_\gamma(\varphi - \I \, \gamma)}
    & =
    \frac{1}{1 + \E^\varphi}         
    \\[2mm]
    \frac{ \Phi_\gamma(\varphi + \I \, \pi)}{
      \Phi_\gamma(\varphi - \I \, \pi)}
    & =
    \frac{1}{1 + \E^{\frac{\pi}{\gamma}\varphi}} 
  \end{align}
\end{subequations}
Due to these relations, the integral $\Phi_\gamma(\varphi)$ defined
in~\eqref{Faddeev_integral}
is
analytically continued to $\varphi \in \mathbb{C}$,
and we see that
  \begin{equation}
    \text{zeros of $\bigl( \Phi_\gamma(\varphi) \bigr)^{\pm 1}$}
    =
    \Bigl\{
    \varphi =
    \mp \I \,
    \bigl(
    (2 \, m +1 ) \, \gamma +
    ( 2 \, n + 1 ) \, \pi
    \bigr)
    \ \big| \
    m, n \in \mathbb{Z}_{\geq 0}
    \Bigr\}
  \end{equation}

The most important properties of the Faddeev integral is that
it fulfills the pentagon
identity~\cite{FaddKash94,LFadd99b}
\begin{equation}
  \label{pentagon_Phi}
  \Phi_\gamma(\Hat{p}) \, \Phi_\gamma(\Hat{q})
  =
  \Phi_\gamma(\Hat{q}) \, \Phi_\gamma(\Hat{p} + \Hat{q}) \,
  \Phi_\gamma(\Hat{p})  
\end{equation}
where $\Hat{p}$ and $\Hat{q}$ are the canonically conjugate
operators satisfying the Heisenberg commutation relation,
\begin{equation}
  \label{canonical}
  [\Hat{p} ~,~ \Hat{q} ]
  = \Hat{p} \, \Hat{q} - \Hat{q} \, \Hat{p}
  = - 2 \, \I \, \gamma   
\end{equation}
By this commuting relation, we call a limit $\gamma\to 0$  a
classical limit.

We use
$\mathbf{V}$ as the momentum space $|p\rangle$ with
$p\in\mathbb{R}$ which is  an eigenstate of the momentum operator;
\begin{equation}
  \Hat{p} \, | p \rangle
  =
  p \, | p \rangle .
\end{equation}

A reason of the \emph{quantum} dilogarithm function reveals when
we take a \emph{classical} limit $\gamma\to 0$.
In this limit, the Faddeev integral reduces to
\begin{equation}
  \label{Phi_gamma_0}
  \Phi_\gamma(\varphi)
  \sim
  \exp
  \left(
    \frac{1}{2 \,  \I \, \gamma} \,
    \Li(- \E^\varphi)
  \right) 
%
\end{equation}
where $\Li(x)$ denotes the Euler dilogarithm function defined by
(see, \emph{e.g.}, Refs.~\citen{LLewi81a,Kiril94})
\begin{equation}
  \label{Euler_dilogarithm}
  \Li(x)
  =
  \sum_{n=1}^\infty \frac{x^n}{n^2} 
\end{equation}
where $|x| \leq 1$.
For $x \in \mathbb{C}$, we use the integral form
\begin{equation*}
  \Li(x) =
  -
  \int_0^x  \log(1-s)  \, \frac{\mathrm{d} s}{s} 
\end{equation*}
where the branch of $\log(1-s)$ is on
 $\mathbb{C}\setminus [ 1, \infty)$
for which $\log(1-0)=0$.
See that the inversion relation~\eqref{inversion_Phi} of the quantum
dilogarithm function $\Phi_\gamma(\varphi)$ gives that of the Euler
dilogarithm function as
  \begin{equation*}
    \Li(- \E^x) + \Li(-\E^{-x})
    +
    \frac{x^2}{2} + \frac{\pi^2}{6} = 0 
  \end{equation*}

We note that
the Fourier transformation of the Faddeev integral can be computed as
follows~\cite{SLWoron00a,FaddKashVolk00a,BytskTesch02a,KharLebeSeme01a};
\begin{gather}
  \label{Fourier_transform_1}
  \frac{1}{\sqrt{4 \, \pi \, \gamma}}
  \int\limits_{\mathbb{R}}
  \mathrm{d} y \ \Phi_\gamma(y) \,
  \E^{\frac{1}{2 \, \I \, \gamma} \, x \, y}
  =
  \Phi_\gamma(-x + \I \pi + \I \gamma)
  \,
  \E^{
    \frac{1}{2 \, \I \, \gamma}
    \left(
      \frac{x^2}{2}
      - \frac{1}{2} \pi \gamma
      -
      \frac{\pi^2 + \gamma^2}{6}
    \right)
  } 
  \\[2mm]
  \label{Fourier_transform_2}
  \frac{1}{\sqrt{4 \, \pi \, \gamma}}
  \int\limits_{\mathbb{R}} \mathrm{d} y \
  \frac{1}{\Phi_\gamma(y)} \,
  \E^{
    \frac{1}{2 \, \I \, \gamma} x \, y
  }
  =
  \frac{1}{
    \Phi_\gamma(x - \I\, \pi - \I\, \gamma)
  } \,
  \E^{
    -\frac{1}{2 \, \I \, \gamma}
    \left(
      \frac{x^2}{2}
      - \frac{1}{2} \pi \gamma
      - \frac{\pi^2+\gamma^2}{6}
    \right)
  } 
\end{gather}

To see a relationship between the integral $\Phi_\gamma(\varphi)$ with
geometry,
we   define the $S$-operator acting  on $\mathbf{V}\otimes\mathbf{V}$ by
\begin{equation}
  \label{define_S_operator}
  S_{1,2}
  =
  \E^{
    \frac{1}{2 \,  \I\, \gamma} \,
    \Hat{q}_1 \, \Hat{p}_2
    } \,
  \Phi_\gamma (\Hat{p}_1 + \Hat{q}_2 - \Hat{p}_2) 
\end{equation}
Here  the Heisenberg operators $\Hat{p}_j$
and $\Hat{q}_j$ act on the $j$-th  vector
space of $\mathbf{V}\otimes\mathbf{V}$,
\emph{i.e.}
$\Hat{p}_1 = \Hat{p} \otimes 1$,
$\Hat{p}_2 = 1 \otimes \Hat{p} $, and so on.
Then the pentagon identity~\eqref{pentagon_Phi} can be
rewritten  in a compact form;
\begin{equation}
  \label{pentagon_S}
  S_{2,3} \, S_{1,2}
  =
  S_{1,2} \, S_{1,3} \, S_{2,3} 
\end{equation}
where $S_{j,k}$ acts on the $j$- and $k$-th spaces of
$\mathbf{V}\otimes \mathbf{V} \otimes \mathbf{V}$.
Matrix elements can be computed by use of~\eqref{Fourier_transform_1}
and~\eqref{Fourier_transform_2};
\begin{subequations}
  \label{S_element_2}
  \begin{gather}
    \begin{split}
      & \langle p_1, p_2  \ | \ S_{1, 2}  \ | \ p_1^\prime , p_2^\prime \rangle
      \\
      &  \qquad =
      \frac{1}{\sqrt{4 \, \pi \, \gamma}} \,
      \delta(p_1 + p_2 - p_1^\prime) \cdot
      \Phi_\gamma(p_2^\prime - p_2
      + \I \, \pi + \I \, \gamma) \,
      \E^{\frac{1}{2 \, \I \, \gamma}
        \left(
          - \frac{\pi^2 + \gamma^2}{6} - \frac{\gamma \, \pi}{2}
          + p_1 \, (p_2^\prime - p_2)
        \right)
        } 
    \end{split}
    \\[4mm]
    \begin{split}
      & \langle p_1, p_2 \ | \ S_{1, 2}^{~-1} \ | \ p_1^\prime , p_2^\prime
      \rangle
      \\
      &  \qquad =
      \frac{1}{\sqrt{4 \, \pi \, \gamma}} \,
      \delta(p_1 - p_1^\prime - p_2^\prime) \,
      \frac{1}{
        \Phi_\gamma(p_2 - p_2^\prime - \I \, \pi
        - \I \, \gamma)
        } \,
      \E^{
        \frac{1}{2 \I \gamma}
        \left(
        \frac{\pi^2 + \gamma^2}{6}
        + \frac{\gamma \, \pi}{2}
        - p_1^\prime ( p_2 - p_2^\prime)
        \right)
        } 
    \end{split}
  \end{gather}
\end{subequations}

In the classical limit $\gamma \to 0$, we find 
by use of \eqref{Phi_gamma_0}
that the $S$-operators~\eqref{S_element_2} 
reduce to
\begin{subequations}
  \label{S_asymptotics}
  \begin{align}
    \langle p_1 , p_2 \ | \ S_{1,2}  \ | \ p_1^\prime , p_2^\prime \rangle
    & \sim
    \delta(p_1 + p_2 - p_1^\prime) \cdot
    \exp
    \left(
      - \frac{1}{2 \, \I \, \gamma}
      \, V(p_2^\prime - p_2 , p_1)
    \right) 
  \\[2mm]
%
  \langle p_1 , p_2 \  | \ S_{1,2}^{-1} \ | \ p_1^\prime , p_2^\prime \rangle
  &
  \sim
  \delta(p_1 - p_1^\prime - p_2^\prime) \cdot
  \exp
  \left(
    \frac{1}{ 2 \, \I \, \gamma} \,
    V(p_2 - p_2^\prime , p_1^\prime)
  \right) 
\end{align}
\end{subequations}
where we have defined  the function $V(x,y)$ by
\begin{equation}
  V(x,y)
  =
  \frac{\pi^2}{6} - \Li(\E^x) - x \, y .
\end{equation}
We see that the function $V(x,y)$ satisfies partial differential
equations;
\begin{gather}
  \label{V_partial}
  V(x,y)
  =
  L(1- \E^x)
    +
    \frac{1}{2}
    \left(
      x \, \frac{\partial V(x,y)}{\partial x}
      +
      y \, \frac{\partial V(x,y)}{\partial y}
    \right) 
    \\[3mm]
    \label{imaginary_V_saddle}
    \Im \, V(x,y)
    =
    D(1-\E^x)
    +
    \log | \E^x | \cdot
    \Im \, \biggl( \frac{\partial}{\partial x} V(x,y) \biggr) 
    +
    \log | \E^y | \cdot
    \Im \, \biggl( \frac{\partial}{\partial y} V(x,y) \biggr) 
  \end{gather}
Here
the Rogers dilogarithm $L(z)$ and the
Bloch--Wigner  function $D(z)$
are respectively defined  in terms of the Euler dilogarithm
function~\eqref{Euler_dilogarithm} by
\begin{align}
  L(z)
  & =
  \Li(z) + \frac{1}{2} \log z \, \log(1-z) 
  \\[3mm]
  \label{Bloch_Wigner}
  D(z)
  & =
  \Im \, \Li(z) + \arg(1-z) \cdot \log |z| 
\end{align}
both of which fulfill the pentagon identity
(see, \emph{e.g.}, Ref.~\citen{LLewi81a});
\begin{gather}
  \label{Rogers_pentagon}
  L(z) - L(w)
  +
  L \left(\frac{w}{z} \right)
  - L\left( \frac{1-z^{-1}}{1- w^{-1}} \right)
  + L\left( \frac{1-z}{1-w} \right)
  =
  \frac{\pi^2}{6}
  \\[2mm]
  \label{Bloch-Wigner_pentagon}
   D(z) - D(w)
   +
   D \left(\frac{w}{z} \right)
   - D\left( \frac{1-z^{-1}}{1- w^{-1}} \right)
   + D\left( \frac{1-z}{1-w} \right)
   =0
\end{gather}

\section{Quantum Invariant and Potential Function}
\label{sec:invariant}

The result that the $S$-operator satisfies the pentagon
identity~\eqref{pentagon_S} indicates that we can assign 
(oriented) tetrahedron to each $S$-operator.
See Refs.~\citen{ChekFock99a,RKasha01a}
(also Ref.~\citen{ChekPenn03a})
for another interpretation as a
quantization of the Teichm{\"u}ller theory.
In the classical limit $\gamma\to 0$,
imaginary part of
the $S$-operator gives the Bloch--Wigner function
$D(z)$ which
denotes  the
hyperbolic volume of the ideal tetrahedron
$\Delta(z)$
with totally geodesic faces
which has a modulus $z$
(see, \emph{e.g.}, Refs.~\citen{BenedPetro92a,JMiln82a,WPThurs80Lecture}).
Extra terms in~\eqref{imaginary_V_saddle}
including partial differentials of $V(x,y)$
can be neglected when we assume a
saddle point condition as will be discussed later,
and 
we can assign the oriented hyperbolic ideal
tetrahedra to these $S$-operators as follows~\cite{Hikam00d};
\begin{align}
  \left\langle
    a, b \  \middle| \ S \  \middle| \  c , d 
  \right\rangle
  & =
  \mbox{
    \raisebox{-1.6cm}{
      \includegraphics[scale=1.6]{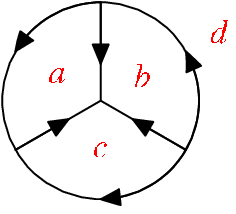}
    }
  }
  =
  \mbox{
    \raisebox{-1.6cm}{
      \includegraphics[scale=1.0]{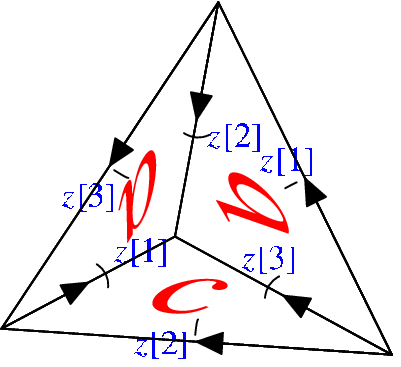}
    }
  }
  \label{S_tetrahedron}
  \\[4.8mm]
%
  \left\langle
    a, b \  \middle| \ S^{-1} \  \middle| \  c, d
  \right\rangle
  & =
  \mbox{
    \raisebox{-1.6cm}{
      \includegraphics[scale=1.6]{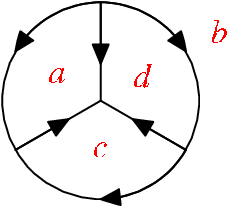}
    }
  }
  =
  \mbox{
    \raisebox{-1.6cm}{
      \includegraphics[scale=1.0]{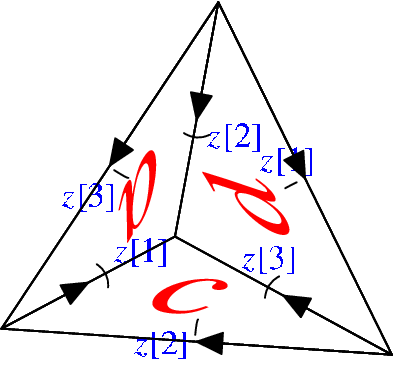}
    }
  }
  \label{S_tetrahedronDual}
\end{align}
Namely in the momentum space  representation 
we assign four
momentum parameters $\{ a , b, c , d\}$ to faces of tetrahedron.
This differs from the (quantum) $6j$ symbol as a solution of the pentagon
identity~\cite{PonzaTRegg68a,KirilReshe89}, which was used to define
the topological (quantum) gravity.

\begin{figure}
  \centering
  \includegraphics[scale=1.2]{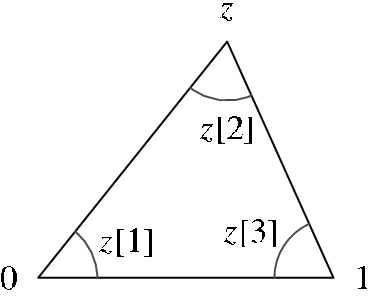}
  \caption{Triangle with vertices $0$, $1$, and $z$ in
    $\mathbb{C}$. Here we set
    $ z[1] = z$, $z[2] = 1 - \frac{1}{z}$, and $z[3] = \frac{1}{1-z}$.
  }
  \label{fig:triangle}
\end{figure}

Here both ideal tetrahedra have the modulus $z = \E^{d-b}$,
and we mean that the dihedral angles are given by
\begin{align*}
  z[1] & = z=  \E^{d-b} 
  \\[2mm]
  z[2] & = 1- \frac{1}{z} 
  \\[2mm]
  z[3] & = \frac{1}{1-z} 
\end{align*}
The cross sections by the horosphere
is similar to the triangle
in $\mathbb{C}$ with vertices $0$, $1$, and $z$ (see
Fig.~\ref{fig:triangle}),
and we have
$z[1] \, z[2] \, z[3] = -1$.
See that the opposite edges of tetrahedra have the same dihedral angles.
Then  the hyperbolic volume of the ideal tetrahedron
$\Delta(z)$ with modulus $z$ is
given by~\cite{WPThurs80Lecture}
\begin{align}
  \Vol \left( \Delta(z) \right)
  & =
  D(z)
  \nonumber
  \\
  & =
  \Lob \left(\arg(z[1]) \right)
  + \Lob \left(\arg(z[2]) \right)
  + \Lob \left( \arg(z[3]) \right)
\end{align}
where $\Lob(\theta)$  is the Lobachevsky function defined by
\begin{equation}
  \Lob(\theta)
  =
  \frac{1}{2} \sum_{n=1}^\infty \frac{\sin(2 \, n \, \theta)}{n^2}
\end{equation}

Any hyperbolic cusped $3$-manifold $\mathcal{M}$
(for simplicity we assume that the manifold $\mathcal{M}$ has one cusp
in this paper)
can be ideally
triangulated, and
it is  constructed from finite
number of the oriented ideal tetrahedra in~\eqref{S_tetrahedron}
and~\eqref{S_tetrahedronDual}.
Other types of oriented ideal tetrahedra,
such as one that has a face with loop,
are prohibited.
So   our triangulation sometimes differ from
the canonical triangulation used in computer programs, such as
SnapPea~\cite{SnapPea99}, Knotscape~\cite{Knotscape99}, and
Snap~\cite{CouGooHodNeu00a}.

To construct
3-manifolds from these ideal  tetrahedra, we glue each face to
match orientations of edges.
Every faces have no loops, so there is no ambiguity in gluing faces
together once two faces to be glued are fixed.
Furthermore,
one finds that in-state $|p \rangle$ can glue only with
out-states
$\langle p|$,
and that we cannot glue two in-states or out-states together.
As suggested by the momentum representation of the pentagon identity,
gluing  two faces, which have same label $p$,
can be formulated by an
integration with respect to $p$,
$\int\limits_\mathbb{R} \mathrm{d} p \, | p \rangle \, \langle p |$.

\begin{figure}
  \centering
  \includegraphics[scale=0.48]{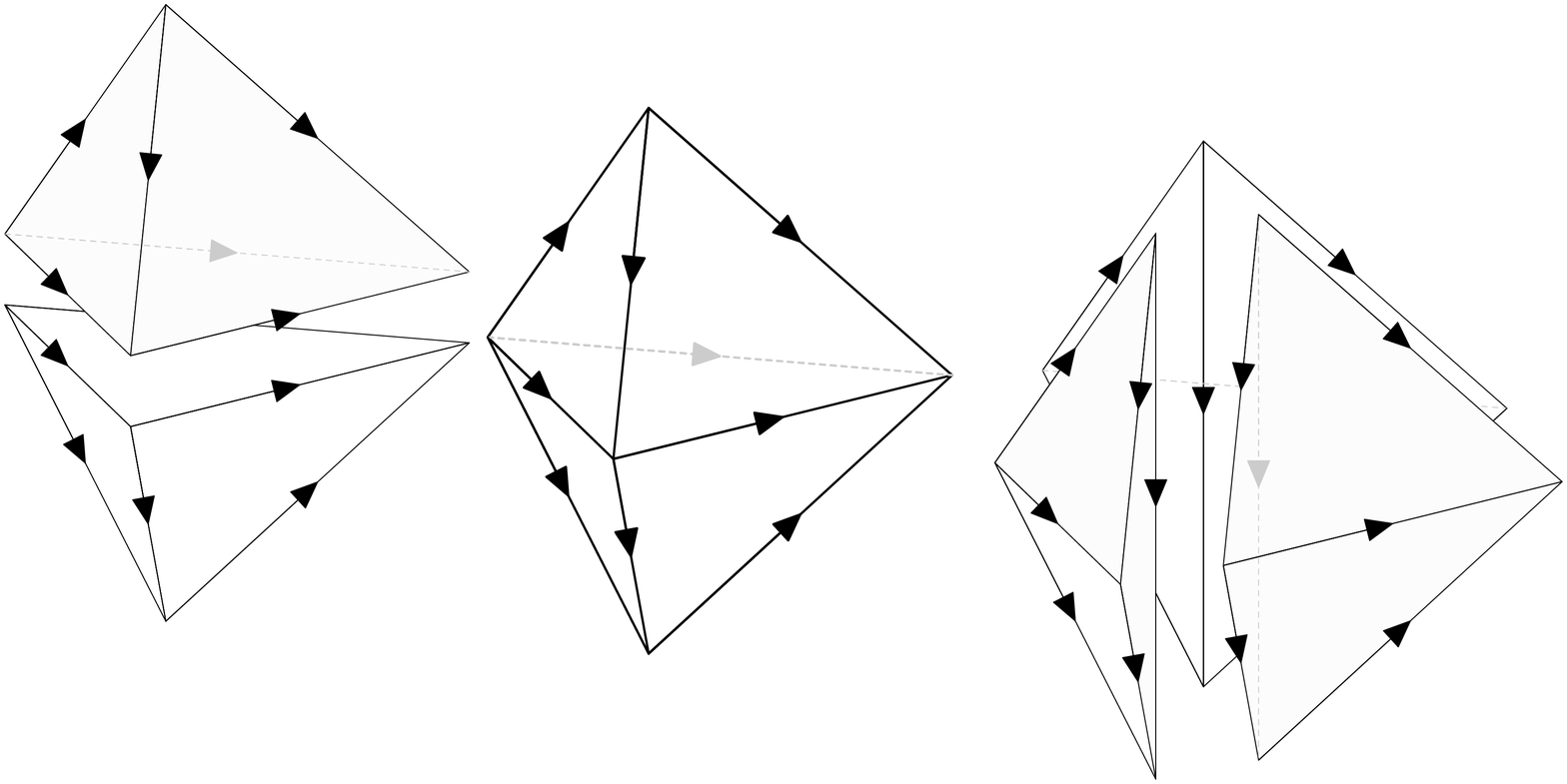}
  \caption{Pentagon identity~\eqref{pentagon_S} is interpreted as the
    Pachner move.}
  \label{fig:pentagon_S}
\end{figure}

By these identifications of the $S$-operators as the oriented ideal
tetrahedra, the pentagon identity~\eqref{pentagon_S} is identified with
the Pachner move as is depicted in
Fig.~\ref{fig:pentagon_S};
oriented polytope with 5 vertices can be decomposed into
2 tetrahedra with a face in common, or into 3 tetrahedra with an edge
in common.
In case of decomposing into 3 tetrahedra, we need a hyperbolic
consistency condition around the common edge.
This condition coincides with a saddle point condition of the right
hand side of~\eqref{pentagon_S} after substituting~\eqref{S_element_2}
and taking a classical limit $\gamma \to 0$.
One can see a coincidence  for any other type of gluings
such as other choices of orientations in the Pachner move
between the hyperbolic consistency condition
and the saddle point equations~\cite{Hikam00d}.

Once triangulation is given and we know how to glue faces together,
we can naturally define the quantum invariant for hyperbolic
cusped $3$-manifold $\mathcal{M}$
based on the
$S$-operator by
\begin{equation}
  \label{define_Z}
  Z_\gamma(\mathcal{M})
  =
  \int\limits_\mathbb{R}
  \mathrm{d} \boldsymbol{p} \,
  \delta_C(\boldsymbol{p}) \,
  \delta_G(\boldsymbol{p}) \,
  \prod_{i=1}^M 
  \left\langle
    p_{2i-1}^{(-)}, p_{2 i}^{(-)}
    ~ \middle| ~ S^{\varepsilon_i}
    ~ \middle| ~ 
    p_{2i-1}^{(+)}, p_{2 i}^{(+)}
  \right\rangle 
\end{equation}
where $\boldsymbol{p}$ denotes a set of variables
$(p_1^{(\pm)}, p_2^{(\pm)}, \dots , p_{2M}^{(\pm)})$,
and
$\varepsilon_i = \pm 1$ depending on an orientation of
tetrahedron.
We set
$M$ as  the number of ideal tetrahedra.
The condition 
$  \delta_G(\boldsymbol{p})$
determines how to glue faces together
(``G'' stands for ``gluing'').
Every faces with same momentum  mean to be glued together, and the
fact that in-states can be glued only with out-states indicates that 
the gluing condition
$  \delta_G(\boldsymbol{p})$
is  a product of
$\delta\left(p_{j}^{(-)} - p_{k}^{(+)}\right)$

We need another geometrical condition
$ \delta_C(\boldsymbol{p})$
to define invariant
besides a way how to glue faces
(``C'' stands for ``completeness'')~\cite{KHikami01e}.
We can draw a developing map from an ideal triangulation of
3-manifold,
and we need to read off  a hyperbolic complete condition.
This condition
can be written as a constraint for $\boldsymbol{p}$ by
identifications~\eqref{S_tetrahedron} and~\eqref{S_tetrahedronDual}.
By construction, the partition function
$Z_\gamma(\mathcal{M})$
is invariant under the Pachner move (Fig.~\ref{fig:pentagon_S}),
and it is essentially the invariant constructed in
Ref.~\citen{Hikam00d}.

As was studied by Thurston~\cite{WPThurs86a},
\emph{deformation}   of the hyperbolic structure on the manifold
$\mathcal{M}$ can be holomorphically parametrized by a parameter $u$
in a neighborhood of completeness condition $u=0$.
The parameter $u$
is the logarithm of the
eigenvalue of the meridian by the holonomy representation, and
we set
\begin{equation}
  \label{meridian_m_u}
  m= e^u
\end{equation}
Correspondingly  we call such  manifold $\mathcal{M}_u$ which is no
more complete,
and define the invariant by
\begin{equation}
  \label{define_deform_Z}
  Z_\gamma(\mathcal{M}_u)
  =
  \int\limits_\mathbb{R}
  \mathrm{d} \boldsymbol{p} \,
  \delta_C(\boldsymbol{p}; u) \,
  \delta_G(\boldsymbol{p}) \,
  \prod_{i=1}^M 
  \left\langle
    p_{2i-1}^{(-)}, p_{2 i}^{(-)}
    ~ \middle| ~ S^{\varepsilon_i}
    ~ \middle| ~ 
    p_{2i-1}^{(+)}, p_{2 i}^{(+)}
  \right\rangle 
\end{equation}
Here the condition
$\delta_C(\boldsymbol{p} ; u)$
follows from that the meridian has a holonomy~\eqref{meridian_m_u}.
The invariant~\eqref{define_Z} for $\mathcal{M}$ with complete
hyperbolic structure follows
from~\eqref{define_deform_Z} by assuming a completeness 
$u=0$ of the
hyperbolic structure of $\mathcal{M}_u$;
\begin{equation}
  Z_\gamma(\mathcal{M})
  =
  Z_\gamma \left( \mathcal{M}_{u=0} \right)
\end{equation}
because
\begin{equation}
  \delta_C(\boldsymbol{p})
  =
  \delta_C(\boldsymbol{p} ; u=0)
\end{equation}

We now study a classical limit of our quantum invariant
$Z_\gamma(\mathcal{M}_u)$.
When we take a  limit $\gamma \to 0$
by use of~\eqref{S_asymptotics},
the quantum invariant defined by~\eqref{define_deform_Z} becomes
\begin{align}
  Z_\gamma \left( \mathcal{M}_u \right)
  &\sim
  \int\limits_{\mathbb{R}} \mathrm{d} \boldsymbol{p} \,
  \delta_C(\boldsymbol{p}; u) \,
  \delta_G(\boldsymbol{p}) \,
  \left[
    \prod_{i=1}^M
    \delta \left(
      p_{2i-1}^{( \varepsilon_i)} - p_{2i-1}^{( - \varepsilon_i)}
      -p_{2i}^{( - \varepsilon_i)}
    \right)
  \right]
  \nonumber \\
  & \qquad \qquad 
  \times
  \exp
  \left(
    \frac{\I}{2 \, \gamma} 
    \sum_{i=1}^M
    \varepsilon_i \,
    V \left(
      p_{2i}^{(\varepsilon_i)} - p_{2i}^{( - \varepsilon_i)} ,
      p_{2i-1}^{( - \varepsilon_i)}
    \right)
  \right)
  \nonumber
  \\
  & =
  \int\limits_{\mathbb{R}} \mathrm{d} \boldsymbol{x}  \,
  \exp
  \left(
    \frac{1}{2 \, \I \, \gamma} \,
    \Phi_{\mathcal{M}}(\boldsymbol{x}; u )
  \right)
  \label{classical_partition_function}
\end{align}
Here in the last equality for our convention
we have re-parametrized variables
$\boldsymbol{p}$
with
$\boldsymbol{x}
=
\left(
  x_1, x_2, \dots, x_{M-1}
\right)$  after incorporating constraints written in
terms of delta functions.

The integral~\eqref{classical_partition_function} could be evaluated
by the saddle point method as we have taken a 
classical
limit $\gamma\to 0$.
The saddle point condition for
variables $\boldsymbol{x}$ is
\begin{equation}
  \label{hyperbolic_consistency}
  \frac{\partial }{\partial x_i}
  \Phi_\mathcal{M} (\boldsymbol{x}   ; u) = 0 
\end{equation}
As was extensively studied in Ref.~\citen{Hikam00d},
these conditions coincide with 
hyperbolic consistency conditions around 
edges when we glue oriented tetrahedra together.
By construction
the variable $u$ denotes the meridian of cusp in
this classical limit,
and the complete hyperbolic structure is  realized by setting $u=0$.
To conclude, the function $\Phi_{\mathcal{M}}(\boldsymbol{x} ; u)$
defined 
by a classical limit~\eqref{classical_partition_function} of the
quantum invariant $Z_\gamma(\mathcal{M}_u)$
is nothing but
the
Neumann--Zagier
potential function~\cite{NeumZagi85a,TYoshi85a}.

As a result, differential of the potential function with respect to
the
deformation parameter $u$ gives
\begin{equation}
  \label{differential_longitude}
  \frac{\partial}{\partial  u}  \,
  \Phi_{\mathcal{M}}(\boldsymbol{x} ; u)
  =
  2 \, v
\end{equation}
where $v$ is related to the
eigenvalue of the longitude by the holonomy representation
\begin{equation}
  \ell = \E^{v}
\end{equation}
Variables  $x_i$ can be solved from the hyperbolic consistency
equation~\eqref{hyperbolic_consistency} as a function of $u$,
and we can regard the potential function
$\Phi_\mathcal{M}(\boldsymbol{x};u)$ as a function of $u$;
$\Phi_\mathcal{M}(u) =
\left.
  \Phi_{\mathcal{M}}( \boldsymbol{x} ; u)
\right|_{\eqref{hyperbolic_consistency}}$.
Then we can rewrite~\eqref{differential_longitude} by
\begin{equation}
  \label{v_differential_limit}
  \frac{\mathrm{d}}{\mathrm{d} u}
  \lim_{\gamma \to 0 }
  \I \, \gamma \,
  \log Z_\gamma\left( \mathcal{M}_u \right)
  =  v
\end{equation}
We note that our variables $(u,v)$ differs from those in
Ref.~\citen{NeumZagi85a};
when we denote their variables $(u_{\text{NZ}}, v_{\text{NZ}})$, we
have
\begin{equation}
  \label{uv_NZ}
  (u, v) =
  \left(
    \frac{u_{\text{NZ}}}{2} ,
    \frac{v_{\text{NZ}}}{2} + \pi \, \I
  \right)
\end{equation}
Hereafter we  also use  the function
$V_{\mathcal{M}}(\boldsymbol{x}; m)$ defined by
\begin{equation}
  V_{\mathcal{M}}(x_1, x_2, \dots ; m)
  =
  \Phi_{\mathcal{M}}(\log x_1, \log x_2, \dots ;  u=\log m)
\end{equation}

As seen from~\eqref{V_partial} and~\eqref{imaginary_V_saddle}, the
potential function $\Phi_{\mathcal{M}}(\boldsymbol{x};u)$
under  saddle point conditions~\eqref{hyperbolic_consistency} becomes
a sum of the Rogers
dilogarithm functions.
We recall here the Bloch invariant 
studied in 
Refs.~\citen{NeumYang99a,NeumYang95b}.
The Bloch invariant $\beta(\mathcal{M})$ is defined for
finite volume
hyperbolic 3-manifold $\mathcal{M}$ as
\begin{equation}
  \beta(\mathcal{M}) = \sum_{i=1}^M \left[ z_i \right]
\end{equation}
where $\left[ z \right]$ satisfies  the Bloch group
\begin{equation}
  \left[ z \right] -  \left[w \right]
  +
  \left[\frac{w}{z} \right]
  - \left[ \frac{1-z^{-1}}{1- w^{-1}} \right]
  + \left[ \frac{1-z}{1-w} \right]
  =
  0
\end{equation}
The Bloch regulator map $\rho$ gives~\cite{NeumYang99a}
\begin{equation}
  \rho\left( \beta(\mathcal{M}) \right)
  =
  \Vol(\mathcal{M})
  + \I \, \CS(\mathcal{M})
\end{equation}
where $\CS$ denotes the Chern--Simons invariant
(see Ref.~\citen{RMeyer86b} for definition of the Chern--Simons invariant for a
case of cusped manifolds).
As identities~\eqref{V_partial}--~\eqref{imaginary_V_saddle} show that
the $S$-operator reduces to the Rogers dilogarithm function
(or the Bloch--Wigner function) in the saddle
point, and that they satisfy the pentagon
identities~\eqref{Rogers_pentagon} and~\eqref{Bloch-Wigner_pentagon},
we can interpret our invariant $Z_\gamma(\mathcal{M})$ as a
quantization of the Bloch invariant.

Generally we have many saddle points as algebraic solutions of a set of
equations~\eqref{hyperbolic_consistency}.
Among them, a solution which has the  largest absolute value
dominates an asymptotics of the quantum invariant $Z_\gamma(\mathcal{M})$.
Combining this with a fact that
our quantum invariant $Z_\gamma(\mathcal{M})$ may be regarded as a
quantization of the Bloch invariant,
we should have
\begin{equation}
  \lim_{\gamma \to 0}
  2 \, \gamma\,
  \log \left(
    Z_\gamma(\mathcal{M})
  \right)
  =
  \Vol(\mathcal{M})
  + \I \,
  \CS(\mathcal{M})
\end{equation}
as a variant of the volume conjecture~\cite{Hikam00d}.
Although, there still remains an ambiguity of branch in complex plane
in an actual computation.

The Neumann--Zagier
potential function $\Phi_\mathcal{M}(\boldsymbol{x};u)$ has much 
information on geometry of manifold.
One of the properties
is  a relationship with the A-polynomial  defined in
Ref.~\citen{CCGLS94a}.
When the manifold $\mathcal{M}$ is a complement of knot
$S^3 \setminus \mathcal{K}$,
the A-polynomial 
$A_\mathcal{K}(\ell, m)$ of $\mathcal{K}$
as an algebraic equation of $\ell$ and $m$
can be given by using
the Gr\"obner
base or resultant theory to eliminate $x_i$ from a set of
equations,~\eqref{hyperbolic_consistency}
and~\eqref{differential_longitude}.
So the pair $(\ell,m)=(\E^v, \E^u)$ defined
from~\eqref{v_differential_limit}
is a zero locus of the A-polynomial $A_\mathcal{K}(\ell,m)$.
This result should be comparable with a
conjecture~\eqref{Gukov_conjecture}, and the invariant
$Z_\gamma(\mathcal{M}_u)$ is a \emph{non-compact} generalization of the
Jones--Witten invariant.

The A-polynomial has following
properties~\cite{CCGLS94a,CoopLong98a};
\begin{itemize}

\item Polynomial $A_\mathcal{K}(\ell,m)$ is an integer polynomial, and
  it  contains only even powers of $m$.

\item Up to powers of $\ell$ and $m$, we have
  \begin{equation}
    A_\mathcal{K}(\ell, m) =
    A_\mathcal{K}(1/\ell, 1/m)
  \end{equation}

\item If $\mathcal{K}$ and $\mathcal{K}^\prime$ are mirror images,
  then
  \begin{equation}
    A_\mathcal{K}(\ell, m)=
    A_{\mathcal{K}^\prime}(1/\ell, m) 
  \end{equation}

\item With additional property that every closed incompressible
  surface embedded in $S^3 \setminus \mathcal{K}$ is parallel to the
  boundary torus, we have
  \begin{equation}
    A_\mathcal{K}(\ell, \pm 1 )
    =
    n \, \left( \ell+1 \right)^{\alpha_+}
    \,
    \left( \ell - 1 \right)^{\alpha_-} \,
    \ell^\beta
  \end{equation}
  with non-zero integer $n$.

  \item Slopes of edges of the Newton polygon of
    $A_\mathcal{K}(\ell,m)$ are boundary slopes of $\mathcal{K}$.

  \item Coefficients of terms in the corners of the Newton polygons of
    $A_\mathcal{K}(\ell,m)$ are $\pm 1$.

\end{itemize}

Another property of
the Neumann--Zagier potential function $\Phi_{\mathcal{M}}(u)$ associated to
cusped manifold $\mathcal{M}$ is
the volume of the Dehn surgered
manifold.
The $(p,q)$-hyperbolic Dehn surgery of $\mathcal{M}$,
where $(p,q)$ is a pair of coprime integers, is performed by
gluing back a
solid torus with cusp of $\mathcal{M}$, where the surgery data
satisfy~\cite{WPThurs80Lecture}
\begin{equation}
  p \, u+ q \, v =  \pi \, \I
\end{equation}
Then for the core $c$ of solid torus, we have
\begin{equation}
  \length(c) + \I \, \torsion(c)
  =
  - 2\, \left(
    r \, u + s\, v
  \right)
  \mod 2 \, \pi \, \I
\end{equation}
where
$
\begin{pmatrix}
  p & r \\
  q & s
\end{pmatrix}
\in SL(2,\mathbb{Z})$.
We have
\begin{equation}
  \Im \left( u \, \bar{v} \right)
  =
  -\frac{ \pi}{2} \, \length(c)
\end{equation}
According to Refs.~\citen{TYoshi85a,NeumZagi85a}, 
we have for
the hyperbolic $(p,q)$-Dehn surgered manifold $\mathcal{M}_{(p,q)}$
as
\begin{multline}
  \left( \Vol\left(\mathcal{M}_{(p,q)}\right)
    + \I \, \CS \left(\mathcal{M}_{(p,q)} \right) \right)
  -
  \left( \Vol(\mathcal{M}) + \I \, \CS(\mathcal{M}) \right)
  \\
  =
  - \frac{\I}{4} \, \left(
    \Phi_{\mathcal{M}}(u) - 4 \, u \, v
  \right) -
  \frac{\pi}{2} \,
  \left(
    \length(c)
    + \I \, 
    \torsion(c)
  \right)
\end{multline}
which follows from 
\begin{equation}
  Z_\gamma \left( \mathcal{M}_{(p,q)} \right)
  \sim
  \int \mathrm{d} u \,
  \E^{
    \frac{1}{2 \I \gamma}
    \left(
      \frac{p}{q} u^2 +
      \frac{2 ( \pi + \gamma)  \I}{q} u
    \right)
  } \,
  Z_\gamma(\mathcal{M}_u)
\end{equation}

\section{Examples}
\label{sec:example}
We explain our constructions of
the quantum invariant by taking some
concrete examples of cusped hyperbolic 3-manifolds.

\subsection{
  Figure-Eight knot $ 4_1$}
\label{sec:41}

\begin{figure}[htbp]
    \centering
    \includegraphics[scale=0.8]{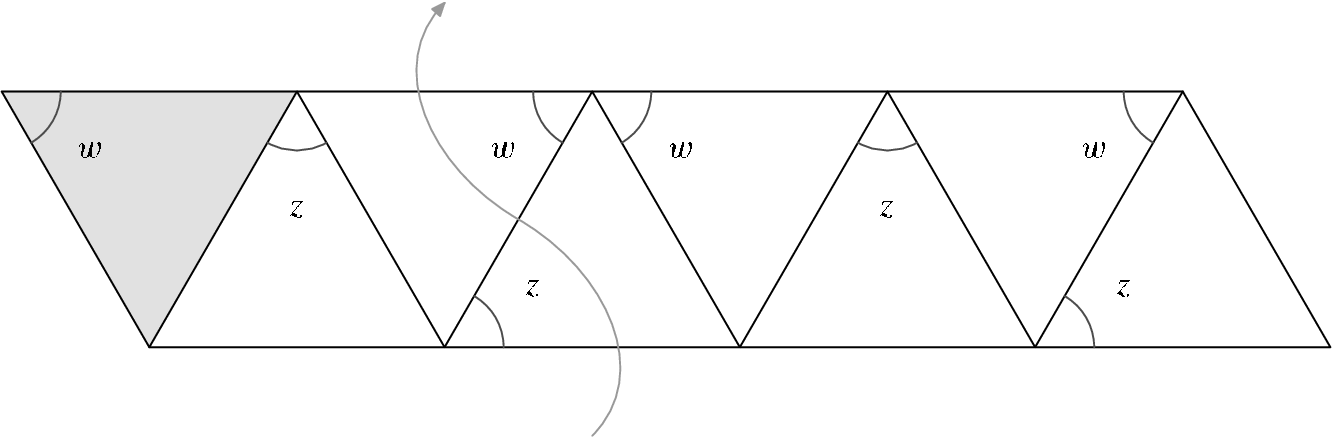}
    \caption{Developing map of the complement of the
      figure-eight knot.
      The gray filled triangle corresponds to top vertex of the
      tetrahedron
      (central vertex of circle in~\eqref{S_tetrahedron})
      with modulus $w$ in projection
      of~\eqref{S_tetrahedron}.
      A curve denotes a meridian of cusp.
  }
    \label{fig:develop41}
  \end{figure}

We set $\mathcal{K}$ as 
the Figure-Eight knot $4_1$ which is depicted
as
\begin{equation*}
  \includegraphics[scale=.16]{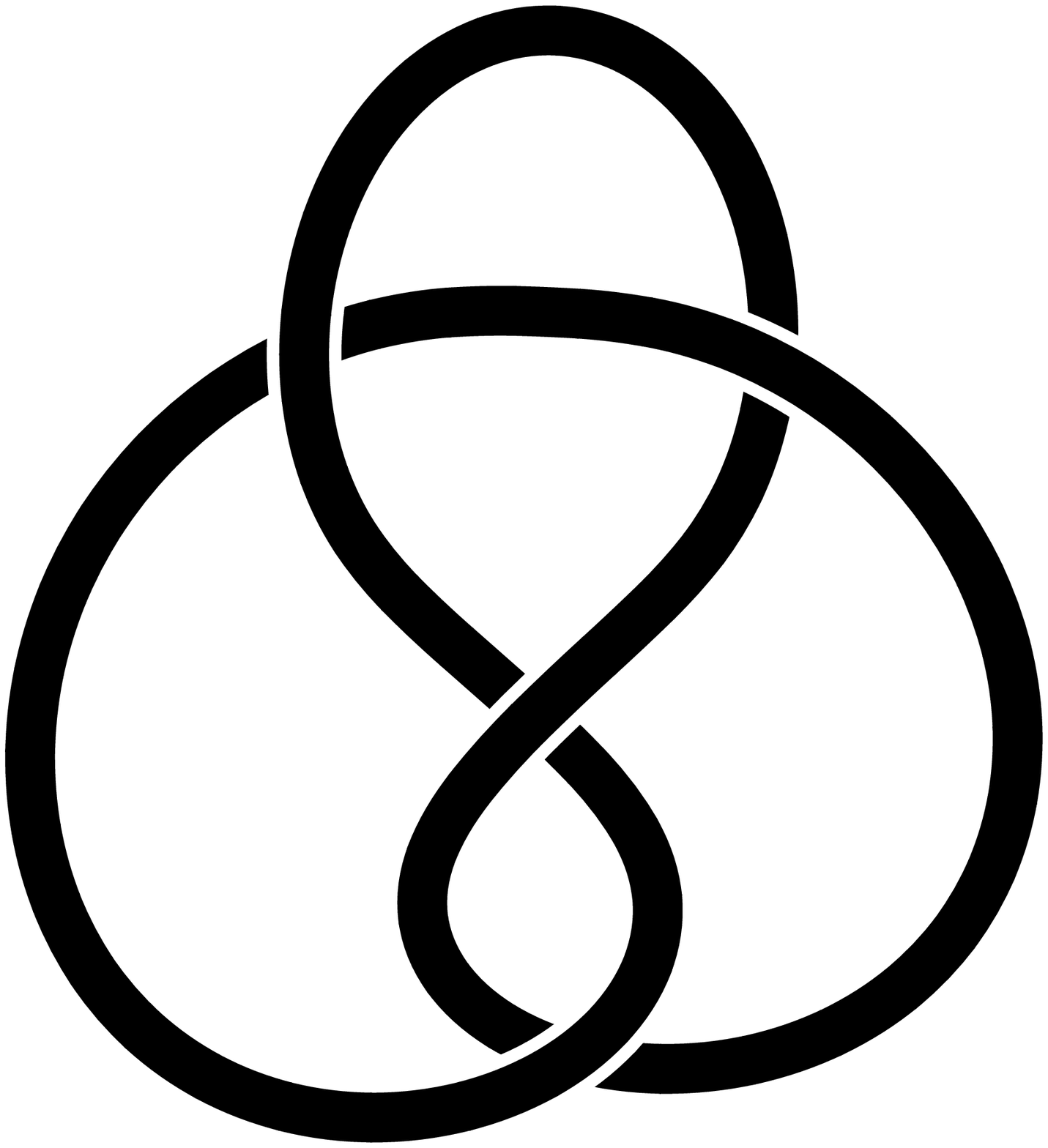}    
\end{equation*}
It is well known that the complement of the figure eight knot,
$\mathcal{M}=S^3 \setminus \mathcal{K}$,
is given by two ideal tetrahedra,
and the triangulation induces the quantum invariant as
\begin{equation}
  Z_\gamma(\mathcal{M}_u)
  =
  \int\limits_{\mathbb{R}} \mathrm{d} \boldsymbol{p} \,
  \delta_C(\boldsymbol{p}; u) \,
  \left\langle
    p_1, p_2  \middle| S \middle| p_3, p_4
  \right\rangle \,
  \left\langle
    p_4, p_3 \middle| S^{-1} \middle| p_2, p_1 \right\rangle 
\end{equation}
Modulus of two tetrahedra are given by
$w=\E^{p_4 - p_2}$
and
$z=\E^{p_1 - p_3}$.
The developing map is drawn in Fig.~\ref{fig:develop41},
and
the meridian  is read to be 
\begin{equation*}
  \frac{w}{z} = \E^{-2 u}
\end{equation*}
which shows that the condition $\delta_C(\boldsymbol{p};u)$ is
\begin{equation*}
  p_4 - p_2 - (p_1 - p_3) = -  2\,  u
\end{equation*}
The complete hyperbolic structure is realized when $u=0$.

We then obtain
  \begin{equation}
    Z_\gamma\left( \mathcal{M}_u \right)
    =
    \frac{1}{ 4 \, \pi \, \gamma}
    \int \mathrm{d} x \,
    \frac{
      \Phi_\gamma(x+ \I \, \pi + \I \, \gamma)
    }{
      \Phi_\gamma(-x-  2 \, u - \I \, \pi - \I \, \gamma)
    } \,
    \E^{\frac{-1}{2 \I \gamma} 4 u(u +x)}
  \end{equation}
  In a limit $\gamma \to 0$, we have
  \begin{align}
    Z_\gamma(\mathcal{M}_u)
    & \sim
    \int \mathrm{d} x \, 
    \exp
    \left(
      \frac{1}{2 \, \I \, \gamma} \,
      \left(
        \Li(\E^x) - \Li(\E^{-x- 2 u}) -  4 \, u \, (u+ x)
      \right)
    \right)
    \nonumber \\
    & =
    \int \mathrm{d}x \,
    \exp \left( \frac{1}{2 \, \I \, \gamma} \,
      V_{\mathcal{M}}( \E^x ; \E^u) \right)
    \label{Z_V_classical_41}
  \end{align}
  where the potential function is set to be
  \begin{equation}
    \label{41_potential}
    V_{\mathcal{M}}(x; m)
    =
    \Li(x) - \Li \left(\frac{1}{x \, m^2 }\right)
    -4   \log m \,   \log (x \, m)
  \end{equation}
  
  To evaluate the integral~\eqref{Z_V_classical_41}, we may apply the
  saddle point method, and a condition~\eqref{hyperbolic_consistency}
  reduces to
  \begin{gather}
    \label{41_saddle}
    - \frac{x}{m^2 \, (1-x) \, (1-m^2 \, x)} = 1
  \end{gather}
  From~\eqref{differential_longitude} the longitude is given by
  \begin{gather}
    \label{41_long}
    \frac{1}{m^2 \, x \, (m^2 \, x - 1)}
    =
    \ell
  \end{gather}

  Completeness condition $m=1$ gives
  $x=\frac{1 \pm \sqrt{3} \, \I}{2}$ from~\eqref{41_saddle}.
  Substituting this solution for the potential
  function~\eqref{41_potential},
  the imaginary part coincides with
  the hyperbolic volume of the figure-eight knot
  $\Vol(S^3 \setminus 4_1)
  = 2.02988 \cdots$.

  For a deformed manifold $\mathcal{M}_u$,
  we get an algebraic equation of $\ell$ and $m$ by eliminating $x$
  from~\eqref{41_saddle} and~\eqref{41_long} as
  \begin{equation}
    A(\ell , m) = 0
  \end{equation}
  Here  $A(\ell, m)$ is the A-polynomial for the figure-eight knot;
  \begin{equation}
    A(\ell,m)
    =
    - m^4 + \ell \, (1 - m^2 - 2 \, m^4 - m^6 + m^8 )
    -
    \ell^2 \, m^4
  \end{equation}
  which can be  expressed as in the Newton polygon as follows;
  \begin{equation*}
    \begin{pmatrix}
      0 & 1 & 0 \\
      0 & -1 & 0 \\
      -1 &-2 & -1 \\
      0 & -1 & 0 \\
      0 & 1 & 0
    \end{pmatrix}
  \end{equation*}

  We note that a set of equations~\eqref{41_saddle}
  and~\eqref{41_long} gives
  \begin{equation*}
    2 \, v
    =
    -2 \, \pi \, \I + 4\, \sqrt{3} \,  \I \, u
    + \frac{16 \, \I}{\sqrt{3}} \, u^3
    +
    \frac{368 \, \I}{15 \sqrt{3}} \, u^5
    +
    \frac{2848 \, \I}{45 \sqrt{3}} \, u^7
    + \cdots
  \end{equation*}
  which coincides with a result in Ref.~\citen{NeumZagi85a}
  under~\eqref{uv_NZ}.

\subsection{
  $5_2$}
\begin{figure}[htbp]
  \centering
  \includegraphics[scale=1.28]{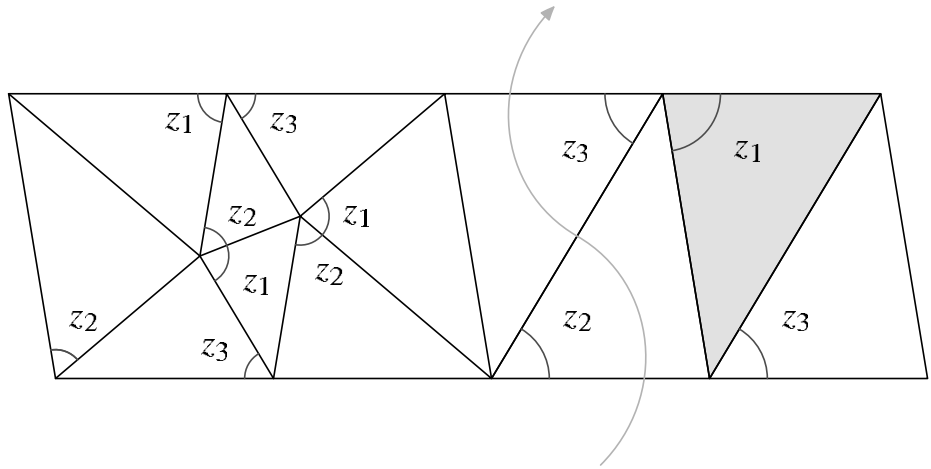}
  \caption{Developing map of the complement of $5_2$.
    Gray filled triangle denotes top vertex of the tetrahedron with
    modulus $z_1$.
    Meridian is denoted by a gray curve.
  }
  \label{fig:develop52}
\end{figure}

We set $\mathcal{K}$ as the knot $5_2$ depicted as
\begin{equation*}
  \includegraphics[scale=.16]{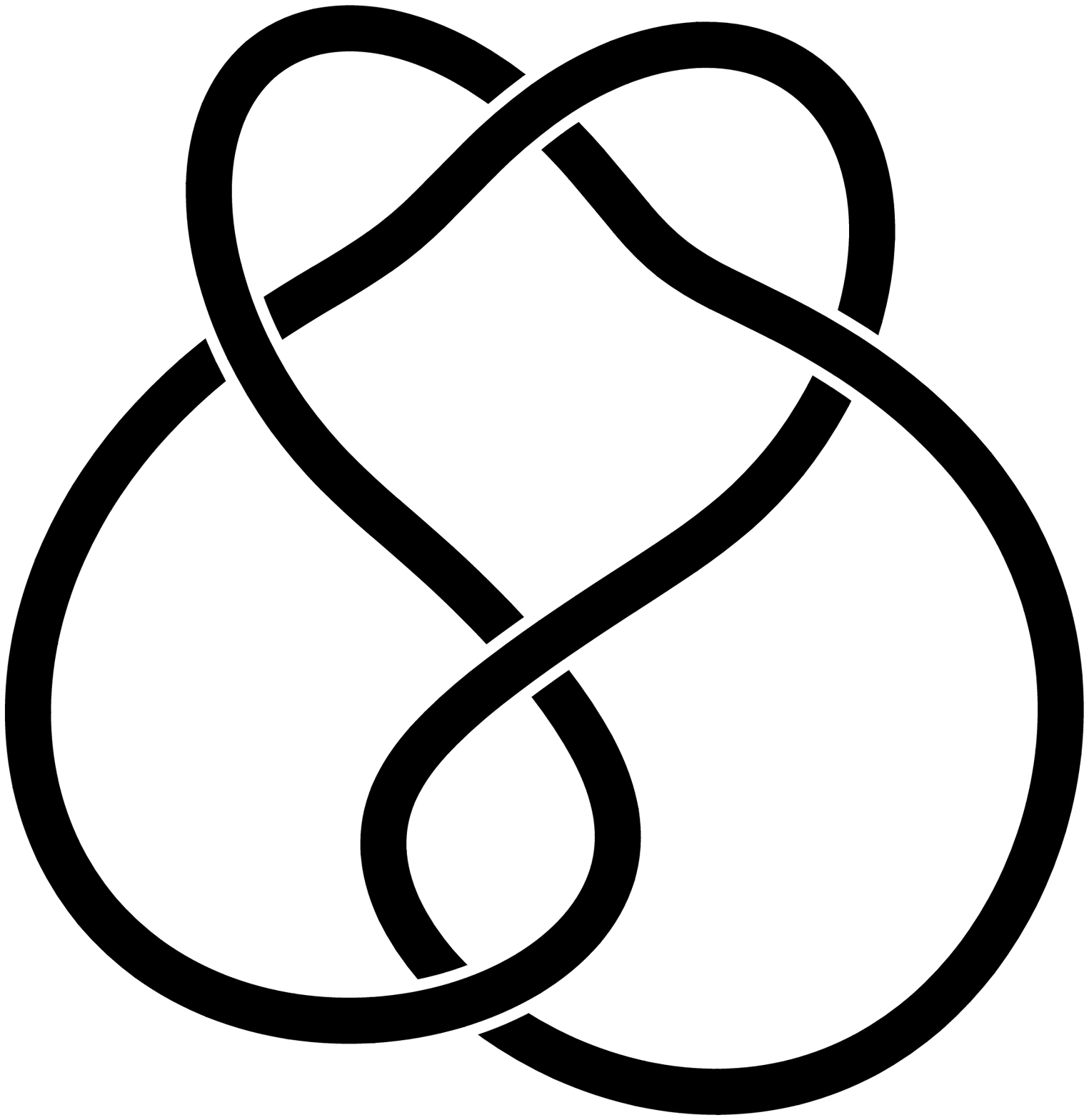}  
\end{equation*}
The knot complement $\mathcal{M}=S^3 \setminus \mathcal{K}$
is constructed from three tetrahedra (see, \emph{e.g.},
Ref.~\citen{MoTaka85a})
and we obtain the invariant as
\begin{equation}
  Z_\gamma\left( \mathcal{M}_u \right)
  =
  \int\limits_\mathbb{R} \mathrm{d}\boldsymbol{p} \,
  \delta_C(\boldsymbol{p}; u ) \,
  \left\langle
    p_1 , p_5 \middle| S^{-1} \middle| p_4 , p_3
  \right\rangle \,
  \left\langle
    p_2 , p_4 \middle| S^{-1} \middle| p_6 , p_5
  \right\rangle \,
  \left\langle
    p_3 , p_6 \middle| S^{-1} \middle| p_1 , p_2
  \right\rangle 
\end{equation}
Modulus of three tetrahedra are respectively
$z_1 =\E^{p_3 - p_5}$,
$z_2 =\E^{p_5 - p_4}$, and
$z_3 =\E^{p_2 - p_6}$.
The developing map is depicted in Fig.~\ref{fig:develop52}.
The meridian is read to be  $\frac{z_3}{z_2}$, and the condition
$\delta_C(\boldsymbol{p}; u)$ is
\begin{equation*}
  p_5 - p_4 + p_6 - p_2
  =  2\, u
\end{equation*}
We then have
\begin{multline}
  Z_\gamma\left( \mathcal{M}_u \right)
  =
  \frac{1}{
    \left( 4 \, \pi \, \gamma \right)^{3/2}
  } \,
  \iint \mathrm{d} x \, \mathrm{d} y 
  \,
  \E^{\frac{1}{2 \I \gamma}
    \left(
      \frac{\pi^2+\gamma^2}{2}
      +
      \frac{3}{2} \pi \gamma 
      +
      ( 2 u-y) (y-x)
    \right)
  }
  \\
  \times
  \frac{1}{
    \Phi_\gamma(-x+y- 2 \,  u-\I \, \pi - \I \, \gamma) \,
    \Phi_\gamma(-y- 2 \, u-\I \, \pi - \I \, \gamma) \,
    \Phi_\gamma(-y- 2 \, u-\I \, \pi - \I \, \gamma)
  }
\end{multline}

In a classical limit $\gamma\to 0$, we obtain
\begin{equation}
  Z_\gamma\left( \mathcal{M}_u \right)
  \sim
  \iint \mathrm{d} x \, \mathrm{d} y \,
  \exp 
  \left(
    \frac{1}{2 \, \I \, \gamma} \,
    V_{\mathcal{M}}(\E^x, \E^y ; \E^u)
  \right)
\end{equation}
where the Neumann--Zagier potential function is
  \begin{equation}
    V_{\mathcal{M}}(x,y; m)
    =
    \frac{\pi^2}{2}
    - \Li\left(\frac{y}{x \, m^2}\right)
    - \Li \left( \frac{1}{ y \, m^2}\right)
    - \Li \left(\frac{1}{y}\right)
    + \log \left( y/x \right) \, \log \left(m^2/y \right)
  \end{equation}
The integral is evaluated by the saddle point method, and  saddle
point conditions~\eqref{hyperbolic_consistency} which denote 
hyperbolic consistency conditions
reduce to
\begin{equation}
  \label{consistency_52}
  \begin{gathered}
    m^2 \, x - y = x \, y
    \\[2mm]
    (m^2 \, y -1) \, (y -1) = m^2 \, (m^2 \, x - y)
  \end{gathered}
\end{equation}
and
the identity for the longitude~\eqref{differential_longitude} 
is given as
\begin{equation}
  \label{longitude_52}
  m^4 \, y^2 =
  (m^2 \, x - y) \, (m^2 \, y - 1) \, \ell
\end{equation}

In the case of complete structure $\mathcal{M}$, \emph{i.e.}, $u=0$,
there exists a solution of~\eqref{consistency_52},
$(x,y)=
(-0.877439 + 0.744862 \, \I ,
0.78492 + 1.30714 \, \I)$,
 such that the
imaginary part of the potential function
\begin{equation*}
  \Im V_\mathcal{M}(x,y;1)
  =
  -D \left( \frac{y}{x} \right)
  - 2 \,D\left( \frac{1}{y} \right)
\end{equation*}
gives the hyperbolic volume of $\mathcal{M}$;
$\Vol(\mathcal{M})
=
2.82812 \cdots$.
.

To get the A-polynomial, we
eliminate variables $x$ and $y$ from a set of
equations~\eqref{consistency_52} and~\eqref{longitude_52}.
After some algebra we obtain
$A_\mathcal{K}(\ell,m)=0$, where the function
$A_\mathcal{K}(\ell,m)$ coincides with the
A-polynomial for $5_2$
given by the following Newton polygon;
  \begin{equation*}
    \begin{pmatrix}
      -1 & 1 &  &  \\
       & -2 &  &  \\
       & -2 & -1 &  \\
        &  & 1 &  \\
        & 1 &  &  \\
        & -1 & -2 &  \\
        &  & -2 &  \\
        &  & 1 & -1 \\
    \end{pmatrix}
  \end{equation*}

\subsection{
  Pretzel knot $(-2,3,7)$}

\begin{figure}[htbp]
  \centering
  \includegraphics[scale=1.04]{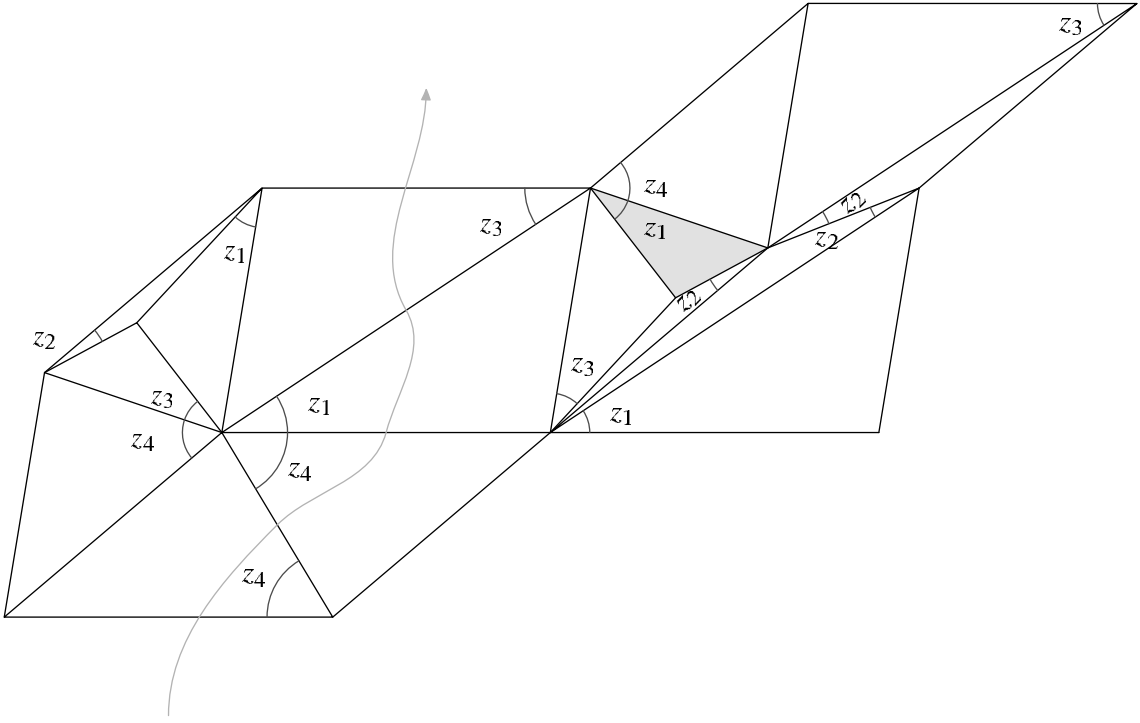}
  \caption{Developing map of  complement of the Pretzel knot
    $(-2,3,7)$.
    Horosphere of the top vertex of the tetrahedron with modulus $z_1$
    is depicted by gray triangle.
    Gray curve denotes a meridian.
  }
  \label{fig:developPretzel}
\end{figure}

Let $\mathcal{K}$ be the $(-2,3,7)$ Pretzel knot depicted as
\begin{equation*}
  \includegraphics[scale=.16]{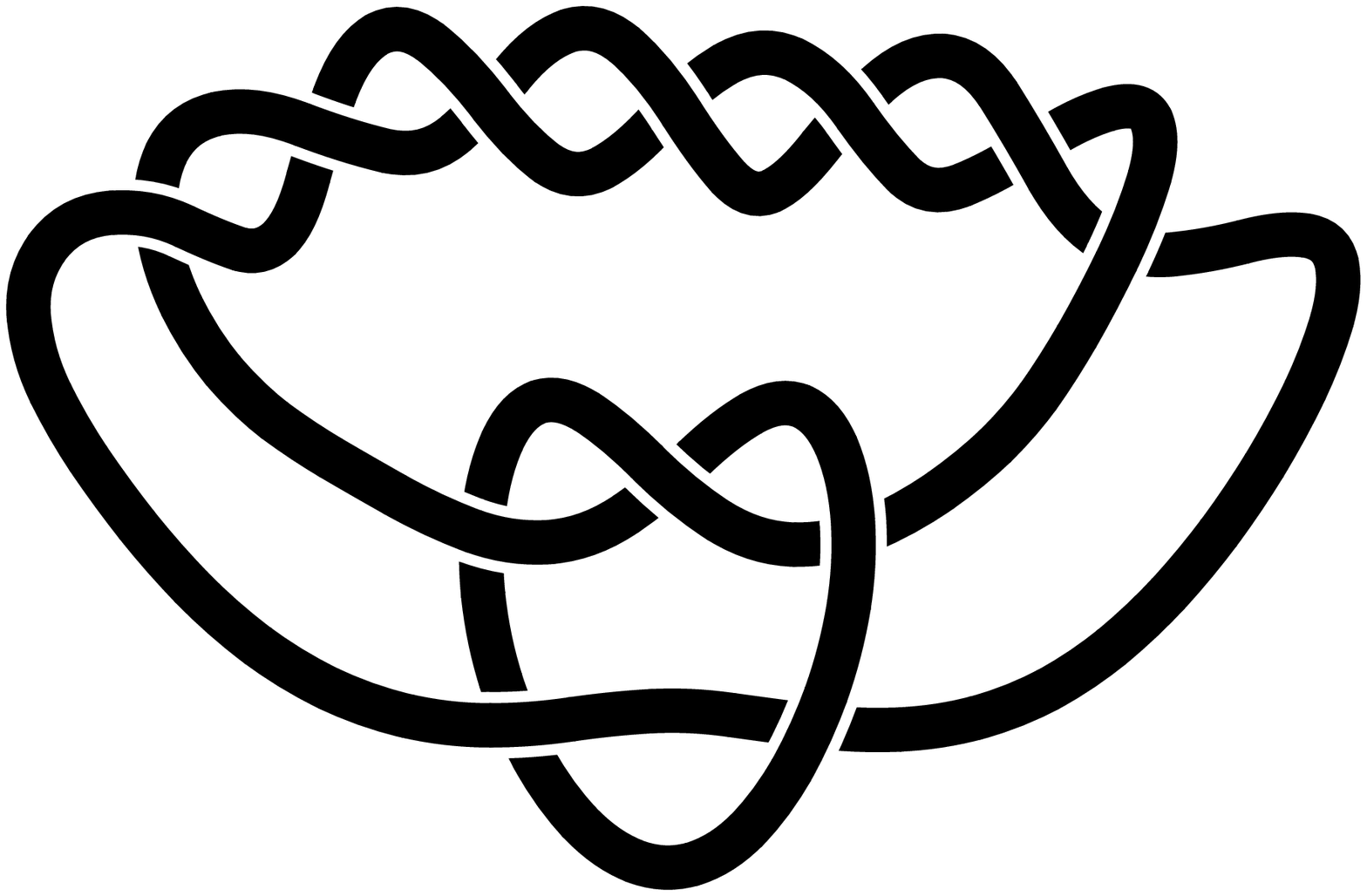}
\end{equation*}
and we set $\mathcal{M}$ as the complement of $\mathcal{K}$.
The Pretzel knot $(-2,3,7)$ has the same hyperbolic volume with $5_2$,
$\Vol(\mathcal{M})
=
2.82812 \cdots$,
but the triangulations  of the complement give the following partition
function;
\begin{multline}
  Z_\gamma\left( \mathcal{M}_u \right)
  = \int\limits_{\mathbb{R}}
  \mathrm{d} \boldsymbol{p} \,
  \delta_C(\boldsymbol{p}; u) \,
  \langle p_4 , p_7 | S | p_6 , p_1 \rangle \,
  \langle p_5 , p_8 | S | p_7 , p_5 \rangle
  \\
  \times
  \langle p_1 , p_6 | S | p_8 , p_2 \rangle \,
  \langle p_3 , p_2 | S | p_4 , p_3 \rangle 
\end{multline}
Note that this triangulation differs from the canonical triangulation
in Ref.~\citen{SnapPea99}.
The developing map is given in Fig.~\ref{fig:developPretzel}.
Here we set the modulus of four ideal tetrahedra as
\begin{align*}
  z_1 & = \E^{p_1 - p_7}
  &
  z_2 &= \E^{p_5-p_8}
  \\
  z_3 & = \E^{p_2 - p_6}
  &
  z_4 & = \E^{p_3 - p_2}
\end{align*}
Then the
meridian is read as
\begin{equation*}
   p_2 - p_6 - p_1 + p_7
   = -  2 \, u
\end{equation*}

In the classical limit of the quantum invariant $Z_\gamma(\mathcal{M}_u)$
we obtain the potential function after some change of variables as
\begin{multline}
  V_{\mathcal{M}}(x,y,z ; m)
  =    -\frac{2 \, \pi^2}{3}+
  \Li\left(\frac{1}{z}\right)
  + \Li\left(\frac{1}{x \, y \, z \, m^4}\right)
  + \Li\left(\frac{1}{z \, m^2}\right)
  + \Li\left(\frac{m^2}{x}\right)
  \\
  + 3 \left( \log \left( m^2 \right) \right)^2 +
  \log \left(\frac{y^5}{x^2} \right) \log \left( m^2 \right)
  + \left( \log x \right)^2 +
  \left( \log y \right)^2
\end{multline}

The hyperbolic consistency
conditions~\eqref{hyperbolic_consistency} give
\begin{equation}
  \label{consistency_pretzel}
  \begin{gathered}
    \left( x \, y \, z \, m^4 - 1 \right) \, \left( x- m^2 \right)
    = y \, z \, m^8
    \\[2mm]
    y \, m^6 \, \left( x \, y \, z \, m^4 - 1 \right)
    = x \, z
    \\[2mm]
    (z-1) \, \left(
      x \, y \, z \, m^4 - 1 \right) \, \left( m^2 \, z -1 \right)
    =
    x \, y \, z^3 \, m^6
  \end{gathered}
\end{equation}
and the longitude~\eqref{differential_longitude} is defined by
\begin{equation}
  \label{long_pretzel}
  \ell
  =
  \frac{
    \left( x \, y \, z \, m^4 -1 \right)^2 \,
    \left( m^2 \, z - 1 \right)^2 \, m^2 \, y^3
  }{
    \left( x - m^2 \right) \, x^3 \, z^3}
\end{equation}
In the complete hyperbolic structure $m=1$,
we have a solution of~\eqref{consistency_pretzel},
\begin{equation*}
  \begin{pmatrix}
    x \\ y \\ z
  \end{pmatrix}
  =
  \begin{pmatrix}
    0.337641 - 0.56228 \, \I \\
    0.122561 + 0.744862 \, \I \\
    0.618504 - 0.410401 \, \I
  \end{pmatrix}
\end{equation*}
s.t.
we recover the hyperbolic volume of $\mathcal{K}$;
\begin{equation*}
  \Im V_\mathcal{M}(x,y,z; m=1)
  =
  2 \, D(1/z) + D\left(\frac{1}{x \, y \, z}\right)
  + D(1/x)
  =
  2.82812 \cdots
\end{equation*}

The A-polynomial is computed by eliminating
$(x,y,z)$ from~\eqref{consistency_pretzel} and~\eqref{long_pretzel},
and we obtain 
$A_\mathcal{K}(\ell, m)=0$, where the A-polynomial for the $(-2,3,7)$
Pretzel knot is
\begin{multline}
  A_{\mathcal{K}}(\ell, m)
  =
  -1 +  \left(
    m^{16} - 2 \, m^{18} + m^{20}
  \right) \, \ell
  +
  \left( 2 \, m^{36} + m^{38} \right) \, \ell^2
  \\
  - \ell^4 \, \left(
    m^{72} + 2 \, m^{74}
  \right)
  -\ell^5 \, \left(
    m^{90} - 2 \, m^{92} + m^{94}
  \right)
  + m^{110} \, \ell^6
\end{multline}

\subsection{Once-Punctured Torus Bundles over the Circle}

One of benefits of the quantum invariant $Z_\gamma(\mathcal{M})$ is
that we can compute it explicitly
for a once-punctured torus bundle over $S^1$.
A once-punctured torus bundle
over $S^1$, which we denote  $\mathcal{M}(\varphi)$,
 is
described by
$F \times [0,1]/(x,0) \sim (\varphi(x),1)$ where
a monodromy matrix 
$\varphi \in SL(2,\mathbb{Z})$
is a homeomorphism from the punctured torus
$F=\mathbb{T}^2 \setminus \{ 0 \}$ to itself.
Thurston's hyperbolization theorem indicates that
$\mathcal{M}(\varphi)$ admits
a complete hyperbolic metric with a finite volume when $\varphi$ has 2
distinct real eigenvalues.
In this case, the monodromy matrix $\varphi$  can be written up to
conjugation as
\begin{equation}
  \varphi
  =
  L^{s_1} \, R^{t_1} \cdots
  L^{s_{n}}
  R^{t_{n}}
\end{equation}
where
\begin{align*}
  L 
  & =
  \begin{pmatrix}
    1 & 1 \\
    0 & 1
  \end{pmatrix}
  &
  R
  &=
  \begin{pmatrix}
    1 & 0 \\
    1 & 1
  \end{pmatrix}
\end{align*}
with
$n>0$,
and $s_j$ and $t_j$  are positive integers.
Note that 
the complement of the figure-eight knot studied in Sec.~\ref{sec:41}
corresponds to
$\varphi=L \, R =
\begin{pmatrix}
  2 & 1 \\
  1 & 1
\end{pmatrix}
$.

\begin{figure}[htbp]
  \centering
  \includegraphics[scale=0.88]{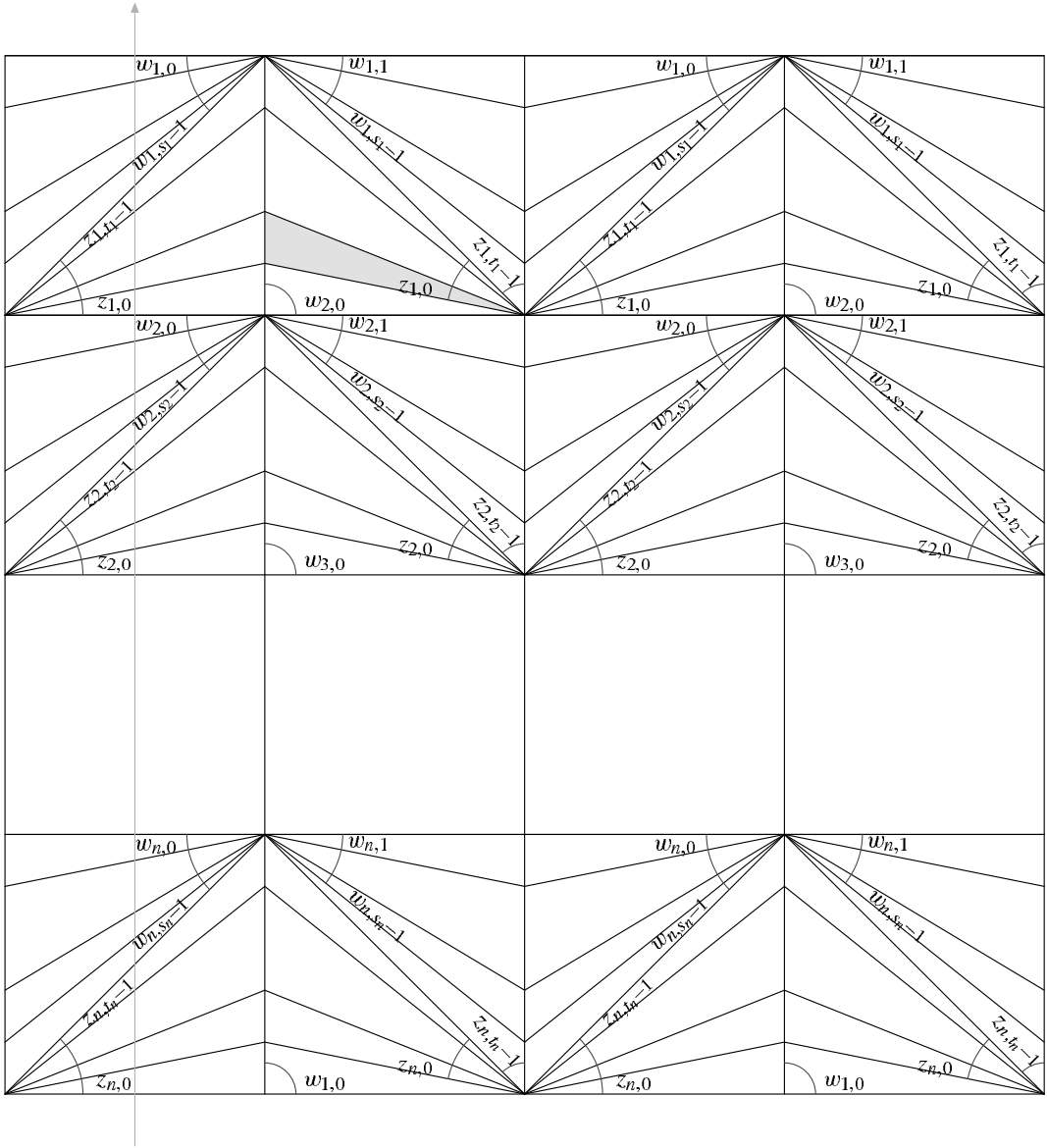}  
  \caption{Schematic developing map for $\mathcal{M}(\varphi)$.
    Top vertex of tetrahedron with modulus $z_{1,0}$ is filled gray.
    Gray straight line denotes meridian.
  }
  \label{fig:developLRgeneral}
\end{figure}

It is known that we can triangulate the manifold $M_\varphi$ with
$\sum_{k=1}^n (s_k+t_k)$ ideal tetrahedra~\cite{FloyHatc82a}, and we have the
quantum invariant as
\begin{multline}
  \label{Z_punctured_torus}
  Z_\gamma\left( \mathcal{M}_u(\varphi) \right)
  =
  \iiiint\limits_{\mathbb{R}}
  \mathrm{d} \boldsymbol{a} \,
  \mathrm{d} \boldsymbol{b} \,
  \mathrm{d} \boldsymbol{c} \,
  \mathrm{d} \boldsymbol{d} \,
  \delta_C(\boldsymbol{a}, \boldsymbol{b} ,
  \boldsymbol{c}, \boldsymbol{d}; u) \,
  \delta_G(\boldsymbol{a}, \boldsymbol{b} ,
  \boldsymbol{c}, \boldsymbol{d})
  \\
  \times
  \prod_{k=1}^n
  \left[
    \prod_{i=0}^{s_{k} -1}
    \left\langle d_{k,i} , c_{k, i+1} \middle| S^{-1} \middle|
      d_{k,i+1}, c_{k,i} \right\rangle \,
    \prod_{j=0}^{t_{k} -1}
    \left\langle b_{k,j} , a_{k,j+1} \middle| S \middle|
      b_{k,j+1}, a_{k,j} \right\rangle \,
  \right]
\end{multline}
where gluing condition $
  \delta_G(\boldsymbol{a}, \boldsymbol{b} ,
  \boldsymbol{c}, \boldsymbol{d})
$ means
\begin{align*}
  &
  \begin{cases}
    b_{k,0} = c_{k+1,0} 
    \\
    a_{k,0} = d_{k+1,0}
  \end{cases}
  \text{for $k=1,2, \dots,n-1$}
  & &
  \begin{cases}
    b_{n,0} = c_{1,0}
    \\
    a_{n,0} = d_{1,0}
  \end{cases}
  \\[2mm]
  & \begin{cases}
    b_{k,t_k} = c_{k, s_k}
    \\
    a_{k,t_k} = d_{k, s_k}
  \end{cases}
  \text{for $k=1,2,\dots, n$}
  & 
\end{align*}
When we set the modulus of each oriented tetrahedra as
\begin{align*}
  z_{k,j}
  & =
  \E^{a_{k,j} - a_{k,j+1}}
  &
  w_{k,j}
  & =
  \E^{c_{k,j} - c_{k,j+1}}
\end{align*}
the developing map is drawn schematically
as Fig.~\ref{fig:developLRgeneral}.
We can then read
a condition
$  \delta_C(\boldsymbol{a}, \boldsymbol{b} ,
  \boldsymbol{c}, \boldsymbol{d}; u)
$ for meridian  as
\begin{equation}
  u=
  \frac{1}{2} 
  \sum_{k=1}^n
  \left(
    c_{k,0} - c_{k,s_{k}} - a_{k,0} + a_{k,t_{k}}
  \right)
\end{equation}
With these conditions, we have the quantum invariant
$Z_\gamma \left( \mathcal{M}_u(\varphi) \right)$, and the
Neumann--Zagier potential function
can be given by taking a classical limit $\gamma\to 0$.
As far as we know,
the A-polynomial-type invariant for the once-punctured torus bundle over
$S^1$ is not studied, but we can  obtain
such polynomials from the
Neumann--Zagier function.

Below we give a few examples for concreteness.

\subsubsection{$L^2  \, R$}

\begin{figure}[htbp]
  \centering
    \includegraphics[scale=1.0]{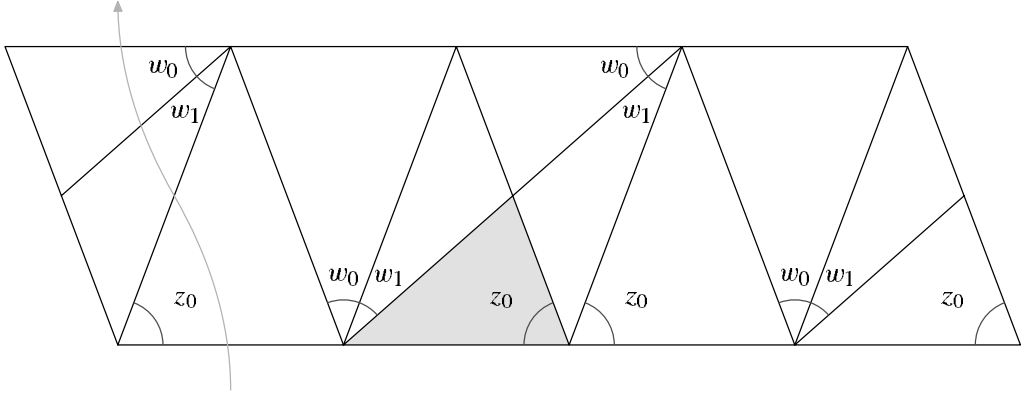}
  \caption{Developing map of $\mathcal{M}(L^2 \, R)$.
    Gray filled triangle is top vertex of the tetrahedron with modulus
    $z_0$.
    Gray curve is meridian.
  }
  \label{fig:developL2R}
\end{figure}

We set $\varphi = L^2 \, R=
\begin{pmatrix}
  3 & 2 \\
  1 & 1
\end{pmatrix}
$.
The hyperbolic Dehn surgery of this manifold is studied in
Ref.~\citen{BetlPrzyZuko86a}.
The invariant~\eqref{Z_punctured_torus} is rewritten as
\begin{multline}
  Z_\gamma
  \left( \mathcal{M}_u(L^2 \, R)
  \right)
  \\
  =
  \int\limits_\mathbb{R} \mathrm{d} \boldsymbol{p} \,
  \delta_C(\boldsymbol{p}; u) \,
  \left\langle
    p_1 , p_5 \middle| S^{-1} \middle| p_6 , p_3
  \right\rangle \,
  \left\langle
    p_6 , p_4 \middle| S^{-1} \middle| p_2 , p_5
  \right\rangle \,
  \left\langle
    p_3 , p_2 \middle| S \middle| p_4 , p_1
  \right\rangle   
\end{multline}
When we set the modulus of tetrahedra as
$w_0=\E^{p_3 - p_5}$,
$w_1 = \E^{p_5 - p_4}$, and
$z_0 = \E^{p_1 - p_2}$,
the developing map $\mathcal{M}(L^2 \, R)$
can be depicted as Fig.~\ref{fig:developL2R}.
Then  the meridian  is read as
\begin{equation}
  p_3 - p_4 - p_1 + p_2 = 2 \, u
\end{equation}

In the classical limit, we have
\begin{equation}
  Z_\gamma \left( \mathcal{M}_u(L^2 \, R) \right)
  \sim
  \iint \mathrm{d} x \, \mathrm{d} y \,
  \exp \left(
    \frac{1}{2 \, \I \, \gamma} \,
    V_{\mathcal{M}(L^2 R)}(\E^x, \E^y; \E^u)
  \right)
\end{equation}
where  the potential function is computed as
\begin{multline}
  V_{\mathcal{M}(L^2 R)}(x,y;m)
  =
  \Li \left( \frac{1}{m^2 \, x} \right)
  - \Li \left( \frac{1}{ m^2 \, x \, y^2} \right)
  - \Li \left( m^2 \, x^2 \, y^2 \right)
  \\
  - \log \left( m^2 \right)
  \log \left( m^2 \, x^2 \, y^4 \right) 
  -
  2 \left(
    \log \left( x \, y \right)
  \right)^2 + \frac{\pi^2}{6}
\end{multline}

The saddle point conditions~\eqref{hyperbolic_consistency},
which correspond to the hyperbolic consistency condition,
reduce to
\begin{equation}
  \label{consistent_L2R}
  \begin{gathered}
    \frac{
      \left( -1 + m^2 \, x \right) \,
      \left( -1 + m^2 \, x^2 \, y^2 \right)^2
    }{
      m^4 \, x^4 \, y^2 \,
      \left( -1 + m^2 \, x \, y^2 \right)
    } = 1 \\[2mm]
    \frac{
      \left(-1 + m^2 \, x^2 \,y^2 \right)
    }{
      m^2 \, x \,  \left( -1 + m^2 \, x \, y^2  \right)
    }
    =1
  \end{gathered}
\end{equation}
and,
as~\eqref{differential_longitude},
the longitude is defined by
\begin{equation}
  \label{longitude_L2R}
  \ell =
  - \frac{
    ( -1 + m^2 \, x) \,
    ( -1 + m^2 \, x^2 \, y^2 )
  }{
    m^4 \, x^2 \, y^2 \, ( -1 + m^2 \, x \, y^2)
  }
\end{equation}

In the case of complete case $m=1$, we 
find that
$
(x,y^2)=
\left(
  \frac{1 \pm \sqrt{7} \, \I}{4},
  \frac{-1 \pm \sqrt{7} \, \I}{2}
\right)
$ solves the consistency condition~\eqref{consistent_L2R}.
Among these, we can check numerically that
the largest value
of the imaginary part of the potential function at the saddle points
\begin{equation*}
  \Im V_{\mathcal{M}(L^2 R)}(x,y;m=1)
  =
  D \left( \frac{1}{ x} \right)
  - D \left( \frac{1}{  x \, y^2} \right)
  - D \left(  x^2 \, y^2 \right)
\end{equation*}
 coincides
with the hyperbolic volume  
$\Vol \left(\mathcal{M}(L^2 \, R) \right)
= 2.66674 \cdots$.

The A-polynomial is now given by eliminating $x$ and $y$ from the  set of
equations,~\eqref{consistent_L2R} and~\eqref{longitude_L2R},
and we obtain the algebraic curve
$A_{\mathcal{M}(L^2 R)}(\ell,m)=0$ whose  Newton polygon is given by
\begin{equation*}
  \begin{pmatrix}
    0 & -1 & 0
    \\
    0 & 2 & 1
    \\
    1 & 2 & 0
    \\
    0 & -1 & 0
  \end{pmatrix}
\end{equation*}

\subsubsection{$L \, R^3$}

\begin{figure}[htbp]
  \centering
  \includegraphics[scale=0.88]{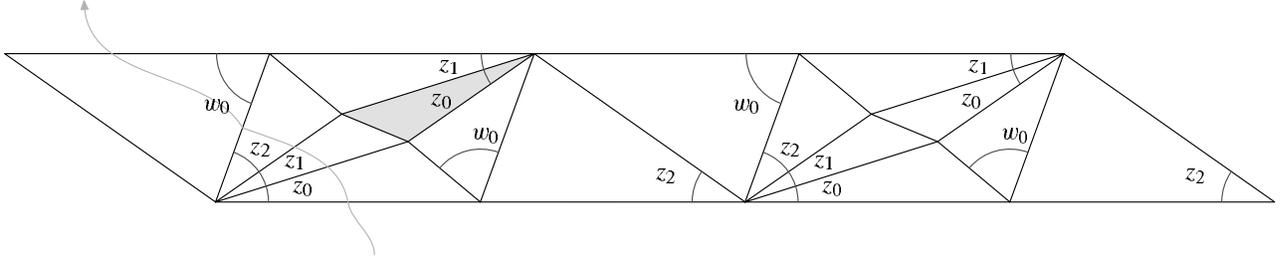}
  \caption{Developing map of $\mathcal{M}(L \, R^3)$.}
  \label{fig:developLR3}
\end{figure}

We take another example $\varphi= L\, R^3
=
\begin{pmatrix}
  4 & 1 \\
  3 & 1
\end{pmatrix}
$.
In this case, the invariant~\eqref{Z_punctured_torus} becomes
\begin{multline}
  Z_\gamma(\mathcal{M}_u(L \, R^3))
  =
  \int\limits_{\mathbb{R}} \mathrm{d} \boldsymbol{p} \,
  \delta_C(\boldsymbol{p};u) \,
  \left\langle p_1, p_4 \middle|  S \middle| p_3, p_2 \right\rangle \,
  \left\langle p_3, p_6 \middle|  S \middle| p_5, p_4 \right\rangle
  \\
  \times
  \left\langle p_5, p_8 \middle|  S \middle| p_7, p_6 \right\rangle \,
  \left\langle p_2, p_7 \middle|  S^{-1} \middle| p_8, p_1 \right\rangle 
\end{multline}
We set the modulus of each tetrahedron as
\begin{align*}
  z_0 & = \E^{p_2 - p_4}
  &
  z_1 & = \E^{p_4 - p_6}
  \\
  z_2 & = \E^{p_6 - p_8}
  &
  w_0 & = \E^{p_1 - p_7}
\end{align*}
The developing map is depicted in Fig.~\ref{fig:developLR3}, and
the meridian is read as $\frac{w_0}{z_0 \, z_1 \, z_2}$.
We thus have
a condition of $\delta_C(\boldsymbol{p} ; u)$ as
\begin{equation}
  p_1 - p_7 - p_2 + p_8 =  2 \, u
\end{equation}

In the classical limit, we obtain
\begin{equation}
  Z_\gamma(\mathcal{M}_u(L \, R^3))
  \sim
  \iiint \mathrm{d} x \, \mathrm{d} y \, \mathrm{d} z \,
  \exp
  \left(
    \frac{1}{2 \, \I \, \gamma} \,
    V_{\mathcal{M}(L R^3)}(\E^x , \E^y , \E^z ; \E^u)
  \right)
\end{equation}
where the potential function is computed as
\begin{multline}
  V_{\mathcal{M}(L R^3)}(x , y , z; m)
  =
  - \Li\left(\frac{1}{m^4 \, y}\right)
  +  \Li\left(\frac{1}{m^2 \, x  \, y \, z^2}\right)
  + \Li\left( m^2 \, x^2 \, z \right)
  \\
  + \Li \left( \frac{m^2 \, y^2 \, z}{x} \right)
  -\left( \log \left( m^2 \right) \right)^2
  + 2 \log\left(\frac{x}{m^2}\right) \log\left(\frac{x}{y}\right)
  \\
  + 2 \log y \log ( y \, z) + 2 \log z \log(x \, z \, m^2)
  -\frac{\pi^2}{3}
\end{multline}
In the case of complete case $m=1$, we can check numerically
that among algebraic solutions of~\eqref{hyperbolic_consistency}
the maximum of the imaginary part of
$V_{\mathcal{M}(L R^3)}(x,y,z;m=1)$,
\begin{equation*}
  \Im V_{\mathcal{M}(L R^3)}(x , y , z; 1)
  =
  - D \left(\frac{1}{ y}\right)
  +  D \left(\frac{1}{ x  \, y \, z^2}\right)
  + D\left(  x^2 \, z \right)
  + D \left( \frac{ y^2 \, z}{x} \right)
\end{equation*}
coincides with the hyperbolic volume
$\Vol(\mathcal{M}(L \, R^3))= 2.98912 \cdots$.
See Ref.~\citen{Franca04a} where algebraic solutions of consistency
equations are investigated in detail.

Correspondingly we obtain the A-polynomial for  $\mathcal{M}(L \, R^3)$
by eliminating $x$, $y$, and $z$ from a set of equations.
Explicitly the polynomial $A_{\mathcal{M}(L R^3)}(\ell,m)$ is given in
the form of Newton
polygon as follows;
\begin{equation*}
  \begin{pmatrix}
    0 & 0 & 0 & 1 & 0 & 0
    \\
    0 & 0 & 0 & -3 & 0 & 0
    \\
    -1 & -2 & -3 & -1 & -2 & 0
    \\
    0 & -3 & -2 & 2 & 3 & 0
    \\
    0 & 2 & 1 & 3 & 2 & 1
    \\
    0 & 0 & 3 & 0 & 0 & 0
    \\
    0 & 0  & -1 & 0 & 0 & 0
  \end{pmatrix}
\end{equation*}

\section{Conclusion and Discussion}

We have constructed  quantum invariant
$Z_\gamma(\mathcal{M})$
for cusped hyperbolic
3-manifolds $\mathcal{M}$
by assigning the Faddeev quantum dilogarithm function to oriented
ideal tetrahedra.
Once the triangulation of cusped 3-manifold $\mathcal{M}$ is given, it
is rather
straightforward to define the invariant $Z_\gamma(\mathcal{M})$ in an
integral form.
In the classical limit, the Faddeev integral reduces to the
dilogarithm function, which denotes the hyperbolic volume of ideal
tetrahedron.
Remarkable is that
the saddle point conditions coincide with the hyperbolic
consistency conditions around edges~\cite{Hikam00d}.
We have discussed as a variant of the volume conjecture
that the invariant
$Z_\gamma(\mathcal{M})$, which can be regarded as a generalization of
the Kashaev invariant (specific value of the colored Jones polynomial),
is dominated by the hyperbolic volume  in the classical limit
$\gamma\to 0$.

We have shown  that the invariant
can be defined for a one-parameter  deformation 
of  manifold $\mathcal{M}_u$ (not complete),
and that the Neumann--Zagier potential function can be given by taking
a classical limit $\gamma\to 0$.
Correspondingly the A-polynomial can be computed from the potential
function~\eqref{v_differential_limit}.
This  supports  the generalized volume conjecture.
We have demonstrated by taking examples
that our method recovers
previously known A-polynomial when $\mathcal{M}$ is a complement of
hyperbolic knots.
We have further  applied our method  in the case of the once-punctured torus bundle over
the circle.
It seems that the A-polynomial-type invariant is not known for this
3-manifold,
and  it will be interesting to study a relationship with the boundary slope.
To conclude,
our results indicate that the $S$-operator~\eqref{define_S_operator}
denotes the quantum Bloch
invariant for oriented ideal tetrahedron.

The A-polynomial may give an interesting insight for mathematical
physics.
For example, the Mahler measure of the A-polynomial is expected to
coincide with
the hyperbolic volume of knot.
On the other hand, the Mahler measure of the determinant of the
Kasteleyn matrix has appeared as the free energy
of the dimer problem.
From the viewpoint of our combinatorial
construction of the quantum invariant,
the quantum invariants based on
oriented ideal triangulation may be interpreted  as a matching of
in-states and out-states, and it might be interesting to investigate a
Kasteleyn matrix interpretation of the A-polynomial.

\section*{Acknowledgments}

We have used computer programs,
SnapPea~\cite{SnapPea99}, Knotscape~\cite{Knotscape99}, and
Snap~\cite{CouGooHodNeu00a}, in studying triangulation of manifolds.
We have also used \texttt{Mathematica} and Pari/GP.
Pictures of knots in this paper are drawn using KnotPlot~\cite{KnotPlot}.
This work is supported in part  by Grant-in-Aid for Young Scientists
from the Ministry of Education, Culture, Sports, Science and
Technology of Japan.


\appendix
\section{More Examples}

We shall study the quantum invariant for other cusped manifolds.
We give 
\begin{itemize}
\item 
  the quantum invariant $Z_\gamma(\mathcal{M})$ in terms of the
  $S$-operators,

\item 
  developing map when the number of the ideal tetrahedron is less than
  four,

\item
  the  condition $\delta_C(\boldsymbol{p};u)$, where $u$
  is a deformation parameter from the completeness $u=0$,

\item 
  the Neumann--Zagier  potential function
  $V_\mathcal{M}(\boldsymbol{x} ; m)$,
  which follows from $Z_\gamma(\mathcal{M})$
  by taking the classical limit $\gamma\to 0$,

\item a solution of the saddle point equations in the complete case
  $u=0$,
  which is checked
  that  the hyperbolic volume coincides with a maximal
  value of the imaginary part of the potential function at this saddle
  point,

\item the A-polynomial, which is given from the potential function
  $V_\mathcal{M}(\boldsymbol{x};m)$ by eliminating parameters
  $\boldsymbol{x}$.

\end{itemize}

In the first part~\ref{sec:first_appendix}, we collect results
for hyperbolic knots up to 7
crossings.
Though the A-polynomial is given in Ref.~\citen{CCGLS94a}, we give an
explicit form for self-completeness.
In the second part~\ref{sec:second_appendix}, we consider
\emph{simple} hyperbolic knots $K \,x_y$ from
Ref.~\citen{CallDeanWeek99a}.
Pictures of knot may be complicated, but the triangulation of the
complement is relatively simple in these cases.
We note that,
though
the number of the ideal tetrahedra of the
complement of knot $K \, x_y$ is $x$ in the case of the canonical
triangulation~\cite{SnapPea99},  our triangulations are different.

\subsection{Complement of Knots up to 7 Crossings}
\label{sec:first_appendix}
\subsubsection{
  $6_1$}

\begin{align*}
  & \mbox{
    \raisebox{-1.6cm}{
      \includegraphics[scale=.16]{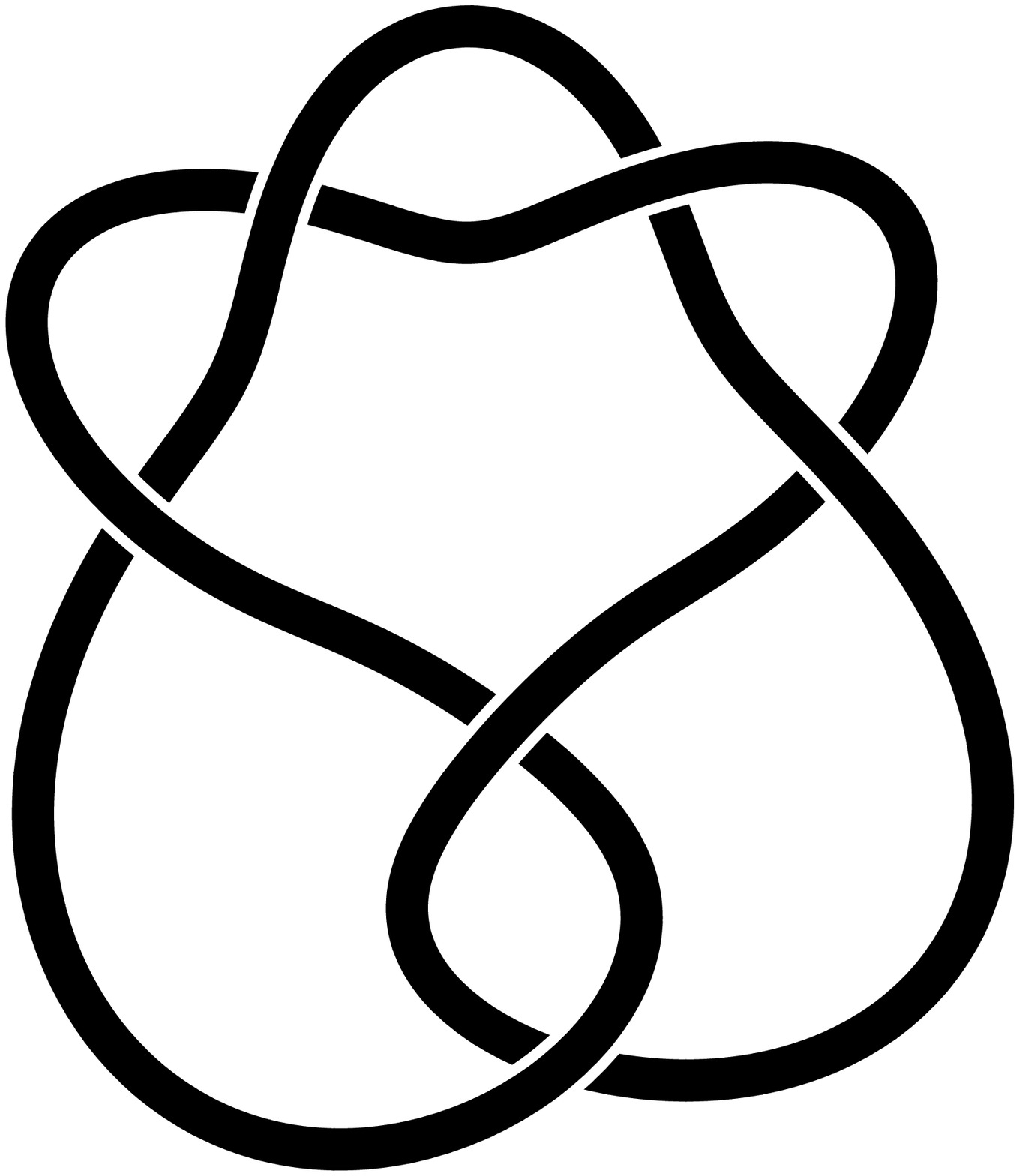}  
    }}
  &
  & \Vol(S^3 \setminus \mathcal{K})= 3.16396 \dots
\end{align*}

\begin{itemize}
\item Quantum invariant
\begin{multline*}
  Z_\gamma\left( \mathcal{M}_u \right)
  =
  \int\limits_\mathbb{R} \mathrm{d} \boldsymbol{p} \,
  \delta_C(\boldsymbol{p}; u ) \,
  \langle p_3 , p_1 | S^{-1} | p_2 , p_4 \rangle \,
  \langle p_6 , p_4 | S | p_5 , p_8 \rangle 
  \\
  \times
  \langle p_7 , p_5 | S^{-1} | p_6 , p_1 \rangle \,
  \langle p_8 , p_2 | S | p_7 , p_3 \rangle 
\end{multline*}

\item Developing map (Fig.~\ref{fig:develop61})

\begin{figure}[htbp]
  \centering
  \includegraphics[scale=1.]{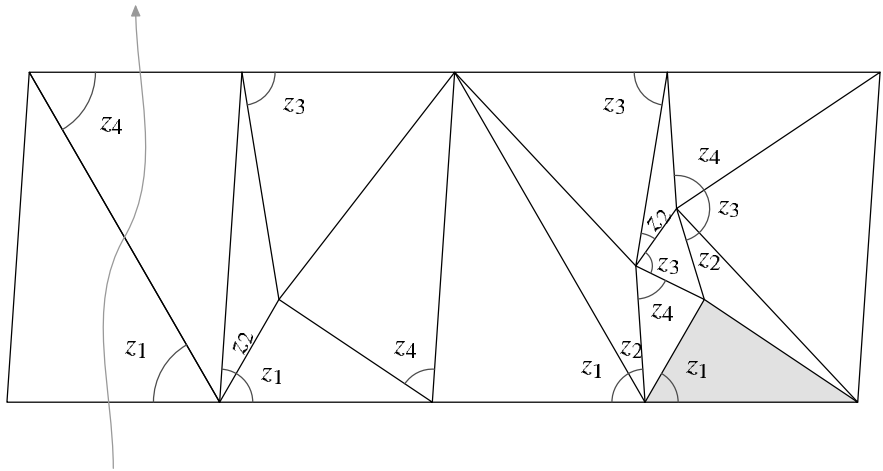}
  \caption{Developing map of the knot complement of $6_1$.
    Gray triangle is a horosphere of the top vertex of the tetrahedron
    with modulus $z_1$, where we have set
    $z_1  = \E^{p_4 - p_1}$ ,
    $z_2  = \E^{p_8 - p_4}$,
    $z_3  = \E^{p_1 - p_5}$, and
    $z_4  = \E^{p_3 - p_2}$.
  }
  \label{fig:develop61}
\end{figure}

\item Condition $\delta_C(\boldsymbol{p};u)$
\begin{equation*}
  p_4 - p_1 - p_3 + p_2
  =
  2 \, u
\end{equation*}

\item Potential function
\begin{multline*}
  V_{\mathcal{M}}(x,y,z ; m)
  =
  \Li\left(\frac{z \, m^2}{x \, y}\right) + \Li(y)
  - \Li\left(\frac{m^2}{y}\right)
  - \Li\left(\frac{y \, z}{m^2}\right)
  \\
  + \log \left( m^2 \right)  \, \log \left(\frac{y \, z^2}{x}\right)
  - \log z \, \log(x \, y) 
\end{multline*}

\item Hyperbolic volume
  \begin{equation*}
    \Vol(\mathcal{M})=
    \Im V_{\mathcal{M}}(x,y,z ; m=1)
    =
    D\left(\frac{z }{x \, y}\right) + D(y)
    - D\left(\frac{1}{y}\right)
    - D\left({y \, z}\right)
  \end{equation*}
  where
  \begin{align*}
    \begin{pmatrix}
      x\\y\\z
    \end{pmatrix}
    =
    \begin{pmatrix}
      -0.851808 + 0.911292 \, \I
      \\
      0.278726 + 0.48342  \, \I
      \\
     -1.50411 - 1.22685 \, \I
   \end{pmatrix}
 \end{align*}

\item   
 A-polynomial
  \begin{equation*}
    \begin{pmatrix}
      0 & 0 & -1 & 1 & 0 \\
      0 & 0 & 3 & -1 & 0 \\
      0 & 2 & 1 & 0 & 0 \\
      0 & -3 & -3 & 0 & 0 \\
      -1 & -3 & -6 & -3 & -1 \\
      0 & 0 & -3 & -3 & 0 \\
      0 & 0 & 1 & 2 & 0 \\
      0 & -1 & 3 & 0 & 0 \\
      0 & 1 & -1 & 0 & 0 
    \end{pmatrix}
  \end{equation*}
\end{itemize}

\subsubsection{
  $6_2$}
\begin{align*}
  & \mbox{
    \raisebox{-1.6cm}{
      \includegraphics[scale=.16]{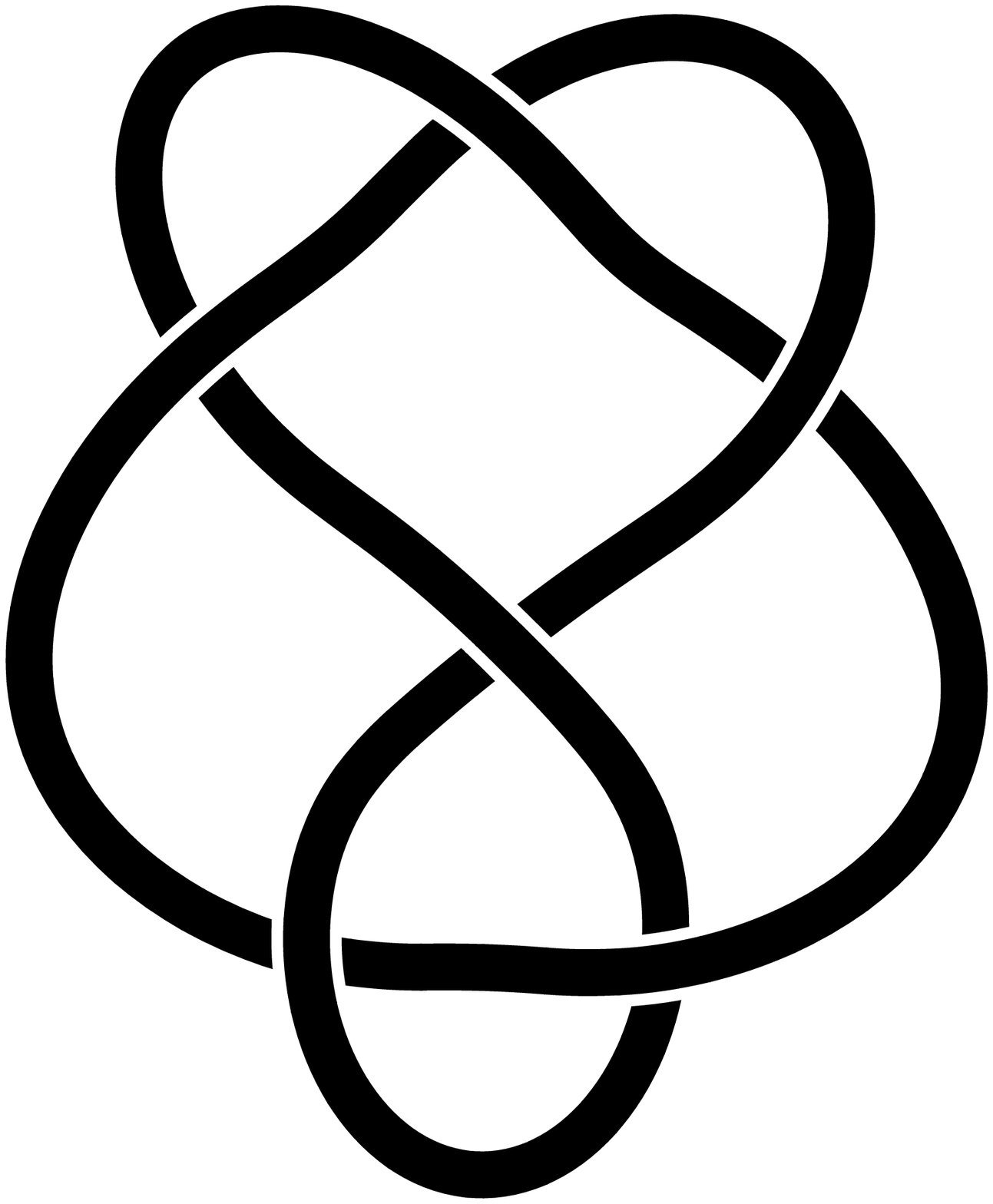}
    }}
  &
  & \Vol(S^3 \setminus \mathcal{K})= 4.40083 \dots
\end{align*}

\begin{itemize}
\item  Quantum invariant
\begin{multline*}
  Z_\gamma\left( \mathcal{M}_u \right)
  = \int\limits_\mathbb{R} 
  \mathrm{d} \boldsymbol{p} \,
  \delta_C(\boldsymbol{p}; u) \,
    \langle p_8 , p_6 | S^{-1} | p_1 , p_2 \rangle \, 
    \langle p_3 , p_7 | S^{-1} | p_8 , p_4 \rangle \,
    \\
    \times
    \langle p_4 , p_5 | S | p_3 , p_9 \rangle \, 
    \langle p_2 , p_{10} | S | p_6 , p_5 \rangle \, 
    \langle p_1 , p_9 | S^{-1} | p_{10} , p_7 \rangle 
  \end{multline*}

\item
Condition $\delta_C(\boldsymbol{p};u)$
  \begin{equation*}
    p_2 - p_6 - p_5+ p_{10}
    = -  2 \, u
  \end{equation*}

\item Potential function
  \begin{multline*}
    V_{\mathcal{M}}(w,x,y,z; m)
    =
    \Li\left(\frac{z}{m^2}\right)
    + \Li\left(\frac{m^4}{y \, w}\right)
    - \Li\left(\frac{y \, w}{m^2}\right)
    - \Li\left(\frac{m^2}{x \, y}\right)
    - \Li\left(\frac{y \, z}{m^2}\right)
    \\
    + \log x \, \log z
    - \log( y \, w) \, \log(y \, z)
    + \log \left( m^2 \right) \, \log\left(\frac{y^2 \, z}{w}\right)
  \end{multline*}

\item Hyperbolic volume
  \begin{equation*}
    \Im V_{\mathcal{M}}(w,x,y,z; 1)
    =
    D\left({z}\right)
    + D\left(\frac{1}{y \, w}\right)
    - D\left(y \, w\right)
    - D\left(\frac{1}{x \, y}\right)
    - D\left(y \, z\right)
  \end{equation*}
  with
  \begin{equation*}
    \begin{pmatrix}
      w \\ x\\ y \\ z\\
    \end{pmatrix}
    =
    \begin{pmatrix}
      -0.455697 + 1.20015  \, \I
      \\
      -0.964913 - 0.621896 \, \I
      \\
      -0.418784 - 0.219165 \, \I
      \\
      0.0904327 + 1.60288 \, \I
    \end{pmatrix}
  \end{equation*}

  \item A-polynomial
  \begin{equation*}
    \begin{pmatrix}
      0 & 1 & 0 & 0 & 0 & 0 \\
      0 & -2 & 1 & 0 & 0 & 0 \\
      -1 & 1 & -3 & 0 & 0 & 0 \\
      0 & 2 & 1 & 0 & 0 & 0 \\
      0 & -5 & 5 & 0 & 0 & 0 \\
      0 & -5 & 3 & -3 & 0 & 0 \\
      0 & 3 & -12 & 8 & 0 & 0 \\
      0 & 0 & -13 & 3 & 0 & 0 \\
      0 & 0 & 3 & -13 & 0 & 0 \\
      0 & 0 & 8 & -12 & 3 & 0 \\
      0 & 0 & -3 & 3 & -5 & 0 \\
      0 & 0 & 0 & 5 & -5 & 0 \\
      0 & 0 & 0 & 1 & 2 & 0 \\
      0 & 0 & 0 & -3 & 1 & -1 \\
      0 & 0 & 0 & 1 & -2 & 0 \\
      0 & 0 & 0 & 0 & 1 & 0 
    \end{pmatrix}
  \end{equation*}
\end{itemize}

\subsubsection{
  $6_3$}
\begin{align*}
  & \mbox{
    \raisebox{-1.6cm}{
      \includegraphics[scale=.16]{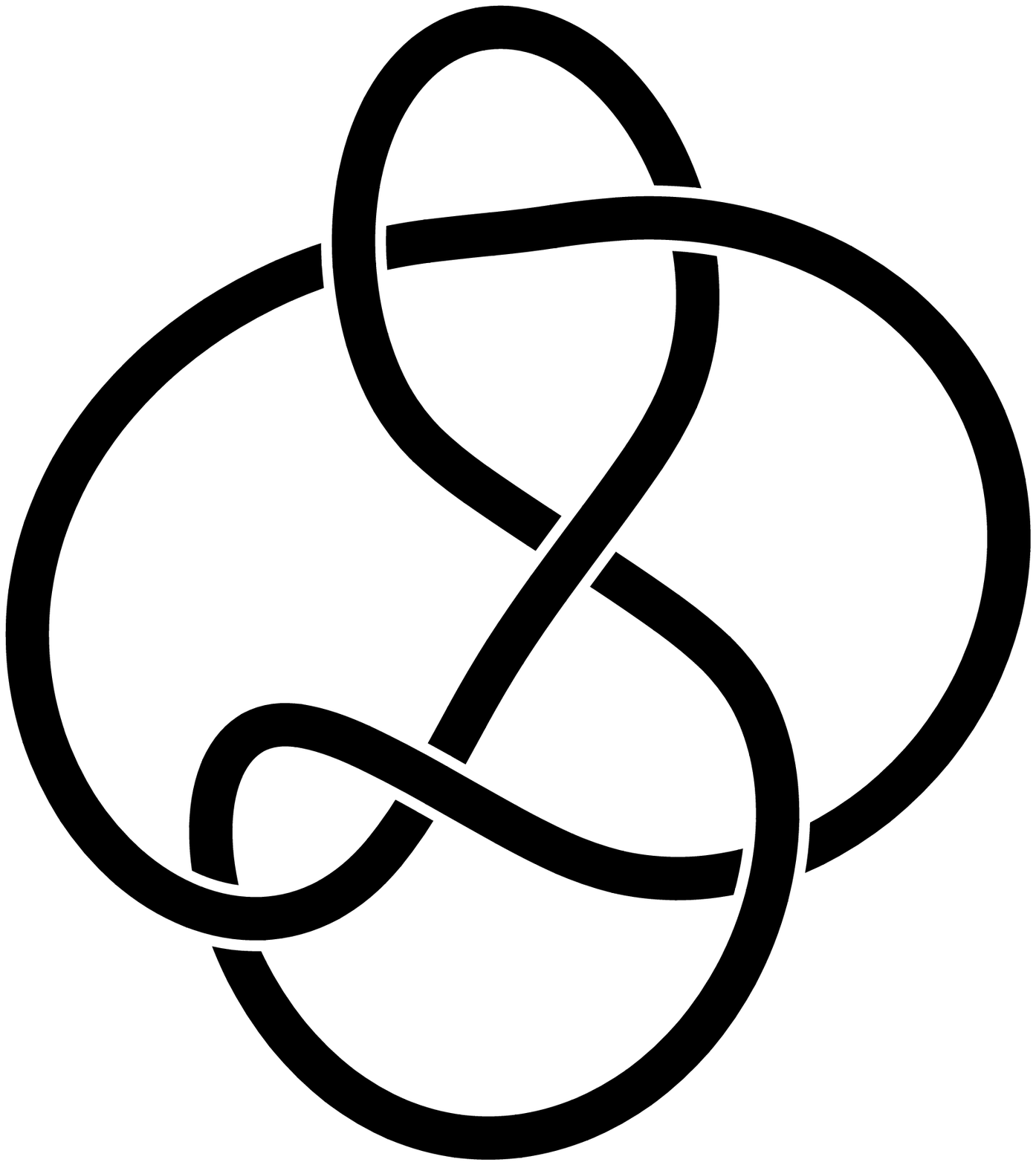}
    }}
  &
  & \Vol(S^3 \setminus \mathcal{K})= 5.69302 \dots
\end{align*}

\begin{itemize}
\item Quantum invariant
  \begin{multline}
    Z_\gamma \left( \mathcal{M}_u \right)
    =    \int\limits_\mathbb{R} 
    \mathrm{d} \boldsymbol{p} \,
    \delta_C(\boldsymbol{p}; u) \,
    \langle p_3 , p_2 | S | p_1 , p_9 \rangle \,
    \langle p_7 , p_9 | S^{-1} | p_5 , p_{10} \rangle \,
    \langle p_{10} , p_4 | S^{-1} | p_6 , p_8 \rangle 
    \\
    \times
    \langle p_8 , p_1 | S^{-1} | p_4 , p_{11} \rangle \,
    \langle p_{12} , p_{11} | S^{-1} | p_2 , p_3 \rangle \,
    \langle p_6 , p_5 | S | p_{12} , p_7 \rangle 
  \end{multline}

\item
Condition $\delta_C(\boldsymbol{p};u)$
  \begin{equation*}
    p_{10} - p_9 - p_7 + p_5
    = -  2 \, u
  \end{equation*}

\item Potential function
  \begin{multline*}
    V_{\mathcal{M}}(v, w, x,y,z; m)
    \\
    =
    \Li\left(\frac{m^2 \, y}{x}\right)
    + \Li\left(\frac{v}{w}\right)
    - \Li\left(\frac{m^2 \, w}{v}\right)
    -\Li\left(\frac{x \, z}{v}\right)
    - \Li\left(\frac{x^2 \, z}{v}\right)
    -\Li\left(\frac{v}{x \, y\, z}\right)
    \\
    +\frac{\pi^2}{3}
    + \log \left( m^2 \right) \, \log\left(\frac{y}{w}\right)
    + \log x \, \log\left(\frac{v}{y \, z^2}\right)
    + \log \left(\frac{w}{y \, z} \right)
    \log \left(\frac{z}{v}\right)
  \end{multline*}

\item Hyperbolic volume
  \begin{equation*}
    \Im V_{\mathcal{M}}(v, w, x,y,z; 1)
    =
    D\left(\frac{ y}{x}\right)
    + D\left(\frac{v}{w}\right)
    - D\left(\frac{ w}{v}\right)
    -D\left(\frac{x \, z}{v}\right)
    - D\left(\frac{x^2 \, z}{v}\right)
    -D\left(\frac{v}{x \, y\, z}\right)
  \end{equation*}
  with
  \begin{equation*}
    \begin{pmatrix}
      v \\ w\\ x\\ y \\ z
    \end{pmatrix}
    =
    \begin{pmatrix}
      0.0739495 + 0.558752 \, \I
      \\
      0.732786 + 0.381252 \, \I
      \\
      1.0
      \\
      0.108378 + 0.818891 \, \I
      \\
      0.415113 + 0.381252 \, \I
    \end{pmatrix}
  \end{equation*}

\item A-polynomial
  \begin{equation*}
    \begin{pmatrix}
      0 & 0 & 0 & 1 & 0 & 0 & 0 \\
      0 & 0 & 1 & -5 & 1 & 0 & 0 \\
      0 & 0 & -4 & 3 & -4 & 0 & 0 \\
      0 & 0 & 4 & 9 & 4 & 0 & 0 \\
      0 & 2 & 2 & -2 & 2 & 2 & 0 \\
      0 & -5 & -6 & -21 & -6 & -5 & 0 \\
      0 & 1 & 2 & 8 & 2 & 1 & 0 \\
      1 & 10 & 17 & 34 & 17 & 10 & 1 \\
      0 & 1 & 2 & 8 & 2 & 1 & 0 \\
      0 & -5 & -6 & -21 & -6 & -5 & 0 \\
      0 & 2 & 2 & -2 & 2 & 2 & 0 \\
      0 & 0 & 4 & 9 & 4 & 0 & 0 \\
      0 & 0 & -4 & 3 & -4 & 0 & 0 \\
      0 & 0 & 1 & -5 & 1 & 0 & 0 \\
      0 & 0 & 0 & 1 & 0 & 0 & 0 
    \end{pmatrix}
  \end{equation*}
\end{itemize}

\subsubsection{
  $7_2$}
\begin{align*}
  & \mbox{
    \raisebox{-1.6cm}{
      \includegraphics[scale=.16]{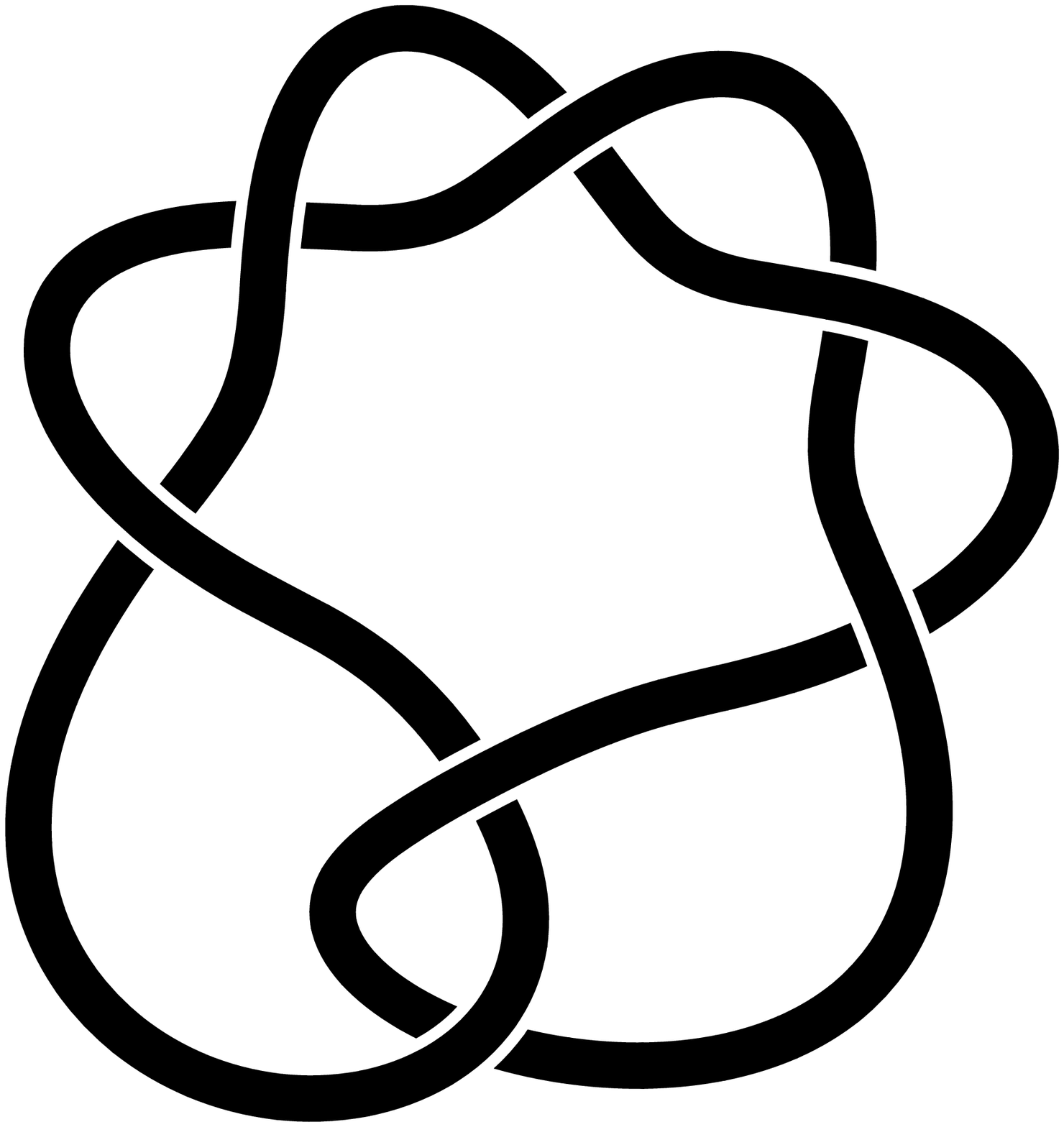}
    }}
  &
  & \Vol(S^3 \setminus \mathcal{K})= 3.33174 \dots
\end{align*}

\begin{itemize}
\item Quantum invariant
\begin{multline*}
  Z_\gamma \left( \mathcal{M}_u \right)
  = \int\limits_\mathbb{R} 
  \mathrm{d} \boldsymbol{p} \,
  \delta_C(\boldsymbol{p}; u) \,
  \langle p_1 , p_4 | S | p_2 , p_3 \rangle \,
  \langle p_5 , p_2 | S | p_8 , p_1 \rangle 
  \\
  \times
  \langle p_8 , p_7 | S^{-1} | p_7 , p_6 \rangle \,
  \langle p_6 , p_3 | S | p_4 , p_5 \rangle 
\end{multline*}

\item Developing map (Fig.~\ref{fig:develop72})

\begin{figure}[htbp]
  \centering
  \includegraphics[scale=1.04]{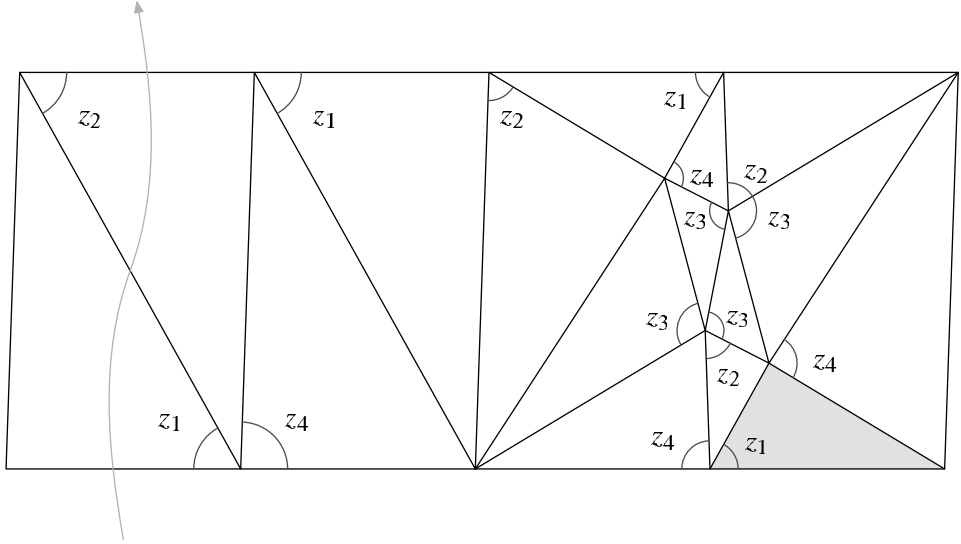}
  \caption{Developing map of the knot complement of $7_2$.
    The gray triangle corresponds to top vertex of the tetrahedron with
    modulus $z_1$, where we have set modulus as
    $z_1  = \E^{p_3 - p_4}$ ,
    $z_2  = \E^{p_1 - p_2}$,
    $z_3  = \E^{p_6 - p_7} $, and
    $z_4  = \E^{p_5 - p_3}$.
  }
  \label{fig:develop72}
\end{figure}

\item Condition $\delta_C(\boldsymbol{p}; u)$
\begin{equation*}
  p_3 - p_4 - p_1 + p_2
  = -  2\, u
\end{equation*}

\item Potential function
  \begin{multline*}
    V_{\mathcal{M}}(x,y,z ; m)
    =
    \Li\left(\frac{1}{x}\right)
    + \Li\left(\frac{m^2}{x}\right)
    + \Li \left( z \, m^2 \right)
    - \Li\left(\frac{y \, z}{x^2}\right)
    -\frac{\pi^2}{3}
    \\
    -    \left( \log(x/y) \right)^2
    + \log z \, \log \left(\frac{x^3 \, m^2}{y^2 \, z}\right)
  \end{multline*}

\item Hyperbolic volume
  \begin{equation*}
    \Im V_{\mathcal{M}}(x,y,z ; 1)
    =
    2 \, D\left(\frac{1}{x}\right)
    + D \left( z  \right)
    - D\left(\frac{y \, z}{x^2}\right)
  \end{equation*}
  with
  \begin{equation*}
    \begin{pmatrix}
      x \\ y\\ z
    \end{pmatrix}
    =
    \begin{pmatrix}
      0.941819 - 1.69128 \, \I  \\
      0.935538  + 0.903908  \, \I \\
      0.0581814 + 1.69128 \, \I
    \end{pmatrix}
  \end{equation*}

\item A-polynomial
  \begin{equation*}
    \begin{pmatrix}
      1 & -2 & 1 & 0 & 0 & 0 \\
      0 & 4 & -4 & 0 & 0 & 0 \\
      0 & 3 & 2 & -2 & 0 & 0 \\
      0 & 0 & 5 & 5 & 0 & 0 \\
      0 & 0 & 6 & 1 & 1 & 0 \\
      0 & 0 & 0 & -4 & -1 & 0 \\
      0 & -1 & -4 & 0 & 0 & 0 \\
      0 & 1 & 1 & 6 & 0 & 0 \\
      0 & 0 & 5 & 5 & 0 & 0 \\
      0 & 0 & -2 & 2 & 3 & 0 \\
      0 & 0 & 0 & -4 & 4 & 0 \\
      0 & 0 & 0 & 1 & -2 & 1 \\
    \end{pmatrix}
  \end{equation*}
\end{itemize}

\subsubsection{
  $ 7_3$}
\begin{align*}
  & \mbox{
    \raisebox{-1.6cm}{
      \includegraphics[scale=.16]{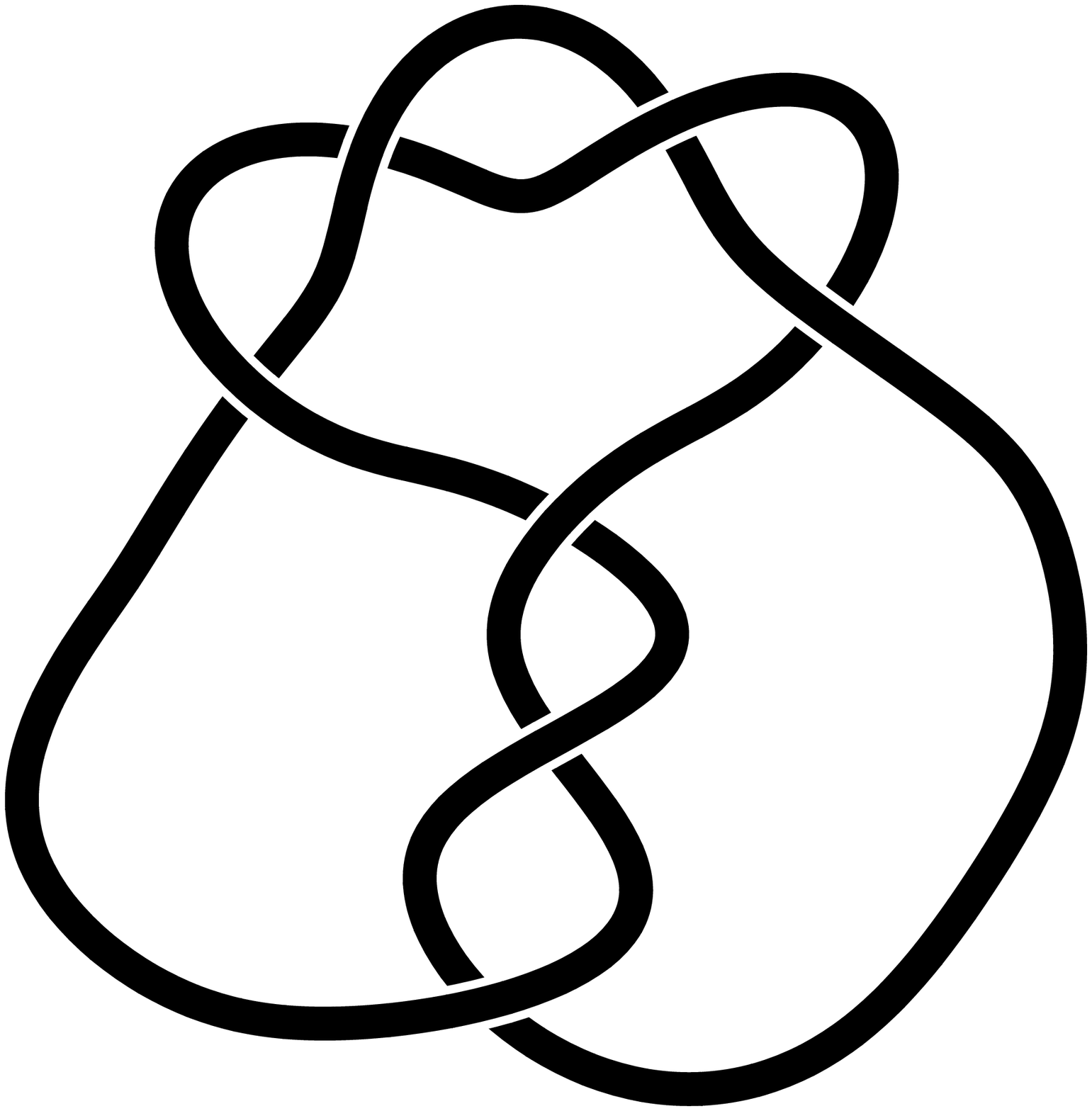}
    }}
  &
  & \Vol(S^3 \setminus \mathcal{K})=4.5921 \dots
\end{align*}

\begin{itemize}
\item Quantum invariant
  \begin{multline*}
    Z_\gamma \left( \mathcal{M}_u \right)
    = \int\limits_\mathbb{R} 
 \mathrm{d} \boldsymbol{p} \,
  \delta_C(\boldsymbol{p}; u) \,
    \langle p_1 , p_3 | S | p_4 , p_8 \rangle \,
    \langle p_2 , p_8 | S | p_5 , p_1 \rangle 
    \\
    \times
    \langle p_5 , p_9 | S^{-1} | p_6 , p_7 \rangle \,
    \langle p_4 , p_6 | S | p_3 , p_{10} \rangle \,
    \langle p_7 , p_{10} | S | p_2 , p_9 \rangle 
  \end{multline*}

\item
Condition $\delta_C(\boldsymbol{p};u)$
  \begin{equation*}
    p_{10} - p_6 - p_1 + p_8
    = -  2 \, u
  \end{equation*}

\item Potential function
  \begin{multline*}
    V_{\mathcal{M}}(w,x,y,z; m)
    =
    \Li \left(m^2 \, w\right) + \Li(m^4 \, w)
    + \Li\left(\frac{1}{m^2 \, x}\right)
    + \Li\left(\frac{z}{w}\right)
    - \Li\left(\frac{z}{y}\right) - \frac{\pi^2}{2}
    \\
    + \log w \, \log x + \log y \, \log z
    + \log \left( m^2 \right) \, \log (y \, w^2 \, z)
  \end{multline*}

\item Hyperbolic volume
  \begin{equation*}
    \Im V_{\mathcal{M}}(w,x,y,z; 1)
    =
    2 \, D \left( w\right)
    + D\left(\frac{1}{ x}\right)
    + D\left(\frac{z}{w}\right)
    - D\left(\frac{z}{y}\right)
  \end{equation*}
  with
  \begin{equation*}
    \begin{pmatrix}
      x \\ y \\ z \\ w
    \end{pmatrix}
    =
    \begin{pmatrix}
      0.645284 - 0.801205 \, \I
      \\
      -0.676708 +0.260961 \, \I
      \\
      -0.87287 + 1.51178 \, \I
      \\
       0.537981  + 1.04357  \, \I
    \end{pmatrix}
  \end{equation*}
\item A-polynomial
  \begin{equation*}
    \begin{pmatrix}
      -1 & 1 & 0 & 0 & 0 & 0 & 0 \\
      0 & -2 & 0 & 0 & 0 & 0 & 0 \\
      0 & 1 & 0 & 0 & 0 & 0 & 0 \\
      0 & 0 & 0 & 0 & 0 & 0 & 0 \\
      0 & -5 & -3 & 0 & 0 & 0 & 0 \\
      0 & -2 & 9 & 0 & 0 & 0 & 0 \\
      0 & 3 & -2 & 0 & 0 & 0 & 0 \\
      0 & -2 & -14 & 0 & 0 & 0 & 0 \\
      0 & 0 & -2 & 3 & 0 & 0 & 0 \\
      0 & 0 & 4 & -10 & 0 & 0 & 0 \\
      0 & 0 & -4 & 3 & 0 & 0 & 0 \\
      0 & 0 & -2 & 12 & 0 & 0 & 0 \\
      0 & 0 & -3 & -6 & -1 & 0 & 0 \\
      0 & 0 & 3 & 24 & 3 & 0 & 0 \\
      0 & 0 & -1 & -6 & -3 & 0 & 0 \\
      0 & 0 & 0 & 12 & -2 & 0 & 0 \\
      0 & 0 & 0 & 3 & -4 & 0 & 0 \\
      0 & 0 & 0 & -10 & 4 & 0 & 0 \\
      0 & 0 & 0 & 3 & -2 & 0 & 0 \\
      0 & 0 & 0 & 0 & -14 & -2 & 0 \\
      0 & 0 & 0 & 0 & -2 & 3 & 0 \\
      0 & 0 & 0 & 0 & 9 & -2 & 0 \\
      0 & 0 & 0 & 0 & -3 & -5 & 0 \\
      0 & 0 & 0 & 0 & 0 & 0 & 0 \\
      0 & 0 & 0 & 0 & 0 & 1 & 0 \\
      0 & 0 & 0 & 0 & 0 & -2 & 0 \\
      0 & 0 & 0 & 0 & 0 & 1 & -1 
    \end{pmatrix}
  \end{equation*}
\end{itemize}

\subsubsection{
  $ 7_4$}
\begin{align*}
  & \mbox{
    \raisebox{-1.6cm}{
      \includegraphics[scale=.16]{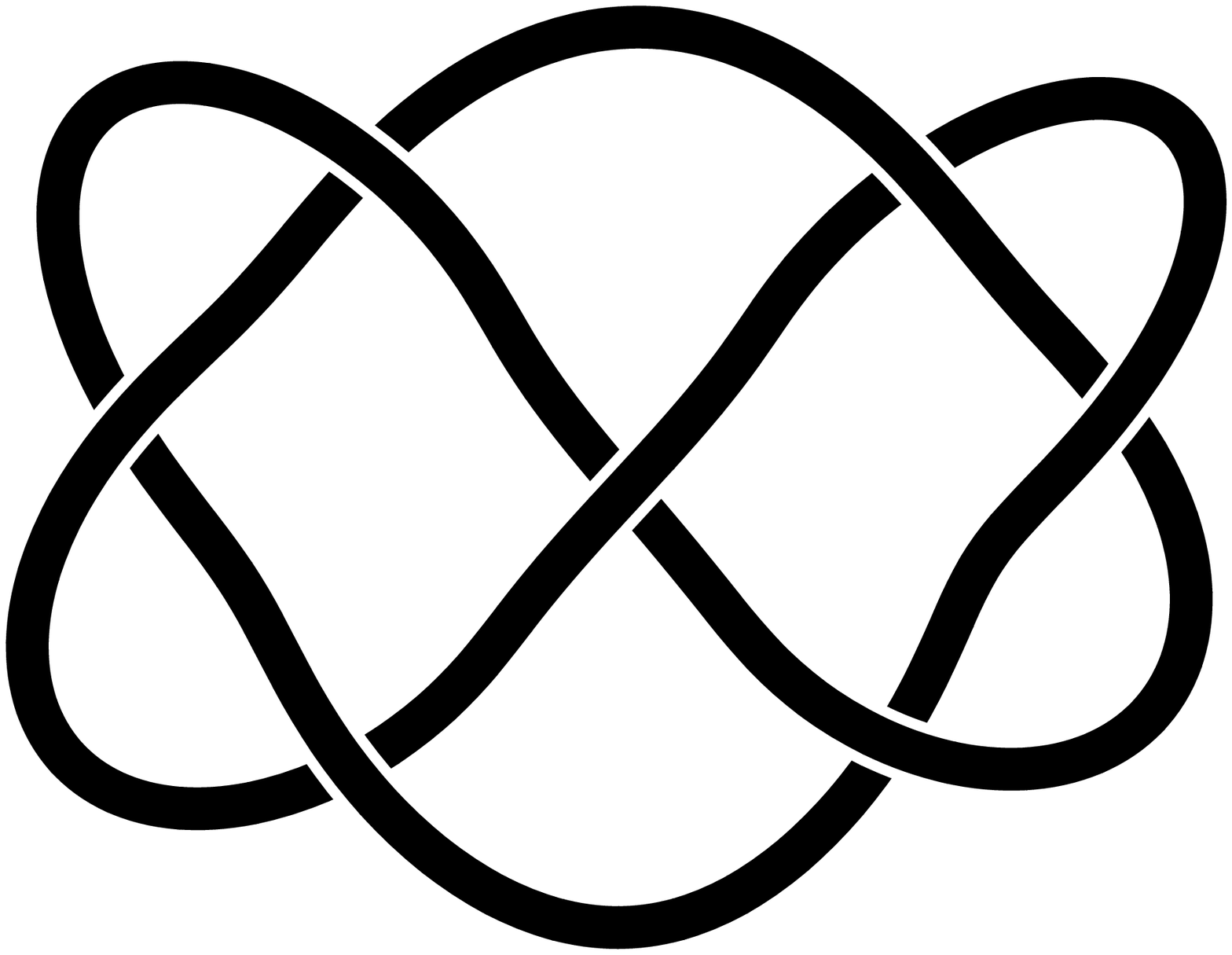}
    }}
  &
  & \Vol(S^3 \setminus \mathcal{K})=5.13794 \cdots
\end{align*}

\begin{itemize}
\item Quantum invariant
  \begin{multline*}
    Z_\gamma \left( \mathcal{M}_u \right)
    =
    \int\limits_{\mathbb{R}}
    \mathrm{d} \boldsymbol{p} \,
    \delta_C(\boldsymbol{p}; u ) \,
    \langle p_1, p_{10} | S| p_3 , p_4 \rangle \,
    \langle p_9, p_{2} | S| p_{12} , p_1 \rangle \,
    \langle p_4, p_{11} | S^{-1} | p_5 , p_6 \rangle
    \\
    \times
    \langle p_3, p_{5} | S^{-1} | p_2 , p_7 \rangle \,
    \langle p_{12}, p_{8} | S| p_{10} , p_9 \rangle \,
    \langle p_6, p_{7} | S| p_8 , p_{11} \rangle 
  \end{multline*}

\item
Condition $\delta_C(\boldsymbol{p};u)$
\begin{equation*}
  p_4-p_{10} + p_7 - p_5 - p_1 + p_2 = - 2 \, u 
\end{equation*}

\item Potential function
\begin{multline*}
  V_{\mathcal{M}}(v,w,x,y,z ; m)
  =
  \Li\left(\frac{w}{x}\right) +   \Li\left(\frac{v \, w}{y}\right)
  -   \Li\left(\frac{v \, x \, m^2}{y}\right)
  \\
  +  \Li\left(\frac{m^4 \, x}{w \, y}\right)
  - \Li\left(m^2 \, z\right)
  + \Li\left(\frac{z\, x}{y}\right)
  -\frac{\pi^2}{3}
  \\
  + \log(v) \log \left( \frac{w}{m^2 \, x \, z}\right)
  - \log x \log w - 2 \log \left(m^2\right) \log z
  +
  \left( \log \left(m^2/w\right) \right)^2
\end{multline*}

\item Hyperbolic volume
  \begin{multline*}
    \Im V_{\mathcal{M}}(v,w,x,y,z ; 1)
    \\
    =
    D\left(\frac{w}{x}\right)
    +   D\left(\frac{v \, w}{y}\right)
    -   D\left(\frac{v \, x }{y}\right)
    +  D\left(\frac{ x}{w \, y}\right)
    - D\left( z\right)
    + D\left(\frac{z\, x}{y}\right)
  \end{multline*}
  with
  \begin{equation*}
    \begin{pmatrix}
      v \\ w \\ x\\ y \\ z
    \end{pmatrix}
    =
    \begin{pmatrix}
      -1.10278 + 0.665457 \, \I
      \\
      -0.102785 + 0.665457 \, \I
      \\
      1.0
      \\
      -0.664742 - 0.401127 \, \I
      \\
      -0.226699 - 1.46771 \, \I
    \end{pmatrix}
  \end{equation*}

\item A-polynomial
\begin{equation*}
  \begin{pmatrix}
    0 & 0 & -1 & 3 & -3 & 1 \\
    0 & 0 & 3 & -10 & 7 & 0 \\
    0 & 0 & -3 & 3 & 4 & 0 \\
    0 & 0 & -2 & 21 & -6 & 0 \\
    0 & 1 & 10 & -3 & 1 & 0 \\
    0 & -2 & 6 & -17 & 3 & 0 \\
    0 & 3 & -17 & 6 & -2 & 0 \\
    0 & 1 & -3 & 10 & 1 & 0 \\
    0 & -6 & 21 & -2 & 0 & 0 \\
    0 & 4 & 3 & -3 & 0 & 0 \\
    0 & 7 & -10 & 3 & 0 & 0 \\
    1 & -3 & 3 & -1 & 0 & 0 \\
  \end{pmatrix}
\end{equation*}
\end{itemize}
\subsubsection{
  $7_5$}
\begin{align*}
  & \mbox{
    \raisebox{-1.6cm}{
      \includegraphics[scale=.16]{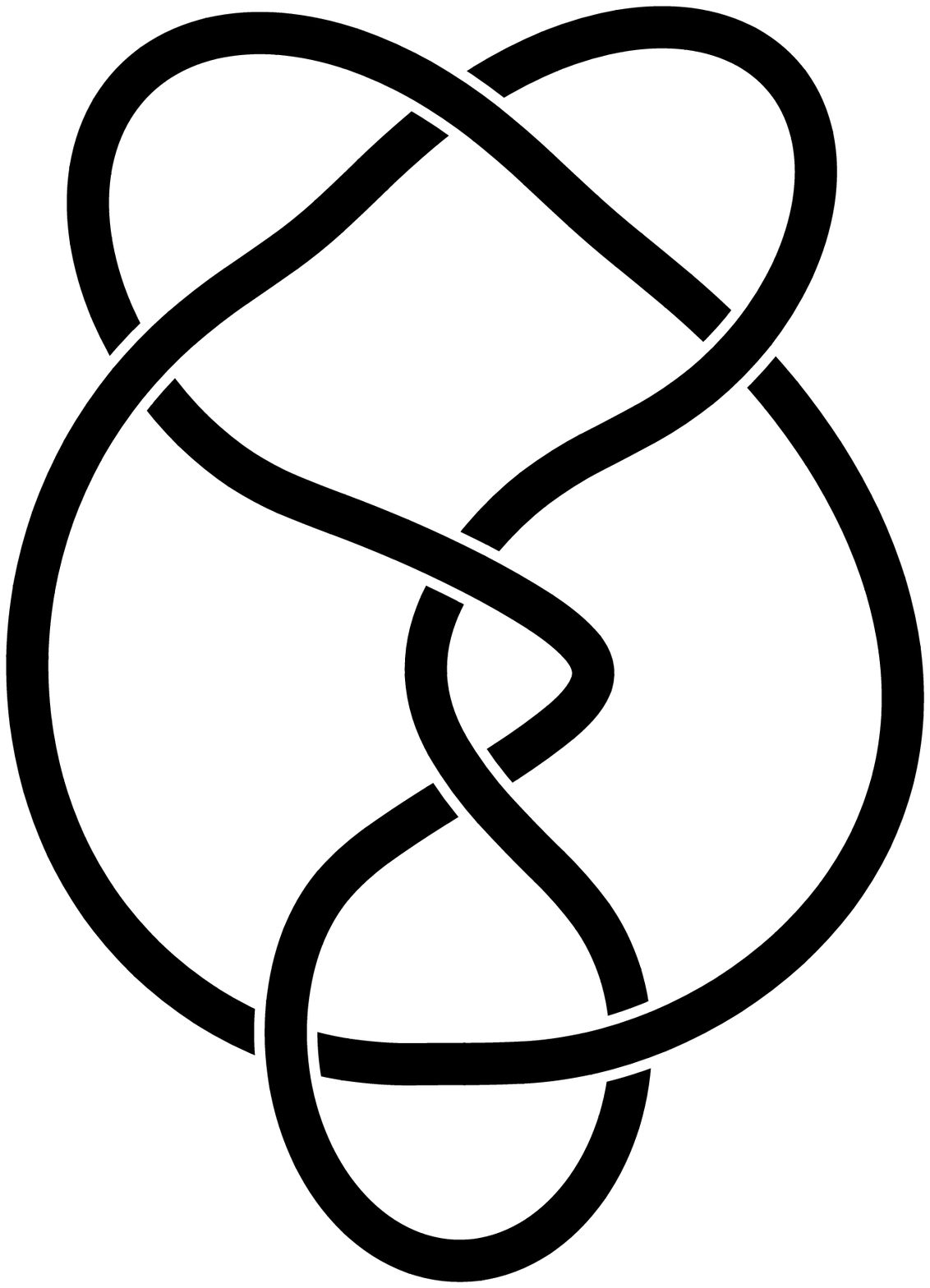}
    }}
  &
  & \Vol(S^3 \setminus \mathcal{K})= 6.443537 \dots
\end{align*}

\begin{itemize}
\item Quantum invariant
  \begin{multline*}
    Z_\gamma \left( \mathcal{M}_u \right)
    = \int\limits_\mathbb{R} 
 \mathrm{d} \boldsymbol{p} \,
  \delta_C(\boldsymbol{p}; u) \,
    \langle p_2 , p_{13} | S | p_1 , p_{12} \rangle \,
    \langle p_4 , p_{1} | S | p_3 , p_{10} \rangle
    \\
    \times 
    \langle p_6 , p_{3} | S | p_2 , p_{13} \rangle \,
    \langle p_5 , p_{11} | S^{-1} | p_4 , p_{7} \rangle \,
    \langle p_8 , p_{12} | S^{-1} | p_6 , p_{5} \rangle
    \\
    \times 
    \langle p_7 , p_{9} | S | p_8 , p_{14} \rangle \,
    \langle p_{10} , p_{14} | S^{-1} | p_9 , p_{11} \rangle 
  \end{multline*}

\item
Condition $\delta_C(\boldsymbol{p};u)$
\begin{equation*}
  p_3 + p_7 + p_{10} = p_1 + p_5 + p_{11} -  2 \, u
\end{equation*}

\item Potential function
  \begin{multline*}
    V_{\mathcal{M}}(u, v,w,x,y,z; m)
    =
    - \Li\left(\frac{u}{v}\right)
    + \Li\left(u \, m^2\right)
    + \Li\left(\frac{x}{m^2}\right)
    \\
    + \Li\left(\frac{m^2}{y}\right)
    + \Li\left(\frac{v}{m^2 \, w \, y}\right)
    - \Li\left(\frac{v}{z}\right)
    - \Li\left(\frac{w \, x}{z}\right)
    -\frac{\pi^2}{6}
    \\
    + \log \left(m^2\right) \, \log \left(\frac{u \, w^2}{v}\right)
    + \log x \, \log y + \log u \, \log z
  \end{multline*}

\item Hyperbolic volume
  \begin{multline*}
    \Im V_{\mathcal{M}}(u, v,w,x,y,z; m=1)
    =
    - D\left(\frac{u}{v}\right)
    + D\left(u \right)
    + D\left(x\right)
    \\
    + D\left(\frac{1}{y}\right)
    + D\left(\frac{v}{ w \, y}\right)
    - D\left(\frac{v}{z}\right)
    - D\left(\frac{w \, x}{z}\right)
  \end{multline*}
  with
  \begin{equation*}
    \begin{pmatrix}
      u \\ v \\ w \\ x \\ y \\ z
    \end{pmatrix}
    =
    \begin{pmatrix}
      0.38762  +   1.0287 \, \I
      \\
      -0.572726 + 0.717749 \, \I
      \\
      -0.259819 + 0.832925  \, \I
      \\
      0.18596 + 0.689115 \, \I
      \\
      0.365014  - 1.35264 \, \I
      \\
      -0.679246 - 0.851242 \, \I
    \end{pmatrix}
  \end{equation*}

\item A-polynomial
  \begin{equation*}
    \begin{pmatrix}
      -1 & 1 & 0 & 0 & 0 & 0 & 0 & 0 & 0 \\
      0 & -4 & 0 & 0 & 0 & 0 & 0 & 0 & 0 \\
      0 & 5 & -1 & 0 & 0 & 0 & 0 & 0 & 0 \\
      0 & 2 & 6 & 0 & 0 & 0 & 0 & 0 & 0 \\
      0 & -13 & -17 & 0 & 0 & 0 & 0 & 0 & 0 \\
      0 & -3 & 10 & -2 & 0 & 0 & 0 & 0 & 0 \\
      0 & 7 & 35 & 12 & 0 & 0 & 0 & 0 & 0 \\
      0 & -3 & -32 & -23 & 0 & 0 & 0 & 0 & 0 \\
      0 & 0 & -56 & -6 & -1 & 0 & 0 & 0 & 0 \\
      0 & 0 & 24 & 48 & 6 & 0 & 0 & 0 & 0 \\
      0 & 0 & 28 & 15 & -11 & 0 & 0 & 0 & 0 \\
      0 & 0 & -22 & -82 & 4 & 0 & 0 & 0 & 0 \\
      0 & 0 & -14 & -28 & 4 & 0 & 0 & 0 & 0 \\
      0 & 0 & 14 & 47 & 4 & -1 & 0 & 0 & 0 \\
      0 & 0 & -3 & -13 & -11 & 7 & 0 & 0 & 0 \\
      0 & 0 & 0 & -46 & 12 & -12 & 0 & 0 & 0 \\
      0 & 0 & 0 & 15 & -16 & -2 & 0 & 0 & 0 \\
      0 & 0 & 0 & 15 & -52 & 15 & 0 & 0 & 0 \\
      0 & 0 & 0 & -2 & -16 & 15 & 0 & 0 & 0 \\
      0 & 0 & 0 & -12 & 12 & -46 & 0 & 0 & 0 \\
      0 & 0 & 0 & 7 & -11 & -13 & -3 & 0 & 0 \\
      0 & 0 & 0 & -1 & 4 & 47 & 14 & 0 & 0 \\
      0 & 0 & 0 & 0 & 4 & -28 & -14 & 0 & 0 \\
      0 & 0 & 0 & 0 & 4 & -82 & -22 & 0 & 0 \\
      0 & 0 & 0 & 0 & -11 & 15 & 28 & 0 & 0 \\
      0 & 0 & 0 & 0 & 6 & 48 & 24 & 0 & 0 \\
      0 & 0 & 0 & 0 & -1 & -6 & -56 & 0 & 0 \\
      0 & 0 & 0 & 0 & 0 & -23 & -32 & -3 & 0 \\
      0 & 0 & 0 & 0 & 0 & 12 & 35 & 7 & 0 \\
      0 & 0 & 0 & 0 & 0 & -2 & 10 & -3 & 0 \\
      0 & 0 & 0 & 0 & 0 & 0 & -17 & -13 & 0 \\
      0 & 0 & 0 & 0 & 0 & 0 & 6 & 2 & 0 \\
      0 & 0 & 0 & 0 & 0 & 0 & -1 & 5 & 0 \\
      0 & 0 & 0 & 0 & 0 & 0 & 0 & -4 & 0 \\
      0 & 0 & 0 & 0 & 0 & 0 & 0 & 1 & -1 \\
    \end{pmatrix}
  \end{equation*}
\end{itemize}

\subsubsection{
  $7_6$}
\begin{align*}
  & \mbox{
    \raisebox{-1.6cm}{
      \includegraphics[scale=.16]{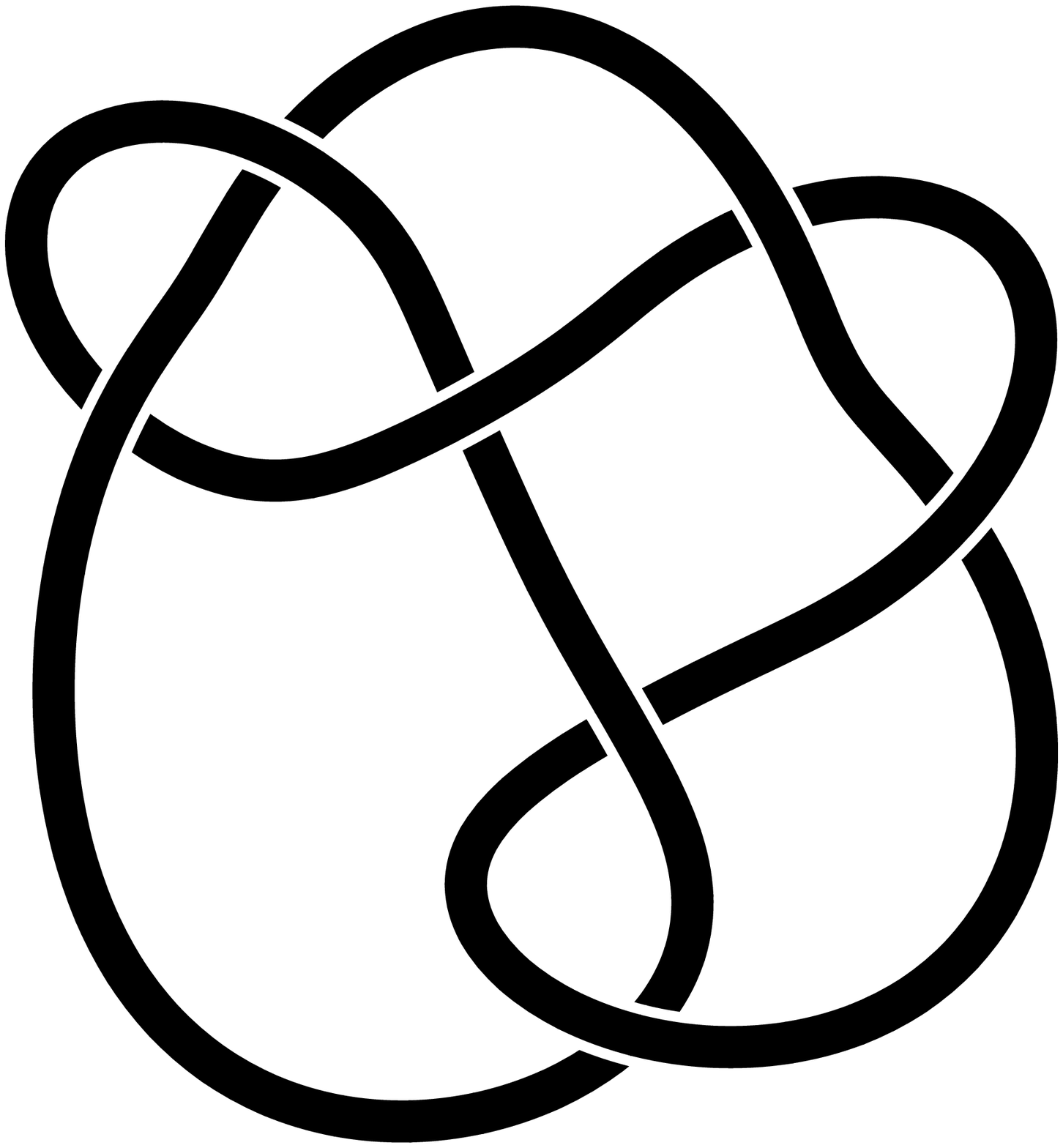}
    }}
  &
  &  \Vol(S^3 \setminus \mathcal{K})=7.08493 \cdots
\end{align*}

\begin{itemize}
\item Quantum invariant
  \begin{multline*}
    Z_\gamma \left( \mathcal{M}_u \right)
    = \int\limits_\mathbb{R} 
    \mathrm{d} \boldsymbol{p} \,
    \delta_C(\boldsymbol{p}; u) \,
    \langle p_1 , p_{14} | S^{-1} | p_3 , p_{13} \rangle \,
    \langle p_4 , p_{13} | S^{-1} | p_1 , p_{2} \rangle \,
    \langle p_9 , p_{5} | S | p_4 , p_{14} \rangle
    \\
    \times
    \langle p_3 , p_{2} | S | p_6 , p_{15} \rangle \,
    \langle p_7 , p_{6} | S | p_5 , p_{16} \rangle \,
    \langle p_{11} , p_{12} | S^{-1} | p_7 , p_{8} \rangle
    \\
    \times
    \langle p_8 , p_{15} | S^{-1} | p_6 , p_{10} \rangle \,
    \langle p_{10} , p_{16} | S^{-1} | p_{11} , p_{12} \rangle
  \end{multline*}

\item
  Condition $\delta_C(\boldsymbol{p};u)$
  \begin{equation*}
    p_6+p_{10} = p_5 + p_8 -  2 \, u
  \end{equation*}

\item Potential function
  \begin{multline*}
    V_{\mathcal{M}}(t,u, v, w,x,y,z;m)
    =
    -\Li(t) - \Li\left(u \, m^2 \right)
    + \Li \left( \frac{v}{x} \right)
    - \Li\left( \frac{m^2}{x} \right)
    \\
    + \Li\left(\frac{t}{y}\right)
    +\Li\left(\frac{u}{w \, y}\right)
    -  \Li\left(\frac{1}{m^2 \, z}\right)
    -\Li\left(\frac{v \, w}{z}\right)
    + \frac{\pi^2}{3}
    \\
    - ( \log \left( m^2 \right) )^2 + 2 \log \left(m^2\right) \log\left(\frac{x}{y \, w} \right)
    + \log(t/v) \log(x/y) - \log u \log z
  \end{multline*}

\item Hyperbolic volume
  \begin{multline*}
    \Im V_{\mathcal{M}}(t,u, v, w,x,y,z;1)
    =
    -D(t) - D\left(u  \right)
    + D \left( \frac{v}{x} \right)
    - D\left( \frac{1}{x} \right)
    \\
    + D\left(\frac{t}{y}\right)
    + D\left(\frac{u}{w \, y}\right)
    -  D\left(\frac{1}{ z}\right)
    - D\left(\frac{v \, w}{z}\right)
  \end{multline*}
  with
  \begin{equation*}
    \begin{pmatrix}
      t \\ u \\ v\\ w \\ x \\ y \\ z
    \end{pmatrix}
    =
    \begin{pmatrix}
       0.558614 - 1.43795 \, \I
       \\
       -0.0892864 - 0.842785 \, \I
       \\
       -0.280101 +1.13004 \, \I
       \\
       0.450985  - 0.808297 \, \I
       \\
       0.234736 +   0.604244 \, \I
       \\
       -0.20665 - 0.833705 \, \I
       \\
       -0.12431 + 1.17337 \, \I
     \end{pmatrix}
   \end{equation*}

\item A-polynomial
  \begin{equation*}
    \begin{pmatrix}
      0 & 0 & 1 & -1 & 0 & 0 & 0 & 0 & 0 & 0
      \\
      0 & 0 & -6 & 7 & 0 & 0 & 0 & 0 & 0 & 0
      \\
      0 & -2 & 11 & -16 & 1 & 0 & 0 & 0 & 0 & 0
      \\
      0 & 6 & 2 & 1 & -9 & 0 & 0 & 0 & 0 & 0
      \\
      1 & -5 & -16 & 34 & 32 & 0 & 0 & 0 & 0 & 0
      \\
      0 & -5 & -7 & 10 & -30 & 2 & 0 & 0 & 0 & 0
      \\
      0 & 16 & 9 & -80 & -68 & -16 & 0 & 0 & 0 & 0
      \\
      0 & 5 & 8 & -9 & 98 & 41 & 0 & 0 & 0 & 0
      \\
      0 & -9 & 42 & 62 & 164 & -18 & 1 & 0 & 0 & 0
      \\
      0 & 3 & 11 & 10 & -212 & -78 & 7 & 0 & 0 & 0
      \\
      0 & 0 & -37 & 34 & -266 & 52 & 19 & 0 & 0 & 0
      \\
      0 & 0 & 8 & 83 & 196 & 158 & -29 & 0 & 0 & 0
      \\
      0 & 0 & 23 & -44 & 377 & -85 & 10 & 0 & 0 & 0
      \\
      0 & 0 & -16 & -48 & 24 & -237 & 47 & 3 & 0 & 0
      \\
      0 & 0 & 3 & 47 & -237 & 24 & -48 & -16 & 0 & 0
      \\
      0 & 0 & 0 & 10 & -85 & 377 & -44 & 23 & 0 & 0
      \\
      0 & 0 & 0 & -29 & 158 & 196 & 83 & 8 & 0 & 0
      \\
      0 & 0 & 0 & 19 & 52 & -266 & 34 & -37 & 0 & 0
      \\
      0 & 0 & 0 & 7 & -78 & -212 & 10 & 11& 3 & 0
      \\
      0 & 0 & 0 & 1 & -18 & 164 & 62 & 42 & -9 & 0
      \\
      0 & 0 & 0 & 0 & 41 & 98 & -9 & 8 & 5 & 0
      \\
      0 & 0 & 0 & 0 & -16 & -68 & -80 & 9 & 16 & 0
      \\
      0 & 0 & 0 & 0 & 2 & -30 & 10 & -7 & -5 & 0
      \\
      0 & 0 & 0 & 0 & 0 & 32 & 34 & -16 & -5 & 1
      \\
      0 & 0 & 0 & 0 & 0 & -9 & 1 & 2 & 6 & 0
      \\
      0 & 0 & 0 & 0 & 0 & 1 & -16 & 11 & -2 & 0
      \\
      0 & 0 & 0 & 0 & 0 & 0 & 7 & -6 & 0 & 0
      \\
      0 & 0 & 0 & 0 & 0 & 0 & -1 & 1 & 0 & 0
    \end{pmatrix}
  \end{equation*}
\end{itemize}

\subsubsection{
  $7_7$}

\begin{align*}
  & \mbox{
    \raisebox{-1.6cm}{
      \includegraphics[scale=.16]{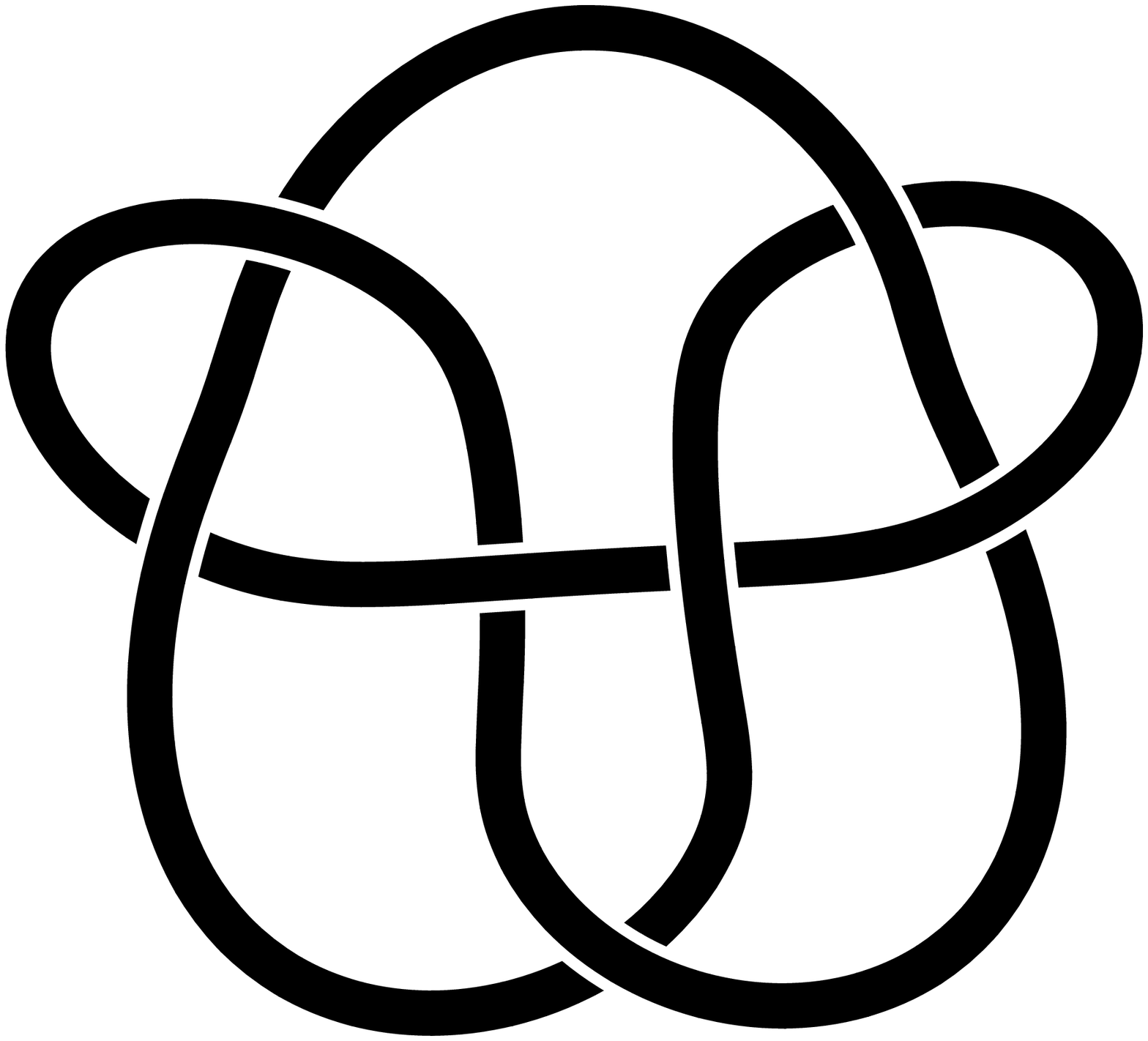}
    }}
  &
  & \Vol(S^3 \setminus \mathcal{K})= 7.64338 \cdots
\end{align*}

\begin{itemize}
\item Quantum invariant
\begin{multline*}
  Z_\gamma \left( \mathcal{M}_u \right)
  = \int\limits_\mathbb{R} 
  \mathrm{d} \boldsymbol{p} \,
  \delta_C(\boldsymbol{p}; u) \,
  \langle p_1 , p_{3} | S | p_{15} , p_{16} \rangle \,
  \langle p_2 , p_{4} | S | p_1 , p_{14} \rangle \,
  \langle p_5 , p_{13} | S^{-1} | p_8 , p_{4} \rangle
  \\
  \times
  \langle p_6 , p_{15} | S^{-1} | p_3 , p_{2} \rangle \,
  \langle p_7 , p_{16} | S^{-1} | p_5 , p_{6} \rangle \,
  \langle p_{10} , p_{14} | S^{-1} | p_7 , p_{9} \rangle
  \\
  \times
  \langle p_8 , p_{9} | S | p_{11} , p_{12} \rangle \,
  \langle p_{11} , p_{12} | S | p_{10} , p_{13} \rangle
\end{multline*}

\item
Condition $\delta_C(\boldsymbol{p};u)$
\begin{equation*}
  p_6+p_{15} = p_2 + p_3 -  2\, u
\end{equation*}

\item Potential function
\begin{multline*}
  V_{\mathcal{M}}(t,u, v,w,x,y,z;m)
  =
  \Li\left(v  \, m^2 \right) + \Li\left(t \, m^2 \, w\right) -
  \Li\left(\frac{t \, m^2 }{x}\right) 
  -  \Li\left(\frac{1}{m^2 \, w \, x}\right) 
  \\
  - \Li(v \, x) + \Li(u \, x) + \Li\left(\frac{1}{m^4 \, z}\right)
  - \Li\left(\frac{u}{m^2 \, z}\right)
  \\
  + \log \left(m^2 \right) \log \left( \frac{v \, u}{m^2 \, w \, y}\right)
  + \log( t \, u) \log\left(\frac{w \, x}{y} \right)
  + \log v \log z
\end{multline*}

\item Hyperbolic volume
  \begin{multline*}
    \Im V_{\mathcal{M}}(t,u, v,w,x,y,z;1)
    =
    D\left(v   \right)
    + D \left(t \,  w\right) -
    D\left(\frac{t  }{x}\right) 
    \\
    -  D\left(\frac{1}{ w \, x}\right) 
    - D(v \, x)
    + D(u \, x) + D\left(\frac{1}{ z}\right)
    - D \left(\frac{u}{ z}\right)
  \end{multline*}
  with
  \begin{equation*}
    \begin{pmatrix}
      t \\ u \\ v \\ w\\ x\\ y\\ z
    \end{pmatrix}
    =
    \begin{pmatrix}
      -0.899232 + 0.400532 \, \I
      \\
      -0.927958 - 0.413327 \, \I
      \\
      0.0287264  + 0.813859  \, \I
      \\
      -0.351808 - 0.720342  \, \I
      \\
      -0.927958 - 0.413327 \, \I
      \\
       -0.927958 - 0.413327  \, \I
       \\
       0.0433154 - 1.22719 \, \I
    \end{pmatrix}
  \end{equation*}

\item A-polynomial
  \begin{equation*}
    \begin{pmatrix}
      0 & 0 & 0 & 1 & -1 & 0 & 0 & 0  \\
      0 & 0 & 0 & -7 & 8 & 0 & 0 & 0  \\
      0 & 0 & 0 & 15 & -19 & 1 & 0 & 0  \\
      0 & 0 & 3 & 1 & -4 & -6 & 0 & 0  \\
      0 & 0 & -18 & -30 & 59 & 11 & 0 & 0  \\
      0 & 0 & 23 & -1 & 0 & -1 & 0 & 0  \\
      0 & 3 & 27 & 41 & -123 & -8 & -2 & 0  \\
      0 & -11 & -65 & 7 & -2 & 1 & 7 & 0  \\
      0 & 4 & -19 & -29 & 130 & -28 & -7 & 0  \\
      1 & 20 & 84 & -46 & 35 & 16 & -7 & 0  \\
      0 & -7 & 16 & 35 & -46 & 84 & 20 & 1  \\
      0 & -7 & -28 & 130 & -29 & -19 & 4 & 0  \\
      0 & 7 & 1 & -2 & 7 & -65 & -11 & 0  \\
      0 & -2 & -8 & -123 & 41 & 27 & 3 & 0  \\
      0 & 0 & -1 & 0 & -1 & 23 & 0 & 0  \\
      0 & 0 & 11 & 59 & -30 & -18 & 0 & 0  \\
      0 & 0 & -6 & -4 & 1 & 3 & 0 & 0  \\
      0 & 0 & 1 & -19 & 15 & 0 & 0 & 0  \\
      0 & 0 & 0 & 8 & -7 & 0 & 0 & 0  \\
      0 & 0 & 0 & -1 & 1 & 0 & 0 & 0  
    \end{pmatrix}
  \end{equation*}
\end{itemize}

\subsection{Complement of ``Simple'' Hyperbolic Knots}
\label{sec:second_appendix}

\subsubsection{
  $K4_{4}$}

\begin{align*}
  & \mbox{
    \raisebox{-1.6cm}{
  \includegraphics[scale=.16]{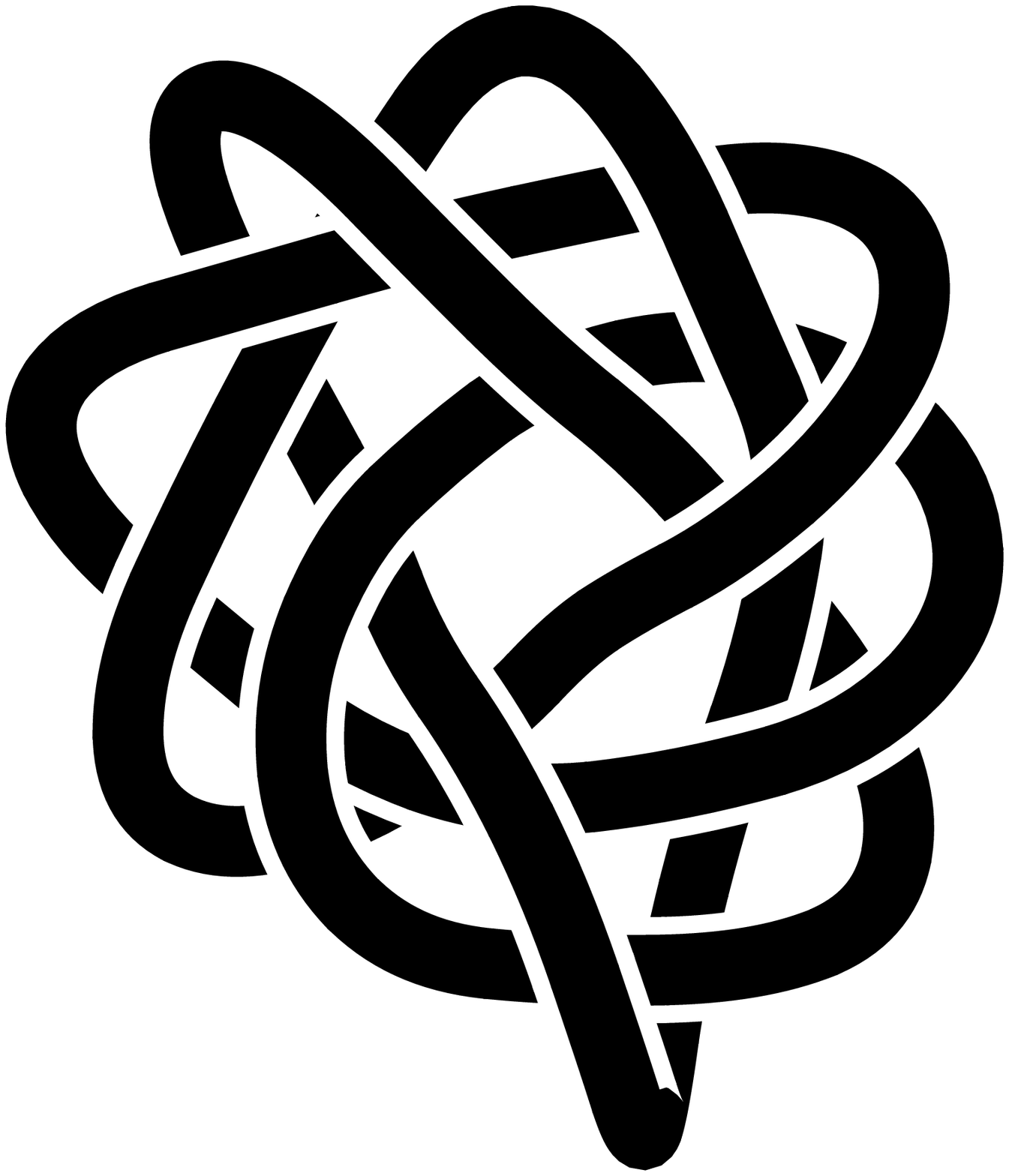}
    }}
  &
  & \Vol(S^3 \setminus \mathcal{K})=  3.608689 \dots
\end{align*}

\begin{itemize}
\item Quantum invariant
  \begin{multline*}
    Z_\gamma \left( \mathcal{M}_u \right)
    = \int\limits_\mathbb{R} 
    \mathrm{d} \boldsymbol{p} \,
    \delta_C(\boldsymbol{p};   u) \,
    \langle p_3 , p_2 | S^{-1} | p_2 , p_1 \rangle \,
    \langle p_{12} , p_5 | S | p_{10} , p_3 \rangle  \,
    \langle p_6 , p_1 | S^{-1} | p_{11} , p_{12} \rangle 
    \\
    \times 
    \langle p_{10} , p_9 | S | p_8 , p_6 \rangle \,
    \langle p_4 , p_7 | S | p_9 , p_{4} \rangle \,
    \langle p_{11} , p_{8} | S | p_7 , p_5 \rangle 
  \end{multline*}

\item
Condition $\delta_C(\boldsymbol{p};u)$
  \begin{equation*}
    p_{12} - p_1 + p_{5} - p_8 - p_6 + p_9 = -  2 \, u
  \end{equation*}

\item Potential function
  \begin{multline*}
    V_{\mathcal{M}}(v,w, x, y,z ; m)
    =
    -\Li\left(\frac{x}{w}\right)
    + \Li\left(\frac{1}{v \, m^2 \, y}\right)
    - \Li\left(\frac{z}{v \, m^2}\right)
    \\
    + \Li\left(m^2 \, x \, y \right)
    +\Li \left( \frac{m^2}{y \, z}\right)
    + \Li\left(\frac{w \, y}{z}\right)
    - \frac{ \pi^2}{3}
    +
    \log \left(\frac{v}{z} \right) \, \log \left( m^4 \, x \, y
    \right)
    \\
    + \log (x) \log\left( \frac{m^2 \, w^2}{x} \right)
    + 2 \log ( w \, y) \log(y)
    + \log\left(m^2\right)  \log \left( m^4 \, w^3 \right)
  \end{multline*}

\item Hyperbolic volume
  \begin{multline*}
    \Im V_{\mathcal{M}}(v,w, x, y,z ; 1)
    =
    -D\left(\frac{x}{w}\right)
    + D\left(\frac{1}{v  \, y}\right)
    - D\left(\frac{z}{v }\right)
    \\
    + D\left( x \, y \right)
    + D \left( \frac{1}{y \, z}\right)
    + D\left(\frac{w \, y}{z}\right)
  \end{multline*}
  with
  \begin{equation*}
    \begin{pmatrix}
      v \\ w\\ x\\ y \\ z
    \end{pmatrix}
    =
    \begin{pmatrix}
      -0.06796 - 1.03267 \, \I
      \\
      -0.597112 + 0.762045  \, \I
      \\
      1.29516  +  0.539127 \, \I
      \\
      0.457778 +  1.02559 \, \I
      \\
      -0.396648 - 0.345221 \, \I
    \end{pmatrix}
  \end{equation*}
\item A-polynomial
  \begin{multline*}
    A(\ell, m )
    =
    1 + ( m^{30} - 2 \, m^{32} + m^{34} ) \,\ell
    + ( -m^{58} + 2 \, m^{60} - 10 \, m^{62} + 4 \, m^{64} - m^{66})
    \, \ell^2
    \\
    +(-2 \, m^{90} + 3 \, m^{92} - m^{96}) \, \ell^3 
    + (m^{120} + 8 \, m^{122} + 6 \, m^{124})  \,\ell^4
    \\
    +(m^{150} - m^{152} - m^{154} + m^{156}) \,  \ell^5
    +
    (2 \, m^{180} - 12 \,  m^{182} - 12 \, m^{186} + 2 \, m^{188})  \,
    \ell^6
    \\
    +(m^{212} - m^{214} - m^{216} + m^{218}) \, \ell^7
    +(6 \,  m^{244} + 8 \, m^{246} + m^{248}) \, \ell^8
    \\
    +(-m^{272} + 3 \,  m^{276} - 2 \, m^{278})  \, \ell^9
    +(-m^{302} + 4 \,  m^{304} - 10 \, m^{406} + 2 \, m^{408} -
    m^{310}) \,  \ell^{10}
    \\
    +(m^{334} - 2 \, m^{336} + m^{338})\, \ell^{11}
    + m^{368}  \, \ell^{12}    
  \end{multline*}
\end{itemize}

\subsubsection{
$K5_1$}

\begin{align*}
  & \mbox{
    \raisebox{-1.6cm}{
      \includegraphics[scale=.16]{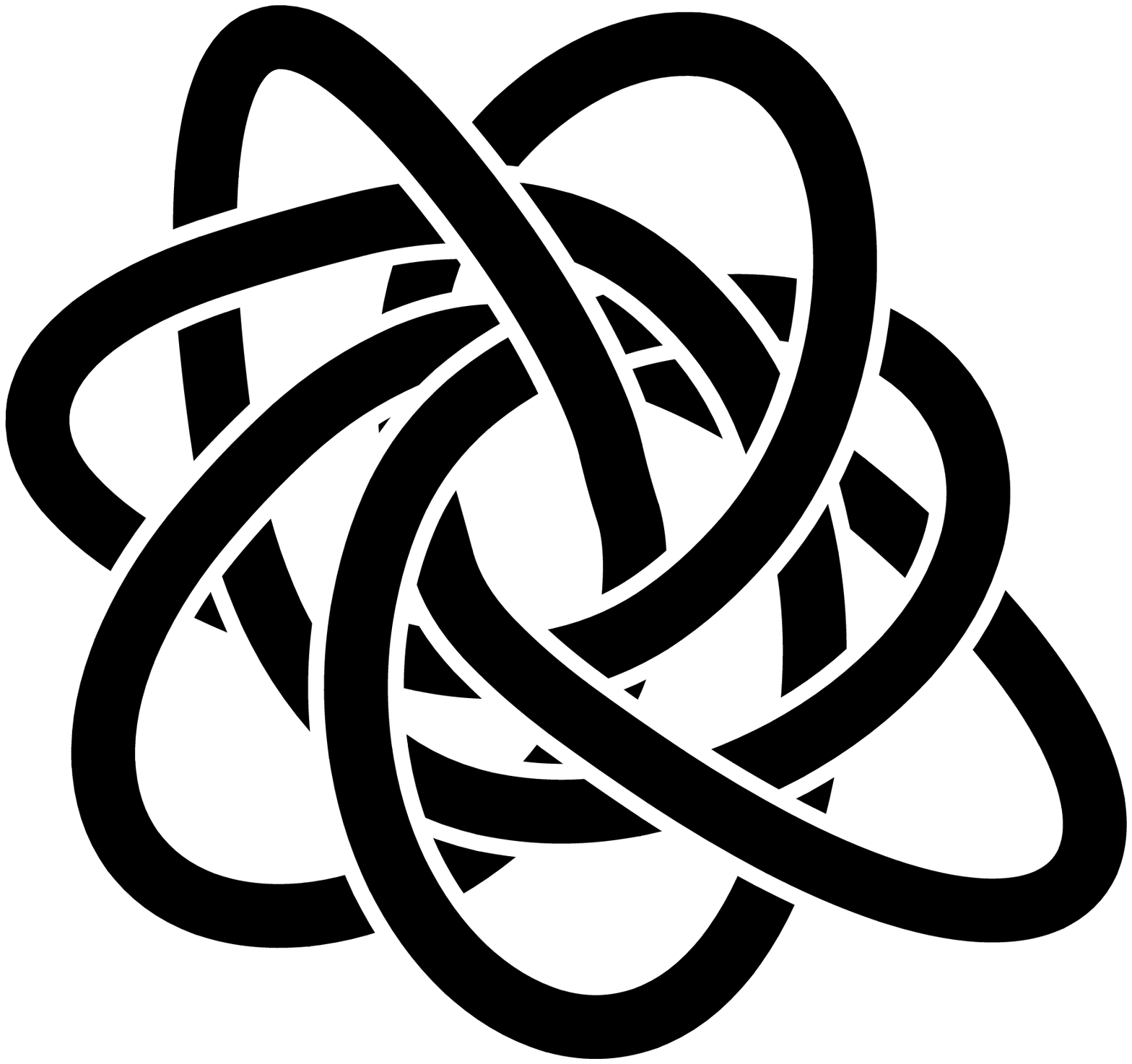}
    }}
  &
  & \Vol(S^3 \setminus \mathcal{K})= 3.4179148 \dots
\end{align*}

\begin{itemize}
\item  Quantum invariant
\begin{multline*}
  Z_\gamma \left( \mathcal{M}_u \right)
  =
  \int\limits_\mathbb{R}
  \mathrm{d} \boldsymbol{p} \,
  \delta_C(\boldsymbol{p}; u ) \,
  \left\langle
    p_5 , p_2 \middle| S^{-1} \middle| p_1 , p_3
  \right\rangle \,
  \left\langle
    p_1 , p_9 \middle| S^{-1} \middle| p_9 , p_2
  \right\rangle \,
  \left\langle
    p_4 , p_3 \middle| S^{-1} \middle| p_5 , p_7
  \right\rangle
  \\
  \times
  \left\langle
    p_{10} , p_6 \middle| S^{-1} \middle| p_4 , p_{11}
  \right\rangle \,
  \left\langle
    p_{12} , p_7 \middle| S^{-1} \middle| p_6 , p_{10}
  \right\rangle \,
  \left\langle
    p_{11} , p_8 \middle| S^{-1} \middle| p_8 , p_{12}
  \right\rangle 
\end{multline*}

\item
Condition $\delta_C(\boldsymbol{p};u)$
\begin{equation*}
  p_{11} - p_6 - p_{10}+ p_7 = - 2 \, u
\end{equation*}

\item Potential function
\begin{multline*}
  V_\mathcal{M}(v,w,x,y,z;m)
  =
  - \Li\left( \frac{w}{x^2}\right)
  - \Li\left( \frac{w \, x}{y} \right)
  - \Li\left( \frac{y^2}{m^4 \, v \, z}\right)
  - \Li\left( \frac{y^2}{w \,z} \right)
  \\
  - \Li\left( \frac{z}{v \, y}\right)
  - \Li\left( \frac{m^2 \, z }{v \, y} \right)
  + \pi^2
  + 5 \log \left( m^2 \right) \log\left(\frac{y}{z}\right)
  \\
  + 2 \log w \log \left( \frac{x \, y}{w} \right)
  +2 \log z \log \left( \frac{y^3}{z} \right)
  - 2 \left( \log( m^2) \right)^2
  -2 \left( \log x \right)^2
  - 5 \left( \log y \right)^2
\end{multline*}

\item Hyperbolic volume
\begin{equation*}
  \Im V_{\mathcal{M}}(v,w,x,y,z; 1)
  =
  - D\left( \frac{w}{x^2}\right)
  - D\left( \frac{w \, x}{y} \right)
  - D\left( \frac{y^2}{ v \, z}\right)
  - D\left( \frac{y^2}{w \,z} \right)
  - 2 \, D\left( \frac{z}{v \, y}\right)
\end{equation*}
with
\begin{equation*}
  \begin{pmatrix}
    v \\ w\\ x\\ y \\ z
  \end{pmatrix}
  =
  \begin{pmatrix}
    0.465534 - 0.473866 \, \I
    \\
    0.693244 + 0.159750 \, \I
    \\
    -1.085877 - 0.175545 \, \I
    \\
     -0.952444 - 0.928780 \, \I
     \\
     -0.907927 + 0.840443 \, \I
  \end{pmatrix}
\end{equation*}

\item A-polynomial
  \begin{multline*}
    A_\mathcal{K}(\ell,m)
    =
    -1 + 
    \left(-  m^{32} +  m^{34} \right) \, \ell
    +\left(
      9 \,  m^{64} - 3  \,  m^{66} +    m^{68}
    \right) \,  \ell^2
    \\
    + \left(
      m^{92} - 3 \,  m^{94} + 12 \,  m^{96} - 14 \, m^{98}
      + 5 \, m^{100} - 
      m^{102}
    \right) \, \ell^3
    \\
    + \left(
      m^{124} - 7 \,  m^{126} - 18 \,  m^{128} + 5 \,  m^{130 }
      - 2 \,  m^{132}
    \right) \, \ell^4
    + \left(
      -  m^{156} - 7 \,  m^{158} + 2  \, m^{160} + 6 \,   m^{162}
    \right)
    \, \ell^5
    \\
    +\left(
      -  m^{188} + 17 \,  m^{190} + 20 \,  m^{192} - 2 \, m^{194} +
      m^{196}
    \right)
    \, \ell^6
    \\
    + \left(
      - 2 \,  m^{218} + 14 \, m^{220} - 12 \,  m^{222} + 12 \,  m^{224}
      - 14 \,    m^{226} + 2 \, m^{228}
    \right) \, \ell^7
    \\
    + \left(
      -  m^{250} + 2 \,  m^{252} - 20 \, m^{254} - 17 \, m^{256} +
      m^{258}
    \right) \ell^8
    \\
    +\left(
      - 6 \,  m^{284} - 2 \, m^{286} + 7 \, m^{288} +  m^{290}
    \right) \ell^9
    + \left(
      2  \, m^{314} - 5 \, m^{316} + 18 \, m^{318}
      + 7 \, m^{320} -  m^{322}
    \right)
    \, \ell^{10}
    \\
    + 
    \left(
      m^{344} - 5 \,  m^{346} + 14 \, m^{348} - 12 \,  m^{350} + 3  \, m^{352}
      - m^{354}
    \right) \, \ell^{11}
    \\
    + \left(
      -  m^{378} + 3  \,  m^{380} - 9  \, m^{382}
    \right) \ell^{12}
    + \left(
      -  m^{412} +  m^{414}
    \right) \ell^{13}
    + \ell^{14} \,  m^{446}
  \end{multline*}
\end{itemize}
\subsubsection{
  $K5_9$ or $10_{132}$}

\begin{align*}
  & \mbox{
    \raisebox{-1.6cm}{
      \includegraphics[scale=.16]{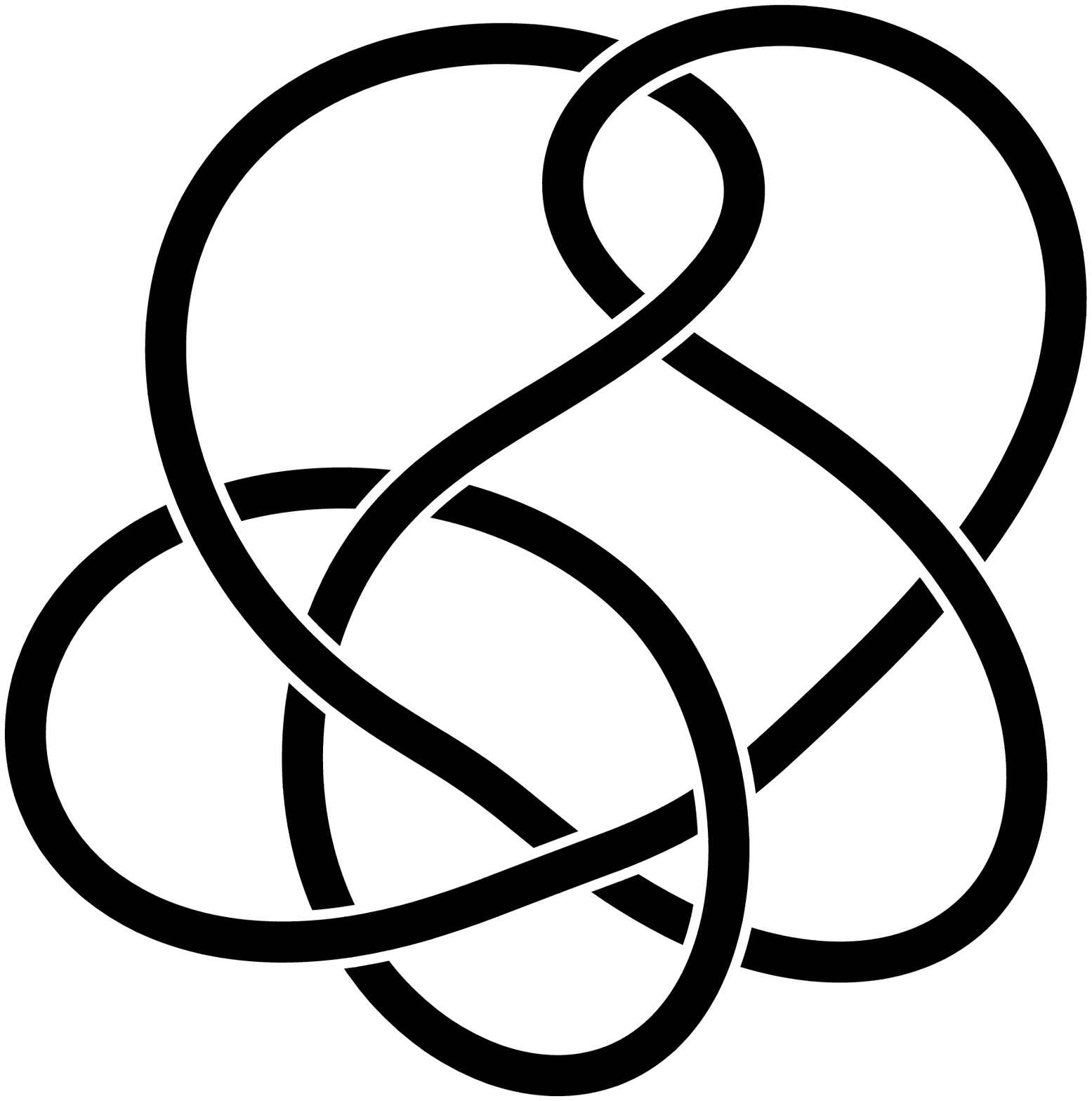}
    }}
  &
  & \Vol(S^3 \setminus \mathcal{K}) = 4.05686 \dots
\end{align*}

\begin{itemize}
\item  Quantum invariant
\begin{multline*}
  Z_\gamma \left( \mathcal{M}_u \right)
  = \int\limits_\mathbb{R} 
  \mathrm{d} \boldsymbol{p} \,
  \delta_C(\boldsymbol{p}; u ) \,
    \langle p_1, p_9 | S^{-1} | p_2 , p_3 \rangle \,
    \langle p_4, p_2 | S^{-1} | p_9 , p_5 \rangle 
    \\
    \times
    \langle p_6, p_7 | S^{-1} | p_7 , p_8 \rangle \,
    \langle p_5, p_3 | S | p_6 , p_{10} \rangle \,
    \langle p_{10}, p_8 | S | p_4 , p_1 \rangle 
  \end{multline*}

\item
Condition $\delta_C(\boldsymbol{p};u)$
\begin{equation*}
  p_2 + p_8 + p_{10} = p_1 + p_5 + p_9 -  2\, u
\end{equation*}

\item Potential function
  \begin{multline*}
    V_{\mathcal{M}}(w,x,y,z; m)
    =
    \Li(w) - \Li\left(\frac{w}{y}\right)
    - \Li\left(\frac{y}{x}\right)
    - \Li(z) + \Li\left(\frac{y \, z}{m^2}\right)
    + \frac{\pi^2}{6} 
    \\
    + \left( \log \left( m^2 \right) \right)^2
    - \log \left( m^2 \right) \, \log(x \, y \, w)
    + \log w \, \log \left(\frac{y}{z}\right)
    +2 \log y  \log(x \, z)
  \end{multline*}

\item Hyperbolic volume
  \begin{equation*}
    \Im V_{\mathcal{M}}(w,x,y,z; 1)
    =
    D(w) - D\left(\frac{w}{y}\right)
    - D\left(\frac{y}{x}\right)
    - D(z) + D\left(y \, z\right)
  \end{equation*}
  with
  \begin{equation*}
    \begin{pmatrix}
      w \\ x\\ y \\ z
    \end{pmatrix}
    =
    \begin{pmatrix}
       -0.0498076 + 0.754729 \, \I
       \\
       -1.54094 - 1.35872 \, \I
       \\
       -0.821578 - 0.131699 \, \I
       \\
       -0.0844626 - 0.905094 \, \I
    \end{pmatrix}
  \end{equation*}

\item A-polynomial
  \begin{equation*}
    \begin{pmatrix}
      0 & 0 & 0 & 0 & 1 & 0 & 0 & 0 & 0 \\
      0 & 0 & 0 & -1 & -7 & 0 & 0 & 0 & 0 \\
      0 & 0 & -1 & 3 & 7 & -4 & -2 & 0 & 0 \\
      0 & 1 & 4 & 2 & 13 & 4 & -6 & -3 & -1 \\
      0 & -1 & -5 & -7 & -6 & 5 & 13 & 1 & 0 \\
      0 & 0 & -5 & -12 & -18 & 9 & 9 & 4 & 0 \\
      0 & 0 & 12 & 13 & -15 & -10 & -7 & -4 & 0 \\
      0 & 4 & 7 & 10 & 15 & -13 & -12 & 0 & 0 \\
      0 & -4 & -9& -9  & 18 & 12 & 5 & 0 & 0 \\
      0 & -1 & -13 & -5  & 6 & 7 & 5 & 1 & 0 \\
      1 & 3 & 6 & -4  & -13 & -2 & -4 & -1 & 0 \\
      0 & 0 & 2 & 4  & -7 & -3 & 1 & 0 & 0 \\
      0 & 0 & 0 & 0  & 7 & 1 & 0 & 0 & 0 \\
      0 & 0 & 0 & 0  & -1 & 0 & 0 & 0 & 0 \\
    \end{pmatrix}
  \end{equation*}
\end{itemize}


\subsubsection{
  $K5_{12}$ or $8_{20}$}

\begin{align*}
  & \mbox{
    \raisebox{-1.6cm}{
      \includegraphics[scale=.16]{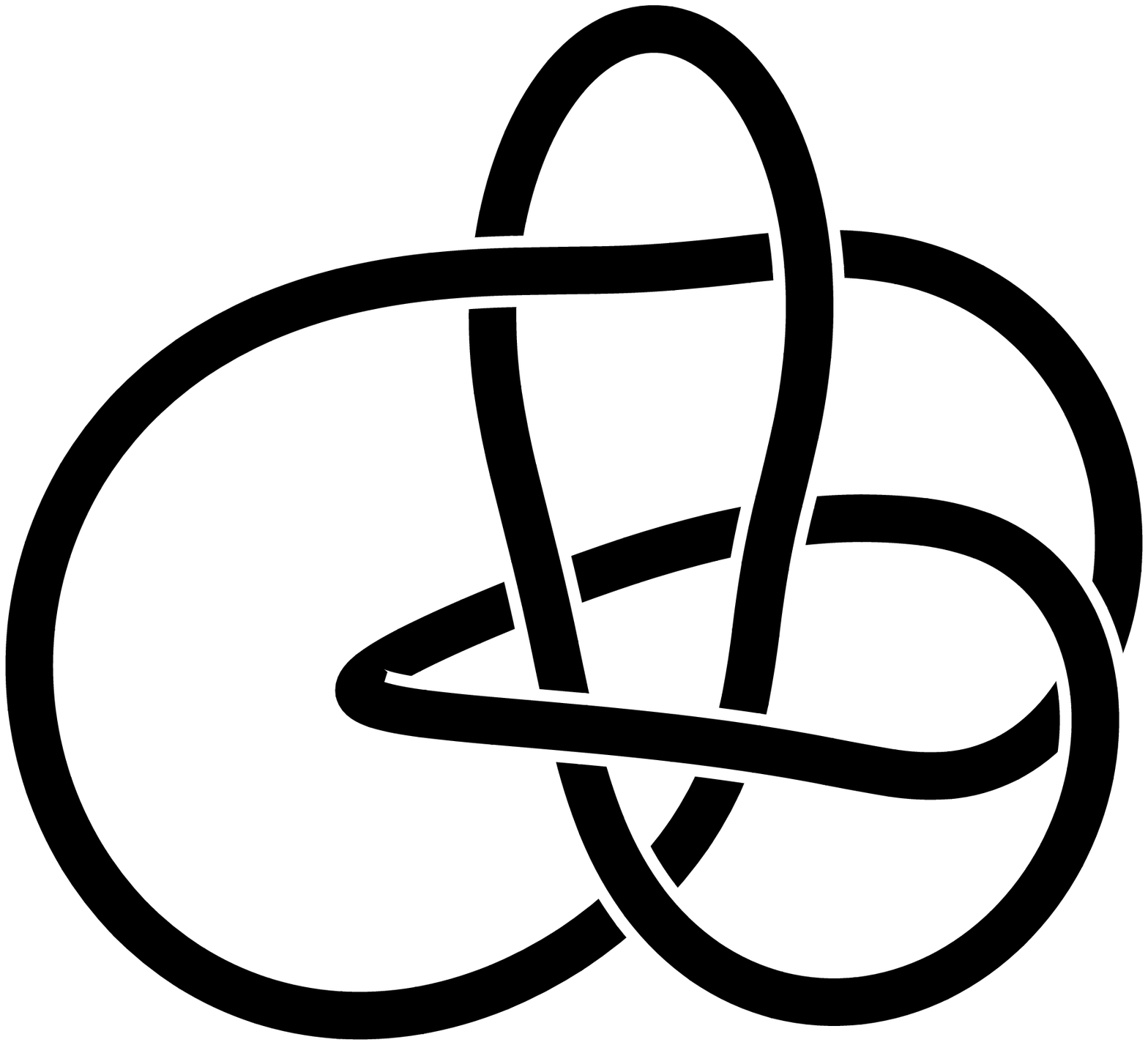}
    }}
  &
  & \Vol(S^3 \setminus \mathcal{K})= 4.12490 \dots
\end{align*}

\begin{itemize}
\item Quantum invariant
  \begin{multline*}
    Z_\gamma \left( \mathcal{M}_u \right)
    = \int\limits_\mathbb{R} 
 \mathrm{d} \boldsymbol{p} \,
  \delta_C(\boldsymbol{p}; u) \,
    \langle p_8 , p_2 | S | p_1 , p_6 \rangle \,
    \langle p_1 , p_{10} | S^{-1} | p_{10} , p_4 \rangle
    \\
    \times 
    \langle p_9 , p_3 | S | p_2 , p_9 \rangle \,
    \langle p_4 , p_7 | S^{-1} | p_5 , p_3 \rangle \,
    \langle p_6 , p_5 | S^{-1} | p_7 , p_8 \rangle 
  \end{multline*}

\item
Condition $\delta_C(\boldsymbol{p};u)$
\begin{equation*}
  p_3 - p_7 + p_6 - p_2 - p_8 + p_5 = -  2 \, u
\end{equation*}

\item Potential function
  \begin{multline*}
    V_{\mathcal{M}}(w,x,y,z; m)
    \\
    =
    - \Li\left(\frac{x}{m^2 \, w}\right)
    - \Li\left(\frac{m^4 \, w}{x \, y}\right)
    + \Li\left(\frac{x}{y}\right)
    - \Li\left(\frac{z}{w \, y}\right)
    + \Li\left(\frac{z}{x \, y}\right) + \frac{\pi^2}{6}
    \\
    -2 \left( \log \left(m^2 \right) \right)^2 - 2 \left( \log w \right)^2
    + \left( \log x \right)^2
    \\
    + \log \left(m^2 \right) \, \log \left( \frac{x \, z^2}{w^5} \right)
    + 2 \log w \log x +2 \log w \log z - 2 \log x \log z
  \end{multline*}

\item Hyperbolic volume
  \begin{equation*}
    \Im V_{\mathcal{M}}(w,x,y,z; 1)
    =
    - D\left(\frac{x}{ w}\right)
    - D\left(\frac{ w}{x \, y}\right)
    + D\left(\frac{x}{y}\right)
    - D\left(\frac{z}{w \, y}\right)
    + D\left(\frac{z}{x \, y}\right) 
  \end{equation*}
  with
  \begin{equation*}
    \begin{pmatrix}
      w \\ x\\ y \\ z
    \end{pmatrix}
    =
    \begin{pmatrix}
      -0.723387 - 0.90034 \, \I
      \\
       -0.637406 + 0.318768 \,\I
       \\
       -0.483596 + 0.741071\, \I
       \\
        1.08906 - 0.727199  \, \I
    \end{pmatrix}
  \end{equation*}

\item A-polynomial
  \begin{equation*}
    \begin{pmatrix}
      0 & -1 & 1 & 0 & 0 & 0 \\
      0 & 1 & -5 & 0 & 0 & 0 \\
      0 & -2 & 0 & -2 & 0 & 0 \\
      0 & 2 & 3 & -4 & -1 & 0 \\
      0 & 1 & -3 & 0 & -5 & -1 \\
      -1 & -5 & 0 & -3 & 1 & 0 \\
      0 & -1 & -4 & 3 & 2 & 0 \\
      0 & 0 & -2 & 0 & -2 & 0 \\
      0 & 0 & 0 & -5 & 1 & 0 \\
      0 & 0 & 0 & 1 & -1 & 0 \\
    \end{pmatrix}
  \end{equation*}
\end{itemize}

\subsubsection{
  $K5_{13}$ 
}

\begin{align*}
  & \mbox{
    \raisebox{-1.6cm}{
      \includegraphics[scale=.16]{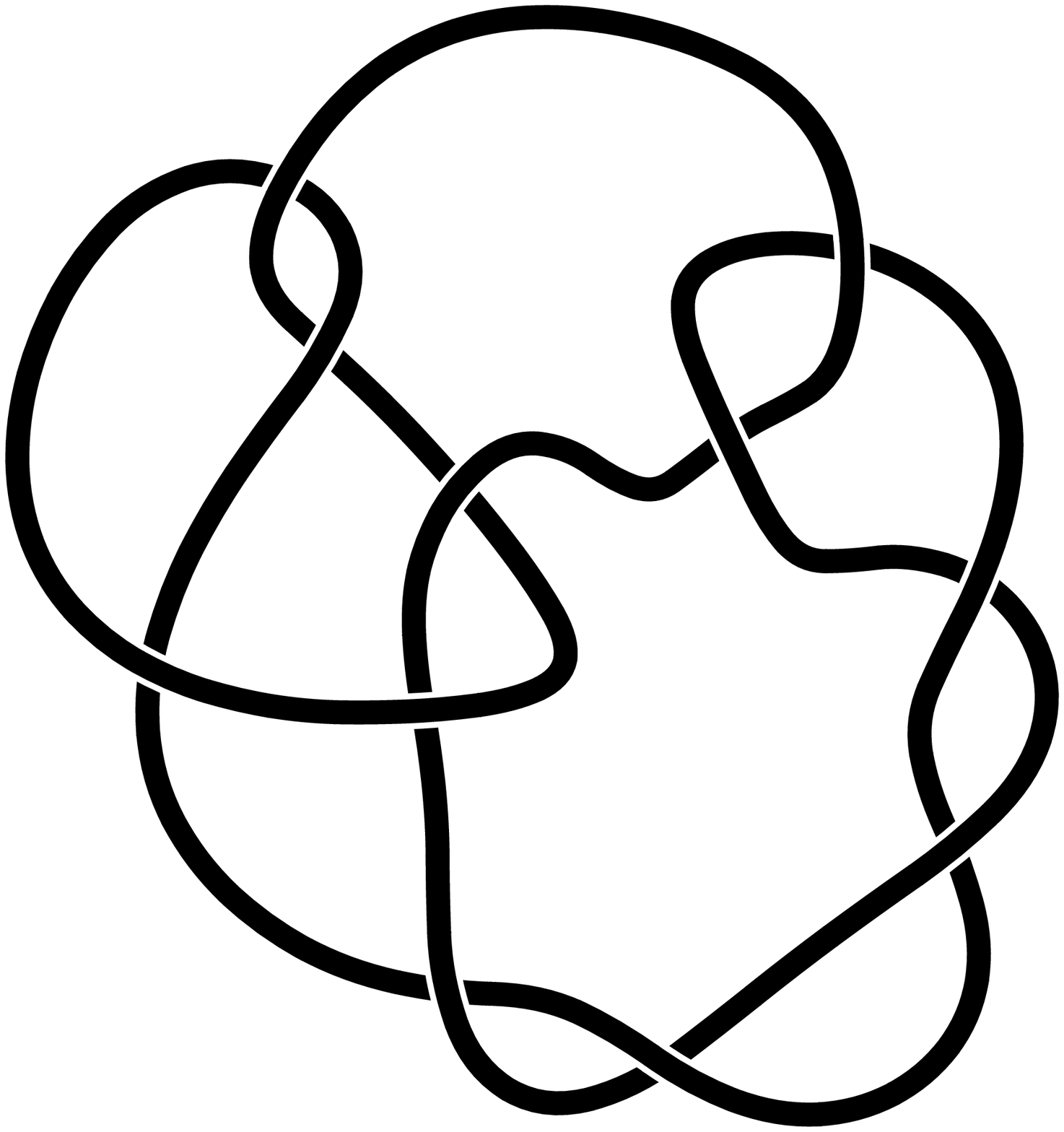}
    }}
  &
  & \Vol(S^3 \setminus \mathcal{K})= 4.1249 \dots
\end{align*}

\begin{itemize}
\item Quantum invariant
\begin{multline}
  Z_\gamma \left( \mathcal{M}_u \right)
  = \int\limits_\mathbb{R} 
  \mathrm{d} \boldsymbol{p} \,
  \delta_C(\boldsymbol{p};  u ) \,
    \langle p_2 , p_9 | S^{-1} | p_1 , p_3 \rangle \,
    \langle p_7 , p_6 | S^{-1} | p_6 , p_2 \rangle 
    \\
    \times 
    \langle p_4 , p_5 | S^{-1} | p_9 , p_7 \rangle \,
    \langle p_8 , p_1 | S | p_4 , p_{10} \rangle \,
    \langle p_{10} , p_3 | S | p_5 , p_8 \rangle 
  \end{multline}

\item
Condition $\delta_C(\boldsymbol{p};u)$
\begin{equation*}
  p_1 + p_5 + p_{8} = p_7 + p_9 + p_{10} - 2 \, u 
\end{equation*}

\item Potential function
\begin{multline*}
  V_{\mathcal{M}}(w,x,y,z ; m)
  =
  - \Li\left(m \, x\right)
  + \Li\left(\frac{x}{w}\right)
  - \Li\left(\frac{1}{w^2 \, y}\right)
  - \Li\left(\frac{z}{y}\right)
  \\
  + \Li\left(\frac{z}{w \, y}\right)
  + \frac{\pi^2}{6}
  - \left( \log \left( m^2 \right) \right)^2
  \\
  - \log \left( m^2 \right) \log(w \, z)
  + \log x \, \log \left( \frac{z}{y \, w} \right)
  - 2 \log w \log z
\end{multline*}

\item Hyperbolic volume
\begin{equation*}
  \Im V_{\mathcal{M}}(w,x,y,z ; 1)
  =
  - D\left( x\right)
  + D\left(\frac{x}{w}\right)
  - D\left(\frac{1}{w^2 \, y}\right)
  - D\left(\frac{z}{y}\right)
  + D\left(\frac{z}{w \, y}\right)
\end{equation*}
with
\begin{equation*}
  \begin{pmatrix}
    w \\ x\\ y \\ z
  \end{pmatrix}
  =
  \begin{pmatrix}
    -0.812447 + 0.173142 \, \I
    \\
    -0.0890598 - 0.727199  \, \I
    \\
    -1.71268 + 1.30259 \, \I
    \\
     1.09977 + 1.12945  \, \I
  \end{pmatrix}
\end{equation*}

\item A-polynomial
  \begin{equation*}
    \begin{pmatrix}
      0 & 1 & 0 & 0 & 0 & 0 & 0 & 0 & 0 &  0 \\
      0 &      -1 & 0 & 0 & 0 & 0 & 0 & 0 & 0 &  0 \\
      0 &      0 & 1 & 2 & 0 & 0 & 0 & 0 & 0 &  0 \\
      0 &      0 & -10 & 1 & 3 & -1 & 0 & 0 & 0 &  0 \\
      0 &      0 & 8 & -14 & -11 & 8 & 0 & 0 & 0 &  0 \\
      0 &      -5 & 11 & 5 & -1 & -12 & 0 & 0 & 0 &  0 \\
      -1 &      5 & -10 & 14 & 29 & -12 & 1 & 0 & 0 &  0 \\
      0 &      0 & -5 & 4 & 5 & 16 & -1 & 0 & 0 &  0 \\
      0 &      -1  & 9 & -13 & -19 & 25 & -7 & 2 & 0 &  0 \\
      0 &      0 & 5 & 3 & -33 & -3 & 9 & -7 & 0 &  0 \\
      0 &      0 & -7 & 9 & -3 & -33 & 3 & 5 & 0 &  0 \\
      0 &      0 & 2 & -7 & 25 & -19 & -13 & 9 & -1 &  0 \\
      0 &      0 & 0 & -1 & 16 & 5 & 4 & -5 & 0 &  0 \\
      0 &      0 & 0 & 1 & -12 & 29 & 14 & -10 & 5 &  -1 \\
      0 &      0 & 0 & 0 & -12 & -1 & 5 & 11 & -5 &  0 \\
      0 &      0 & 0 & 0 & 8 & -11 & -14 & 8 & 0 &  0 \\
      0 &      0 & 0 & 0 & -1 & 3 & 1 & -10 & 0 &  0 \\
      0 &      0 & 0 & 0 & 0 & 0 & 2 & 1 & 0 &  0 \\
      0 &      0 & 0 & 0 & 0 & 0 & 0 & 0 & -1 &  0 \\
      0 &      0 & 0 & 0 & 0 & 0 & 0 & 0 & 1 &  0 \\
    \end{pmatrix}
  \end{equation*}
\end{itemize}


\subsubsection{
  $K5_{21}$ or  $9_{46}$}

\begin{align*}
  & \mbox{
    \raisebox{-1.6cm}{
      \includegraphics[scale=.16]{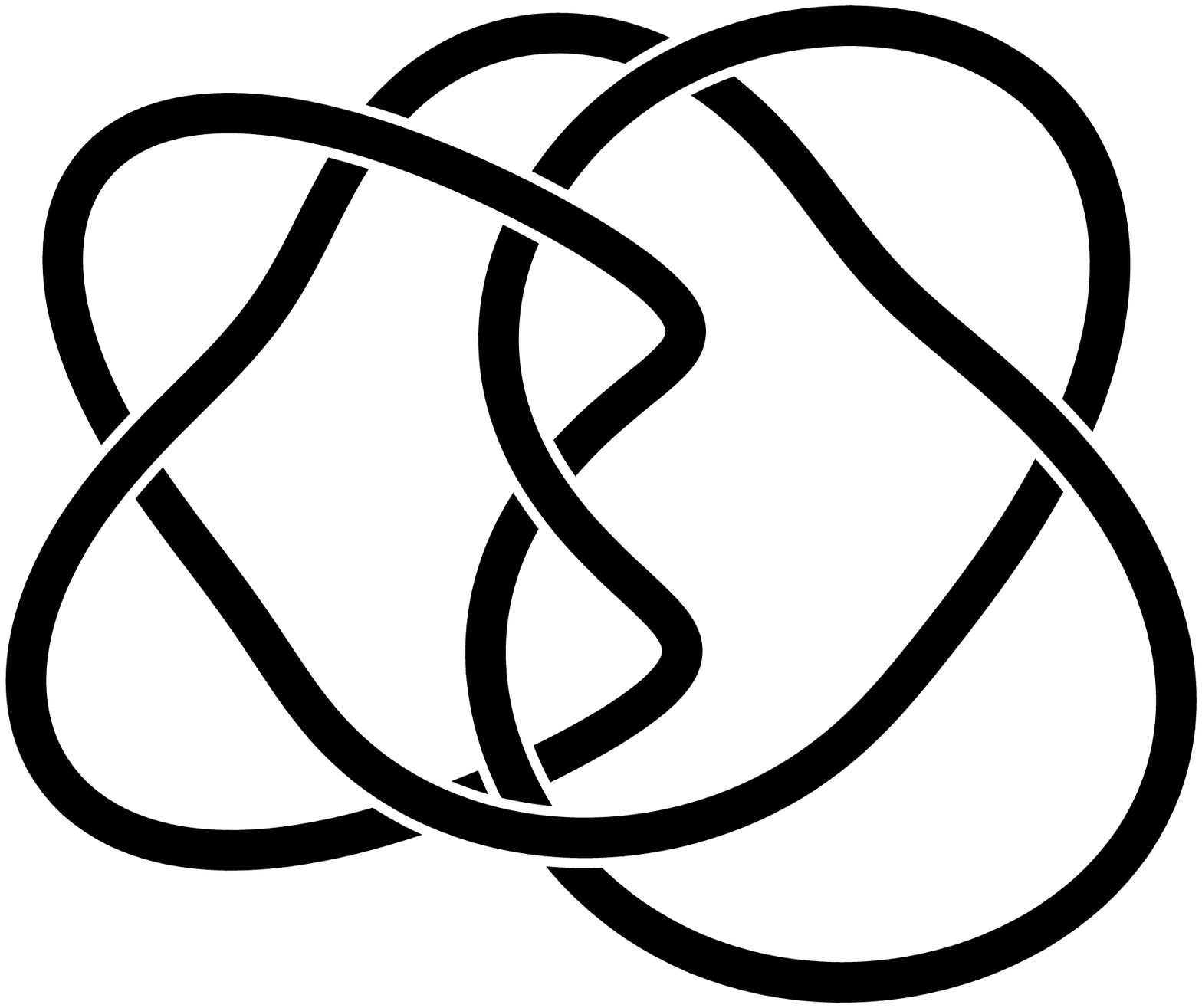}
    }}
  &
  & \Vol(S^3 \setminus \mathcal{K})= 4.7517 \cdots
\end{align*}

\begin{itemize}
\item  Quantum invariant
  \begin{multline*}
    Z_\gamma \left( \mathcal{M}_u \right)
    = \int\limits_\mathbb{R} 
 \mathrm{d} \boldsymbol{p} \,
  \delta_C(\boldsymbol{p}; u) \,
    \langle p_4 , p_1 | S | p_6, p_7 \rangle \,
    \langle p_2 , p_9 | S^{-1} | p_1 , p_8 \rangle 
    \\
    \times 
    \langle p_7 , p_8 | S^{-1} | p_3, p_4 \rangle \,
    \langle p_6 , p_5 | S^{-1} | p_9, p_{10} \rangle \,
    \langle p_3 , p_{10} | S^{-1} | p_5, p_2 \rangle 
  \end{multline*}

\item
Condition $\delta_C(\boldsymbol{p};u)$
  \begin{equation*}
    p_1 + p_4 + p_5 = p_7 + p_9 + p_{10} -  2 \, u 
  \end{equation*}

\item Potential function
  \begin{multline*}
    V_{\mathcal{M}}(w, x, y,z; m)
    =
    - \Li\left(\frac{w^2 \, x}{m^2}\right)
    + \Li\left(\frac{m^2 \, y}{x}\right)
    - \Li\left(\frac{m^2}{w^2 \, x \, z}\right)
    \\ - \Li\left(\frac{z}{w}\right)
    - \Li\left(\frac{w \, y \, z}{m^2}\right)
    + \frac{\pi^2}{2}
    + \left( \log \left( m^2 \right) \right)^2 -2 \, \left( \log w \right)^2
    - \left( \log z \right)^2
    \\
    + \log \left( m^2 \right) \log \left( \frac{y \, z}{x}\right)
    - \log w \, \log x
    - \log (y \, z) \, \log(x \, w)
  \end{multline*}

\item Hyperbolic volume
  \begin{equation*}
    \Im V_{\mathcal{M}}(w, x, y,z; 1)
    =
    - D\left(w^2 \, x\right)
    + D\left(\frac{ y}{x}\right)
    - D\left(\frac{1}{w^2 \, x \, z}\right)
    - D\left(\frac{z}{w}\right)
    - D\left(w \, y \, z\right)
  \end{equation*}
  with
  \begin{equation*}
    \begin{pmatrix}
      w \\ x\\ y \\ z
    \end{pmatrix}
    =
    \begin{pmatrix}
       -1.21844 + 0.168108 \, \I
       \\
       0.640448 - 0.637204 \, \I
       \\
       1.0 \\
       -0.445837 + 0.526085  \, \I
    \end{pmatrix}
  \end{equation*}

\item A-polynomial
  \begin{equation*}
    \begin{pmatrix}
      0 & 0 & -1 & 1 & 0 \\
      0 & 2 & 5 & -1 & 0 \\
      -1 & -5 & -3 & 2 & 0 \\
      0 & 0 & -5 & -5 & 0 \\
      0 & 2 & 2 & 2 & 0 \\
      0 & -5 & -5 & 0 & 0 \\
      0 & 2 & -3 & -5 & -1 \\
      0 & -1 & 5 & 2 & 0 \\
      0 & 1 & -1 & 0 & 0 \\
    \end{pmatrix}
  \end{equation*}
\end{itemize}

\subsubsection{
 $K5_{22}$ or  $10_{139}$}

\begin{align*}
  & \mbox{
    \raisebox{-1.6cm}{
      \includegraphics[scale=.16]{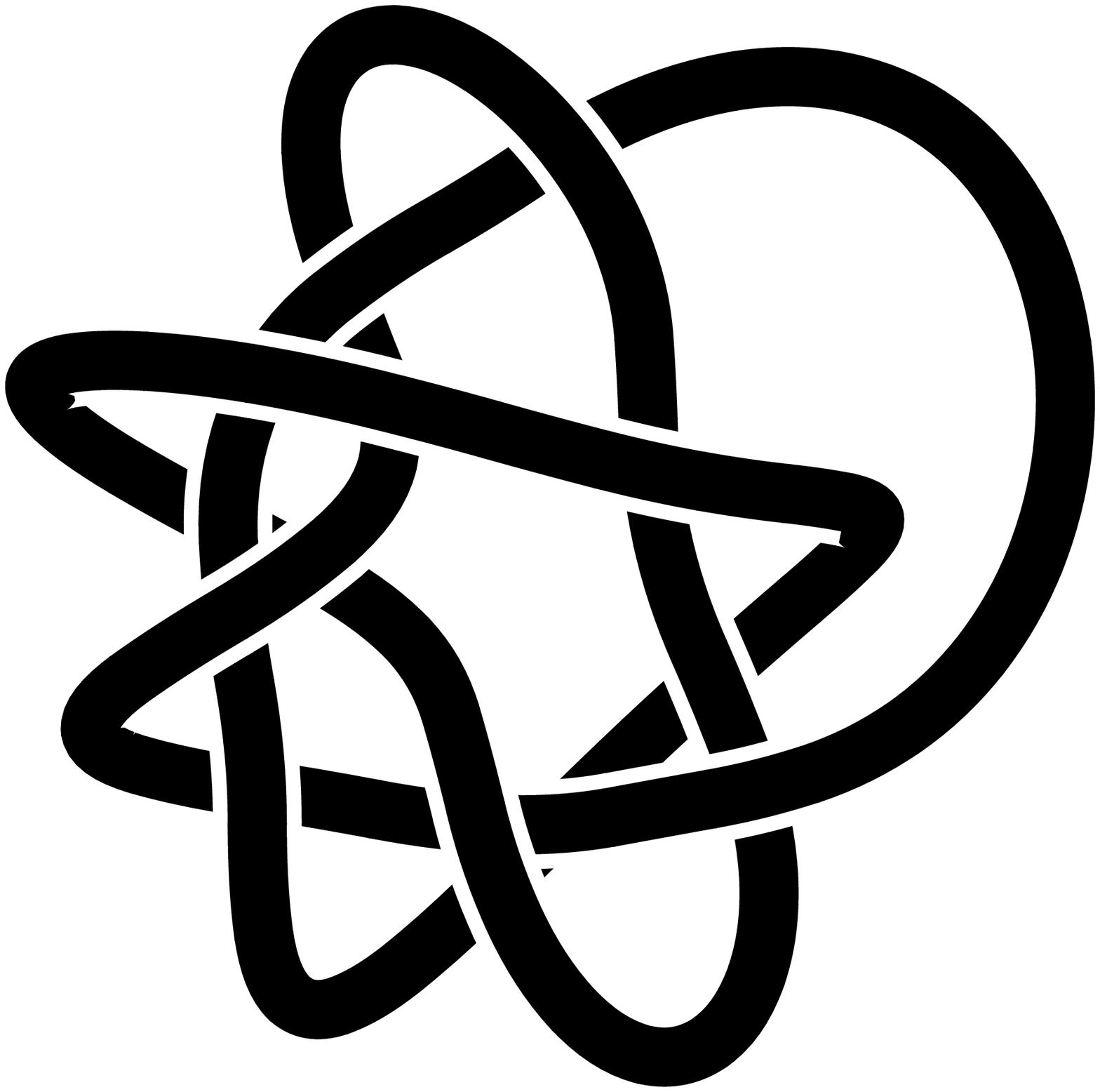}
    }}
  &
  & \Vol(S^3 \setminus \mathcal{K})= 4.85117 \dots
\end{align*}

\begin{itemize}
\item  Quantum invariant
  \begin{multline*}
    Z_\gamma \left( \mathcal{M}_u \right)
    = \int\limits_\mathbb{R} 
 \mathrm{d} \boldsymbol{p} \,
  \delta_C(\boldsymbol{p}; u ) \,
    \langle p_5 , p_9 | S | p_8 , p_1 \rangle \,
    \langle p_1 , p_3 | S | p_2 , p_9 \rangle 
    \\
    \times 
    \langle p_6 , p_7 | S | p_4 , p_3 \rangle \,
    \langle p_2 , p_8 | S | p_7 , p_{10} \rangle \,
    \langle p_4 , p_{10} | S | p_5 , p_6 \rangle 
  \end{multline*}

\item
Condition $\delta_C(\boldsymbol{p};u)$
  \begin{equation*}
    p_7 + p_9 + p_{10} = p_1 + p_3 + p_8 -  2 \, u
  \end{equation*}

\item Potential function
  \begin{multline*}
    V_{\mathcal{M}}(w, x, y,z ; m)
    =
    \Li\left(\frac{1}{m^2 \, x}\right)
    + \Li\left(\frac{m^4 \, w}{y}\right)
    + \Li\left(\frac{w}{x \, z}\right)
    \\
    + \Li\left(\frac{z}{w}\right)
    + \Li\left(\frac{w \, z}{y}\right)
    - \frac{5 \, \pi^2}{6}
    +
    2 \, \left( \log \left(m^2 \right) \right)^2
    \\
    + \log \left(m^2 \right) \, \log\left(\frac{x}{y^2 \, z}\right)
    + \log x \, \log \left(\frac{z}{y}\right)
    + \log w \, \log\left(\frac{w}{z}\right)
  \end{multline*}

\item Hyperbolic volume
  \begin{equation*}
    \Im V_{\mathcal{M}}(w, x, y,z ; 1)
    =
    D\left(\frac{1}{ x}\right)
    + D\left(\frac{ w}{y}\right)
    + D\left(\frac{w}{x \, z}\right)
    + D\left(\frac{z}{w}\right)
    + D\left(\frac{w \, z}{y}\right)
  \end{equation*}
  with
  \begin{equation*}
    \begin{pmatrix}
      w \\ x\\ y \\ z
    \end{pmatrix}
    =
    \begin{pmatrix}
       0.660443 - 0.716885 \, \I
       \\
        -0.0777392 - 0.946923  \, \I
        \\
        -0.460355 - 1.13932 \, \I
        \\
         1.0
    \end{pmatrix}
  \end{equation*}

\item A-polynomial
  \begin{equation*}
    \begin{pmatrix}
      1 & 0  & 0 & 0 & 0 \\
      0 & 0 & 0 & 0 & 0 \\
      0 & 0 & 0 & 0 & 0 \\
      0 & 0 & 0 & 0 & 0 \\
      0 & 0 & 0 & 0 & 0 \\
      0 & 0 & 0 & 0 & 0 \\
      0 & -1 & 0 & 0 & 0 \\
      0 & 7 & 0 & 0 & 0 \\
      0 & -3 & 0 & 0 & 0 \\
      0 & 1 & 0 & 0 & 0 \\
      0 & 0 & 0 & 0 & 0 \\
      0 & 0 & 0 & 0 & 0 \\
      0 & 0 & 0 & 0 & 0 \\
      0 & 0 & 0 & 0 & 0 \\
      0 & 0 & 6 & 0 & 0 \\
      0 & 0 & 0 & 0 & 0 \\
      0 & 0 & 0 & 0 & 0 \\
      0 & 0 & 0 & 0 & 0 \\
      0 & 0 & 0 & 0 & 0 \\
      0 & 0 & 0 & 1 & 0 \\
      0 & 0 & 0 & -3 & 0 \\
      0 & 0 & 0 & 7 & 0 \\
      0 & 0 & 0 & -1 & 0 \\
      0 & 0 & 0 & 0 & 0 \\
      0 & 0 & 0 & 0 & 0 \\
      0 & 0 & 0 & 0 & 0 \\
      0 & 0 & 0 & 0 & 0 \\
      0 & 0 & 0 & 0 & 0 \\
      0 & 0 & 0 & 0 & 1 \\
    \end{pmatrix}
  \end{equation*}
\end{itemize}

\subsubsection{
  $K6_{10}$}

\begin{align*}
  & \mbox{
    \raisebox{-1.6cm}{
      \includegraphics[scale=.16]{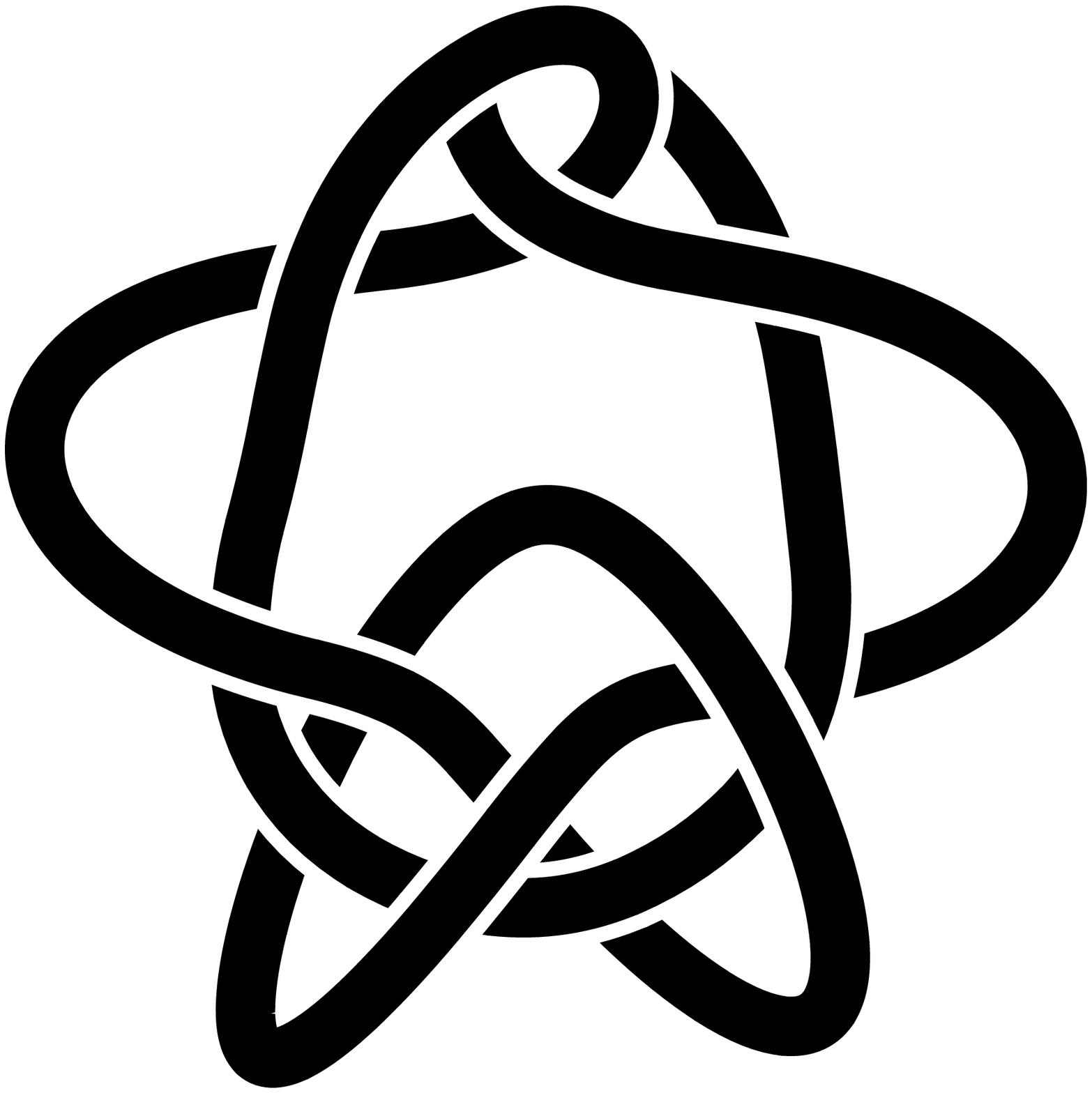}
    }}
  &
  & \Vol(S^3 \setminus \mathcal{K})= 4.40083 \cdots
\end{align*}

\begin{itemize}
\item Quantum invariant
  \begin{multline*}
    Z_\gamma \left( \mathcal{M}_u \right)
    = \int\limits_\mathbb{R} 
    \mathrm{d} \boldsymbol{p} \,
    \delta_C(\boldsymbol{p}; u ) \,
    \left\langle
      p_2 , p_6 \middle| S^{-1} \middle| p_4 , p_1
    \right\rangle \,
    \left\langle
      p_1 , p_{12} \middle| S^{-1} \middle| p_{10} , p_3
    \right\rangle \,
    \left\langle
      p_3 , p_8 \middle| S \middle| p_6 , p_{11}
    \right\rangle
    \\
    \times
    \left\langle
      p_7 , p_9 \middle| S^{-1} \middle| p_9 , p_8
    \right\rangle \,
    \left\langle
      p_{10} , p_4 \middle| S^{-1} \middle| p_{12} , p_5
    \right\rangle \,
    \left\langle
      p_{11} , p_5 \middle| S \middle| p_2 , p_7
    \right\rangle 
  \end{multline*}

\item Condition $\delta_C(\boldsymbol{p};u)$
\begin{equation*}
  p_5 + p_6 + p_{11} = p_1 + p_4 + p_8 -  2\, u
\end{equation*}

\item Potential function
  \begin{multline*}
    V_{\mathcal{M}}(v,w,x,y,z; m)
    =
    -\Li\left(\frac{w}{v}\right) - \Li\left(\frac{m^4 \, v}{y} \right)
    +\Li\left(\frac{x}{y}\right)
    - \Li\left(\frac{x}{z^2} \right)
    \\
    +\Li\left(\frac{w}{z}\right) - \Li\left(\frac{z}{m^2 \, v} \right)
    + \frac{\pi^2}{3}    
    - 2 \left( \log \left( m^2 \right) \right)^2
    \\
    +
    \log \left( m^2 \right)
    \log \left( \frac{y^2 \, z}{v^3 \, w} \right)
    + \log v \log \left( \frac{y \, z}{v} \right)
    + \log x \log \left( \frac{w \, z^3}{x} \right)
    -\log z \log \left( y \, z^2 \, w \right)
  \end{multline*}

\item Hyperbolic volume
  \begin{equation*}
    \Im V_{\mathcal{M}}(v,w,x,y,z; 1)
    =
    -D\left(\frac{w}{v}\right)
    - D\left(\frac{ v}{y} \right)
    + D\left(\frac{x}{y}\right)
    - D\left(\frac{x}{z^2} \right)
    + D\left(\frac{w}{z}\right)
    - D\left(\frac{z}{ v} \right)
  \end{equation*}
  with
  \begin{equation*}
    \begin{pmatrix}
      v \\ w\\ x\\ y \\ z
    \end{pmatrix}
    =
    \begin{pmatrix}
      1.18608 + 0.874646 \, \I
      \\
      1.0 
      \\
      -1.09737 + 0.230836 \, \I
      \\
      -1.23271 +   1.09381 \, \I
      \\
      0.40897 - 0.337176 \, \I
    \end{pmatrix}
  \end{equation*}

\item A-polynomial
  \begin{equation*}
    \left(
    \begin{array}{ccccccccccc}
      0 & 0 & 0 & 0 & 0 & -1 & 1 & 0 & 0 & 0 & 0
      \\
      0 & 0 & 0 & 0 & 0 & 10 & -7 & 3 & -1 & 1 & 0
      \\
      0 & 0 & 0 & 0 & 1 & -24 & 16 & -3 & 6 & -2 & 1
      \\
      0 & 0 & 0 & 2 & -4 & -5 & 5 & -19 & 1 & 1 & 0
      \\
      0 & 0 & 0 & -10 & 7 & 47 & -39 & 13 & -12 & -6 & 0
      \\
      0 & 0 & -1 & 15 & 5 & 18 & -14 & 22 & -5 & 6 & 0
      \\
      0 & -1 & 3 & 3 & -36 & -21 & 45 & 12 & 9 & 0 & 0
      \\
      0 & 1 & 7 & -35 & 16 & -69 & 30 & 2 & -6 & 0 & 0
      \\
      0 & 0 & -19 & 15 & 44 & -45 & 6 & -40 & -7 & 0 & 0
      \\
      0 & 0 & 7 & 40 & -6 & 45 & -44 & -15 & 19 & 0 & 0
      \\
      0 & 0 & 6 & -2 & -30 & 69 & -16 & 35 & -7 & -1 & 0
      \\
      0 & 0 & -9 & -12 & -45 & 21 & 36 & -3 & -3 & 1 & 0
      \\
      0 & -6 & 5 & -22 & 14 & -18 & -5 & -15 & 1 & 0 & 0
      \\
      0 & 6 & 12 & -13 & 39 & -47 & -7 & 10 & 0 & 0 & 0
      \\
      0 & -1 & -1 & 19 & -5 & 5 & 4 & -2 & 0 & 0 & 0
      \\
      -1 & 2 & -6 & 3 & -16 & 24 & -1 & 0 & 0 & 0 & 0
      \\
      0 & -1 & 1 & -3 & 7 & -10 & 0 & 0 & 0 & 0 & 0
      \\
      0 & 0 & 0 & 0 & -1 & 1 & 0 & 0 & 0 & 0 & 0
    \end{array}
    \right)
  \end{equation*}
\end{itemize}
\subsubsection{
$K6_{22}$}

\begin{align*}
  & \mbox{
    \raisebox{-1.6cm}{
      \includegraphics[scale=.16]{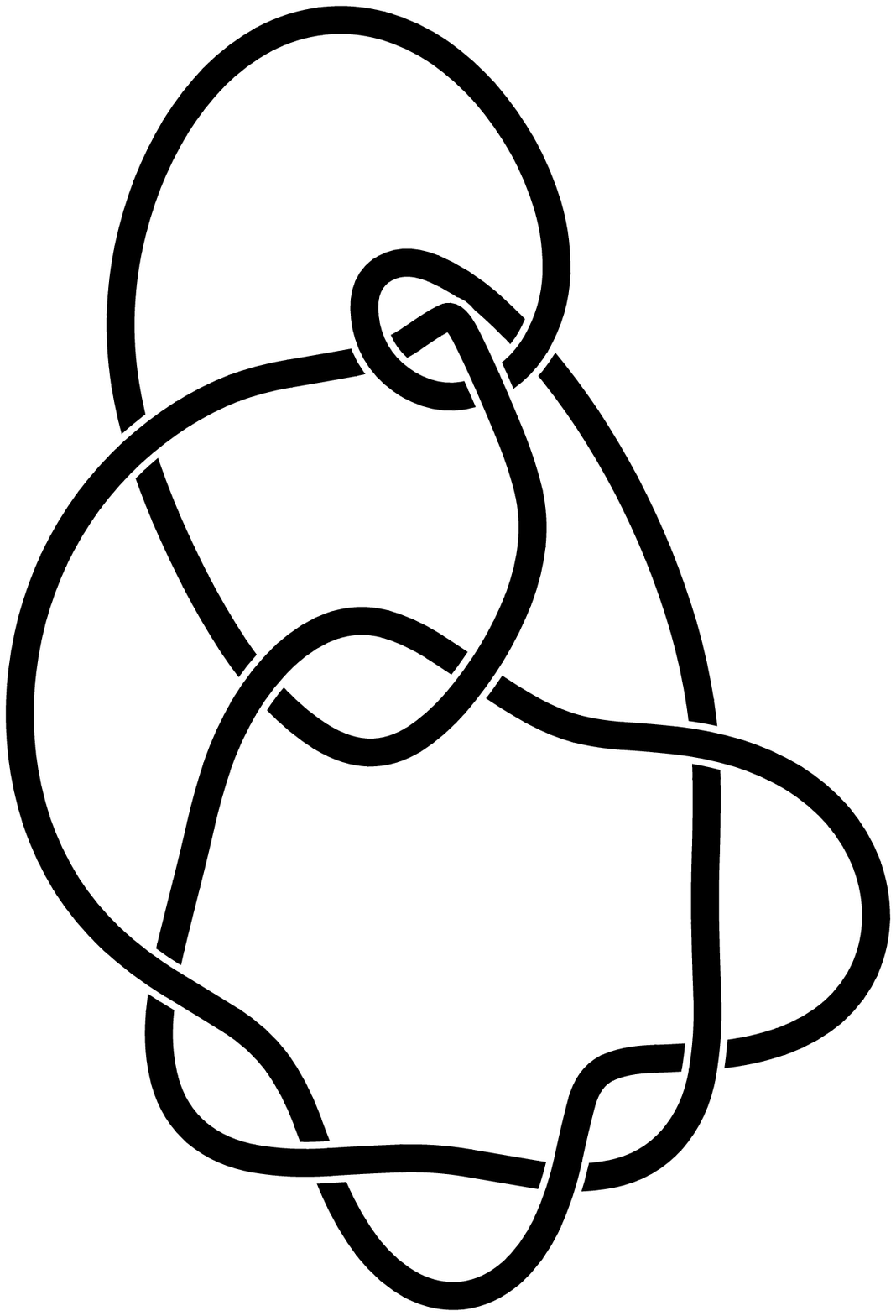}
    }}
  &
  & \Vol(S^3 \setminus \mathcal{K})= 4.7698896 \dots
\end{align*}

\begin{itemize}

\item Quantum invariant
\begin{multline*}
  Z_\gamma \left( \mathcal{M}_u \right)
  =
  \int\limits_\mathbb{R}
  \mathrm{d} \boldsymbol{p} \,
  \delta_C(\boldsymbol{p}; u ) \,
  \left\langle
    p_6 , p_1 \middle| S^{-1} \middle| p_1 , p_2
  \right\rangle \,
  \left\langle
    p_2 , p_{11} \middle| S^{-1} \middle| p_3 , p_7
  \right\rangle \,
  \left\langle
    p_{10} , p_7 \middle| S \middle| p_{12} , p_{10}
  \right\rangle
  \\
  \times
  \left\langle
    p_{3} , p_9 \middle| S^{-1} \middle| p_4 , p_{11}
  \right\rangle \,
  \left\langle
    p_{5} , p_4 \middle| S \middle| p_9 , p_{8}
  \right\rangle \,
  \left\langle
    p_{8} , p_{12} \middle| S \middle| p_6 , p_{5}
  \right\rangle 
\end{multline*}

\item
Condition $\delta_C(\boldsymbol{p};u)$
\begin{equation*}
  p_4 + p_5 + p_7 = p_8 + p_9 + p_{12} - 2 \, u
\end{equation*}

\item Potential function
\begin{multline*}
  V_\mathcal{M}(v,w,x,y,z; m)
  =
  - \Li\left( \frac{w}{v^2} \right)
  - \Li\left( \frac{m^2 \, v}{w^2 \, x^2} \right)
  + \Li\left( \frac{m^2 \, v}{y} \right)
  + \Li\left( \frac{y}{v^2 \, x} \right)
  \\
  - \Li\left( \frac{v \, x \, z}{w} \right)
  + \Li\left( \frac{x \, z}{y} \right)
  - \left( \log \left( m^2 \right) \right)^2
  + \log v \log \left(w \, z^2\right)
  - 3 \left( \log w \right)^2
  \\
  +
  \log\left( m^2 \right) \log \left(\frac{w^3 \, x^4 \, z^2}{y}\right)
  -2 \log y \log \left( v^2 \, x^2 \, y \, z \right)
  + \log x \log \left( \frac{v^8}{w^4 \, x} \right)
\end{multline*}

\item Hyperbolic volume
\begin{equation*}
  \Im V_{\mathcal{M}}(v,w,x,y,z; 1)
  =
  - D\left( \frac{w}{v^2} \right)
  - D\left( \frac{ v}{w^2 \, x^2} \right)
  + D\left( \frac{ v}{y} \right)
  + D\left( \frac{y}{v^2 \, x} \right)
  - D\left( \frac{v \, x \, z}{w} \right)
  + D\left( \frac{x \, z}{y} \right)
\end{equation*}
with
\begin{equation*}
  \begin{pmatrix}
    v \\ w \\ x\\ y \\ z
  \end{pmatrix}
  =
  \begin{pmatrix}
    0.127216 + 0.708358 \, \I
    \\
    -0.534645 + 1.01011 \, \I
    \\
    -0.40328 + 0.755585 \, \I
    \\
    0.66186 - 0.301752 \, \I
    \\
    -0.54977 - 1.03005 \, \I
  \end{pmatrix}
\end{equation*}

\item A-polynomial
  \begin{equation*}
    \left(
      \begin{array}{ccccccccccc}
        0 & 0 & 0 & 0 & 0 & 0 & 0 & 0 & 0 & -1 & 0 \\
        0 & 0 & 0 & 0 & 0 & 0 & 0 & 0 & 0 & 1 & 0 \\
        0 & 0 & 0 & 0 & 0 & 0 & 0 & 0 & 0 & 0 & 0 \\
        0 & 0 & 0 & 0 & 0 & 0 & 0 & 0 & 0 & 3 & 0 \\
        0 & 0 & 0 & 0 & 0 & 0 & 0 & 0 & -1 & -3 & 1 \\
        0 & 0 & 0 & 0 & 0 & 0 & 0 & 0 & 1 & 5 & 0 \\
        0 & 0 & 0 & 0 & 0 & 0 & 0 & 0 & 9 & -8 & 0 \\
        0 & 0 & 0 & 0 & 0 & 0 & 0 & -2 & -16 & 3 & 0 \\
        0 & 0 & 0 & 0 & 0 & 0 & 0 & 3 & 19 & 2 & 0 \\
        0 & 0 & 0 & 0 & 0 & 0 & 0 & 20 & -13 & 0 & 0 \\
        0 & 0 & 0 & 0 & 0 & 0 & -1 & -44 & -5 & 0 & 0 \\
        0 & 0 & 0 & 0 & 0 & 0 & -1 & 20 & -4 & 0 & 0 \\
        0 & 0 & 0 & 0 & 0 & 0 & 30 & 8 & 9 & 0 & 0 \\
        0 & 0 & 0 & 0 & 0 & -1 & -56 & 2 & -2 & 0 & 0 \\
        0 & 0 & 0 & 0 & 0 & -5 & -10 & -36 & 0 & 0 & 0 \\
        0 & 0 & 0 & 0 & 0 & 35 & 30 & 22 & 0 & 0 & 0 \\
        0 & 0 & 0 & 0 & 1 & -32 & 69 & 1 & 0 & 0 & 0 \\
        0 & 0 & 0 & 0 & -9 & -41 & -65 & -2 & 0 & 0 & 0 \\
        0 & 0 & 0 & 0 & 16 & -13 & -9 & 0 & 0 & 0 & 0 \\
        0 & 0 & 0 & 0 & 7 & 126 & 7 & 0 & 0 & 0 & 0 \\
        0 & 0 & 0 & 0 & -9 & -13 & 16 & 0 & 0 & 0 & 0 \\
        0 & 0 & 0 & -2 & -65 & -41 & -9 & 0 & 0 & 0 & 0 \\
        0 & 0 & 0 & 1 & 69 & -32 & 1 & 0 & 0 & 0 & 0 \\
        0 & 0 & 0 & 22 & 30 & 35 & 0 & 0 & 0 & 0 & 0 \\
        0 & 0 & 0 & -36 & -10 & -5 & 0 & 0 & 0 & 0 & 0 \\
        0 & 0 & -2 & 2 & -56 & -1 & 0 & 0 & 0 & 0 & 0 \\
        0 & 0 & 9 & 8 & 30 & 0 & 0 & 0 & 0 & 0 & 0 \\
        0 & 0 & -4 & 20 & -1 & 0 & 0 & 0 & 0 & 0 & 0 \\
        0 & 0 & -5 & -44 & -1 & 0 & 0 & 0 & 0 & 0 & 0 \\
        0 & 0 & -13 & 20 & 0 & 0 & 0 & 0 & 0 & 0 & 0 \\
        0 & 2 & 19 & 3 & 0 & 0 & 0 & 0 & 0 & 0 & 0 \\
        0 & 3 & -16 & -2 & 0 & 0 & 0 & 0 & 0 & 0 & 0 \\
        0 & -8 & 9 & 0 & 0 & 0 & 0 & 0 & 0 & 0 & 0 \\
        0 & 5 & 1 & 0 & 0 & 0 & 0 & 0 & 0 & 0 & 0 \\
        1 & -3 & -1 & 0 & 0 & 0 & 0 & 0 & 0 & 0 & 0 \\
        0 & 3 & 0 & 0 & 0 & 0 & 0 & 0 & 0 & 0 & 0 \\
        0 & 0 & 0 & 0 & 0 & 0 & 0 & 0 & 0 & 0 & 0 \\
        0 & 1 & 0 & 0 & 0 & 0 & 0 & 0 & 0 & 0 & 0 \\
        0 & -1 & 0 & 0 & 0 & 0 & 0 & 0 & 0 & 0 & 0 
      \end{array}
    \right)
  \end{equation*}
\end{itemize}

\subsubsection{
$K6_{33}$ or  $10_{140}$}

\begin{align*}
  & \mbox{
    \raisebox{-1.6cm}{
      \includegraphics[scale=.16]{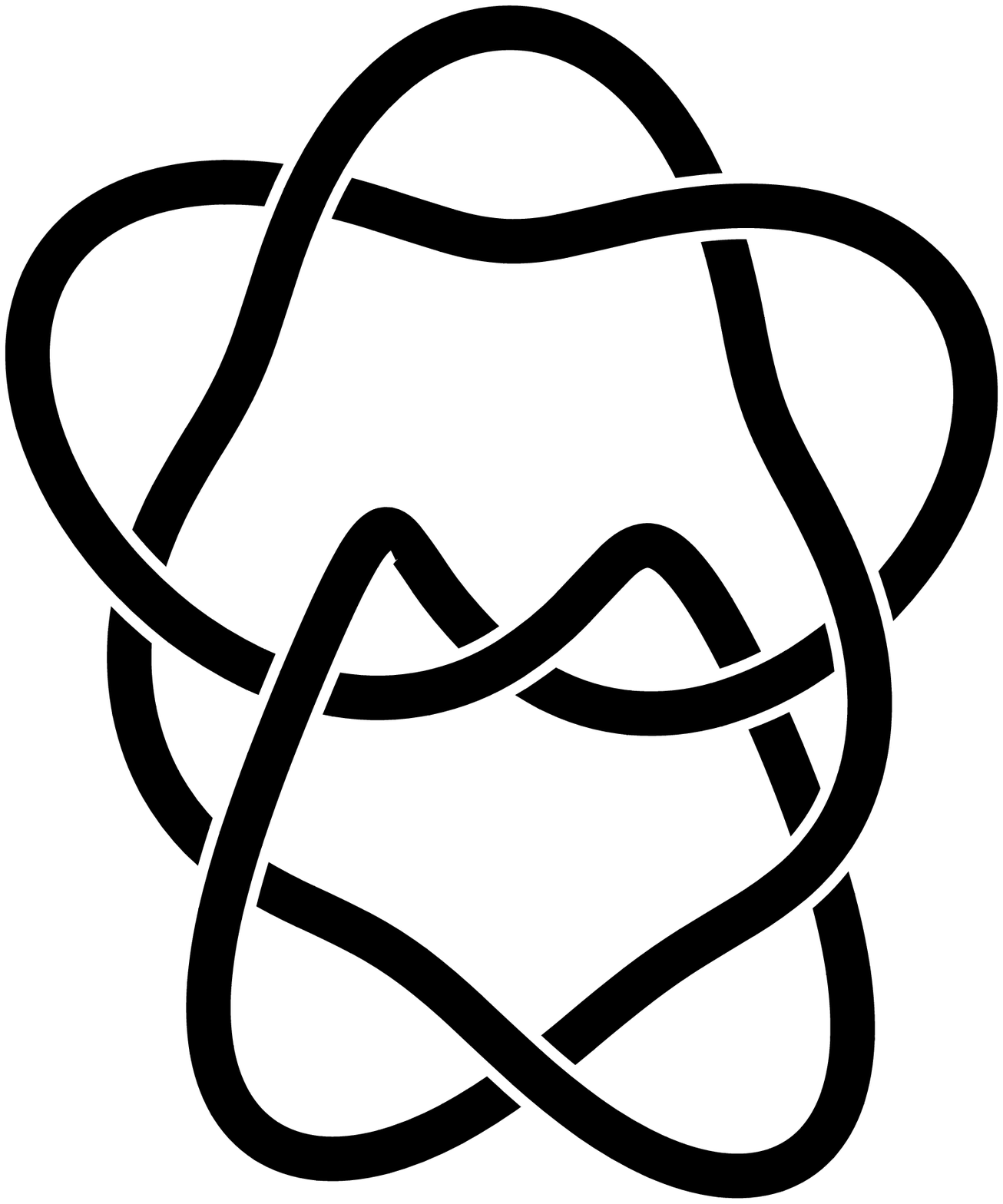}
    }}
  &
  & \Vol(S^3 \setminus \mathcal{K})= 5.212567 \dots
\end{align*}

\begin{itemize}
\item  Quantum invariant
\begin{multline*}
  Z_\gamma \left( \mathcal{M}_u \right)
  = \int\limits_\mathbb{R} 
 \mathrm{d} \boldsymbol{p} \,
  \delta_C(\boldsymbol{p}; u ) \,
  \langle p_1 , p_7 | S | p_5 , p_3 \rangle \,
  \langle p_2 , p_9 | S^{-1} | p_1 , p_4 \rangle 
  \langle p_3 , p_{12} | S^{-1} | p_{10} , p_2 \rangle 
  \\
  \times 
  \langle p_5 , p_4 | S^{-1} | p_6 , p_{11} \rangle \,
  \langle p_8 , p_6 | S^{-1} | p_7 , p_{12} \rangle \,
  \langle p_{10} , p_{11} | S^{-1} | p_9 , p_8 \rangle 
\end{multline*}

\item
Condition $\delta_C(\boldsymbol{p};u)$
\begin{equation*}
  p_3 + p_6 + p_{11} = p_2 + p_7 + p_9 -  2\, u
\end{equation*}

\item Potential function
\begin{multline*}
  V_{\mathcal{M}}(v,w,x,y,z; m)
  =
  - \Li\left(\frac{1}{v \, y} \right)
  - \Li\left(\frac{1}{m^2 \, w \, y} \right)
  + \Li\left(\frac{w^2 \, x}{y} \right)
  \\
  - \Li\left(v \, m^2 \, w \, y\right) 
  - \Li\left( v \, m^2 \, w \, x \, z \right)
  - \Li\left(m^2 \, y \, z\right)
  +  \frac{2 \, \pi^2}{3}
  \\
  - 
  \left( \log(v \, m^2 ) \right)^2
  - 2 \log v \log w
  - 2 \log \left( m^2 \right) \log(w \, y)
  \\
  - \log v \log ( x \, y)
  - \log \left( m^2 \right) \log z
  - \log y \, \log(z \, w \, y)
\end{multline*}

\item Hyperbolic volume
\begin{multline*}
  \Im V_{\mathcal{M}}(v,w,x,y,z; 1)
  =
  - D\left(\frac{1}{v \, y} \right)
  - D\left(\frac{1}{ w \, y} \right)
  + D\left(\frac{w^2 \, x}{y} \right)
  \\
  - D\left(v  \, w \, y\right) 
  - D\left( v  \, w \, x \, z \right)
  - D\left( y \, z\right)
\end{multline*}
with
\begin{equation*}
  \begin{pmatrix}
    v \\ w\\ x\\ y \\ z
  \end{pmatrix}
  =
  \begin{pmatrix}
     -1.1238 - 0.998279 \, \I
     \\
     -0.439261 - 0.570751 \, \I
     \\
     -0.836795 + 1.7323  \, \I
     \\
     -0.829546 - 0.0564355 \, \I
     \\
      -0.549394 + 0.740149 \, \I
  \end{pmatrix}
\end{equation*}

\item A-polynomial
\begin{equation*}
  \begin{pmatrix}
    0 & -1 & 1 & 0 & 0 & 0 & 0 & 0 \\
    0 & 1 & -9 & 0 & 0 & 0 & 0 & 0 \\
    0 & -2 & 8 & -12 & -2 & 0 & 0 & 0 \\
    0 & 5 & -4 & 6 & -12 & -3 & 0 & 0 \\
    0 & -6 & 5 & -7 & 4 & -13 & -6 & -1 \\
    0 & 4 & -7 & 8 & -20 & 1 & -2 & 0 \\
    0 & -2 & 1 & -20 & 8 & -7 & 4 & 0 \\
    -1 & -6 & -13 & 4 & -7 & 5 & -6 & 0 \\
    0 & 0 & -3 & -12 & 6 & -4 & 5 & 0 \\
    0 & 0 & 0 & -2 & -12 & 8 & -2 & 0 \\
    0 & 0 & 0 & 0 & 0 & -9 & 1 & 0 \\
    0 & 0 & 0 & 0 & 0 & 1 & -1 & 0 
  \end{pmatrix}
\end{equation*}
\end{itemize}

\bibliographystyle{alphaKH}

\end{document}